\documentclass[12pt, reqno, a4paper]{amsart}
\usepackage{amsthm, amsmath, amssymb, amscd,latexsym}
\usepackage{enumerate}
\usepackage{latexsym}
\usepackage{graphicx, color}
\def\C{\mathbb C}

\def\D{\mathbb D}

\def\N{\mathbb N}

\def\disp{\displaystyle}
\def\wt{\widetilde}

\def\noin{\noindent}
\def\dist{{\rm dist}}

\def\vep{\varepsilon}

\def\pt{\partial}

\theoremstyle{plain}
\newtheorem{thm}{Theorem}[section]
\newtheorem{cor}[thm]{Corollary}
\newtheorem{prop}[thm]{Proposition}
\newtheorem{lem}[thm]{Lemma}

\newtheorem*{thmA*}{Theorem A}
\newtheorem*{thmA'*}{Theorem A'}
\newtheorem*{corA'*}{Corollary A'}
\newtheorem*{thmB*}{Theorem B}
\newtheorem*{thmC*}{Theorem C}
\newtheorem*{corC*}{Corollary C}
\newtheorem*{corD*}{Corollary D}
\newtheorem*{thmD*}{Theorem D}
\newtheorem*{thmE*}{Theorem E}
\newtheorem*{thmF*}{Theorem F}
\newtheorem*{conj*}{Conjecture}
\newtheorem*{thm*}{Theorem}

\theoremstyle{definition}
\newtheorem{defn}{Definition}[section]

\newtheorem{Question*}[thm]{Question}
\newtheorem*{defn*}{Definition}
\newtheorem*{rem*}{Remark}

\makeatletter
\@addtoreset{equation}{section} 
\makeatother

\newcommand{\bs}[1]{\boldsymbol{#1}}
\renewcommand{\Bar}{\overline}
\newcommand{\norm}[1]{{\left\| #1 \right\|}}
\newcommand{\abs}[1]{{\left| #1 \right|}}

\newcommand{\paren}[1]{{\left( #1 \right)}}
\newcommand{\skakko}[1]{{\left\{ #1 \right\}}}
\newcommand{\braces}[1]{{\left\{ #1 \right\}}}
\newcommand{\al}{\alpha}

\newcommand{\cM}{{\mathcal{M}}}

\newcommand{\cK}{{\mathcal{K}}}
\newcommand{\cc}{\circ}
\newcommand{\st}{\ | \ }
\newcommand{\QED}{\hfill $\blacksquare$}

\newcommand{\sminus}{\smallsetminus}
\newcommand{\sig}{\sigma}

\begin{document}

\parindent=5mm
\baselineskip=15pt
\parskip=3pt

\title{Julia sets appear quasiconformally
 in the Mandelbrot set}

\author{Tomoki Kawahira and Masashi Kisaka}
\setcounter{footnote}{-1}

\address{
Graduate School of Economics\\
Hitotsubashi University\\
Tokyo 186-8601 \\
Japan}
\email{t.kawahira@r.hit-u.ac.jp}

\address{
Department of Mathematical Sciences \\
Graduate School of Human and Environmental Studies \\
Kyoto University \\
Kyoto 606-8501 \\
Japan }
\email{kisaka@math.h.kyoto-u.ac.jp}

\begin{abstract}
In this paper we prove the following: 
In the boundary of the Mandelbrot set,
we can find quasiconformal copies of a Cantor Julia set 
which is a small perturbation of the Julia set 
of any given parabolic or Misiurewicz parameter.
Indeed, we can specify the locations of 
such copies near the boundary of any small Mandelbrot set.
If we zoom in the middle part of such a copy,
then we can find a certain nested structure (\lq\lq decoration") and finally another 
\lq\lq smaller Mandelbrot set" appears. A similar nested 
structure exists in the Julia set for any parameter in the \lq\lq smaller 
Mandelbrot set". We can also find images of Julia sets by quasiconformal
maps with dilatation arbitrarily close to 1. This answers a question by
Adrian Douady. All the parameters belonging to 
these images are semihyperbolic and this leads to the fact that the set of
semihyperbolic but non-Misiurewicz and non-hyperbolic parameters is dense 
with Hausdorff dimension 2 in the boundary of the Mandelbrot set.
\end{abstract}

\maketitle

\section{Introduction}

Let $P_c(z) := z^2 +c \ (c \in \C)$ and recall that its 
{\it filled Julia set} $K(P_c)$ is defined by
$$
  K(P_c) := \{ z \in \C \ | \ \{ P_c^n(z) \}_{n=0}^\infty \ \text{is bounded} \}
$$
and its {\it Julia set} $J(P_c)$ is the boundary of $K(P_c)$, that is, 
$J(P_c) := \partial K(P_c)$. It is known that $J(P_c)$ is connected if and
only if the critical orbit $\{ P_c^n(0) \}_{n=0}^\infty$ is bounded
and if $J(P_c)$ is disconnected, then it is a Cantor set. The connectedness locus
of the quadratic family $\{ P_c \}_{c \in \mathbb C}$ is the famous 
{\it Mandelbrot set} and we denote it by $M$:
$$
M := \{ c \in \C \ | \ J(P_c) \  \text{is connected} \} 
= \{ c \in \C \ | \ \{ P_c^n(0) \}_{n=0}^\infty  \  \text{is bounded} \}.
$$
A parameter $c$ is called a {\it Misiurewicz parameter} if the 
critical point $0$ is strictly preperiodic, that is,  
$$
  P_c^k(P_c^l(0)) = P_c^l(0) \quad \text{and} \quad 
  P_c^k(P_c^{l-1}(0)) \ne P_c^{l-1}(0)
$$
for some $k, \ l \in \N = \{ 1,2,3,\cdots \}$. A parameter $c$ is called a 
{\it parabolic parameter}
if $P_c$ has a parabolic periodic point. Here, a periodic point $z_0$ with 
period $m$ is called {\it parabolic} if $P_c^m(z_0) = z_0$ and its multiplier
$(P_c^m)'(z_0)$ is a root of unity. For the basic knowledge of complex dynamics,
we refer to \cite{Beardon 1991} and \cite{Milnor 2006}.

Douady et al. (\cite{Douady 2000}) proved the following: At a small neighborhood
of the cusp point $c_0 \ne 1/4$ in $M$, which is in a \lq\lq primitive
small Mandelbrot set", there is a sequence $\{ M_n \}_{n \in \N}$ of
small quasiconformal copies of $M$ tending to $c_0$. Moreover each $M_n$ is
encaged in a nested sequence of sets which are homeomorphic to the preimage
of $J(P_{1/4 + \eta})$ (for $\eta > 0$ small) by $z \mapsto z^{2^m}$ for 
$m \geq 0$ and accumulate on $M_n$. 

\begin{figure}[htbp]
\hskip -40mm
{\small (1)}
\hskip 43mm
{\small (2)}
\hskip 43mm
{\small (3)}

\includegraphics[scale=0.19, bb = 0 0 640 480]{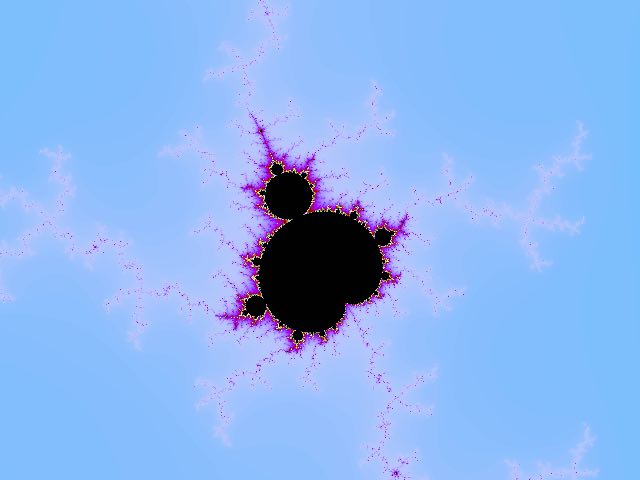} \hskip 5mm
\includegraphics[scale=0.19, bb = 0 0 640 480]{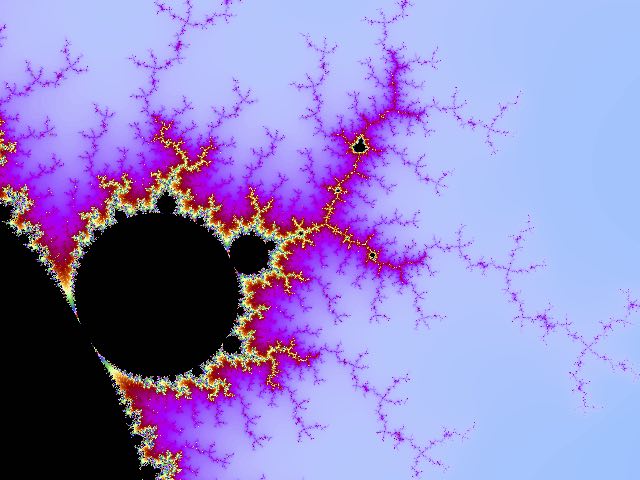} \hskip 5mm
\includegraphics[scale=0.19, bb = 0 0 640 480]{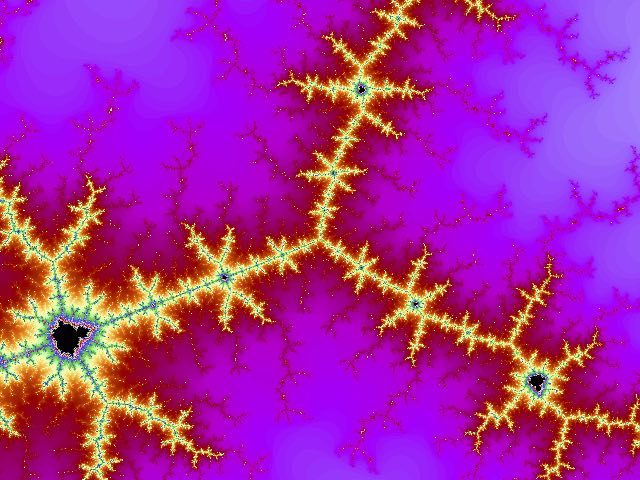} \hskip 5mm

\hskip -40mm
{\small (4)}
\hskip 43mm
{\small (5)}
\hskip 43mm
{\small (6)}

\includegraphics[scale=0.19, bb = 0 0 640 480]{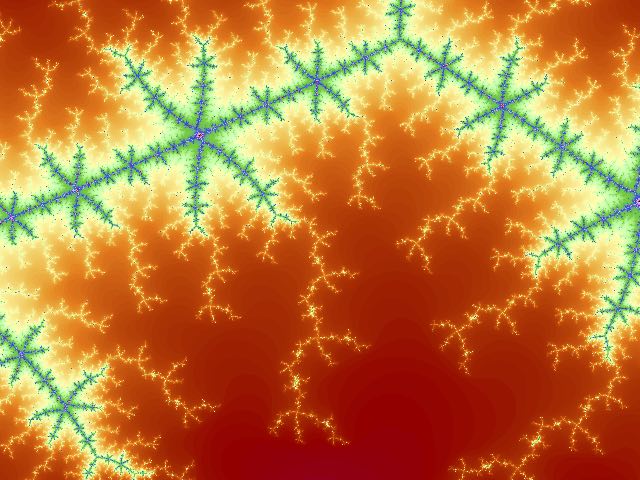} \hskip 5mm
\includegraphics[scale=0.19, bb = 0 0 640 480]{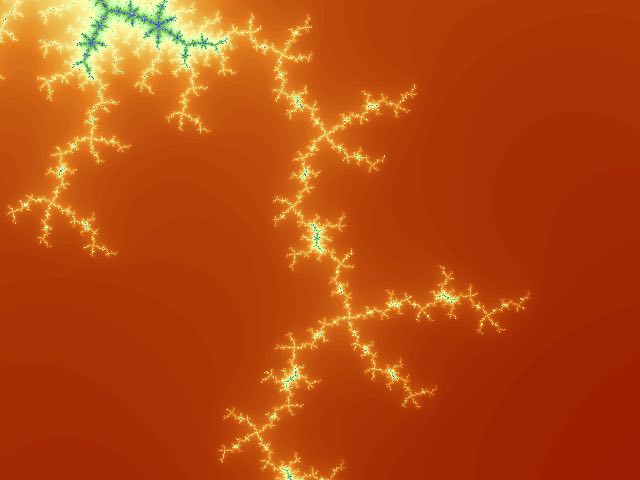} \hskip 5mm
\includegraphics[scale=0.19, bb = 0 0 640 480]{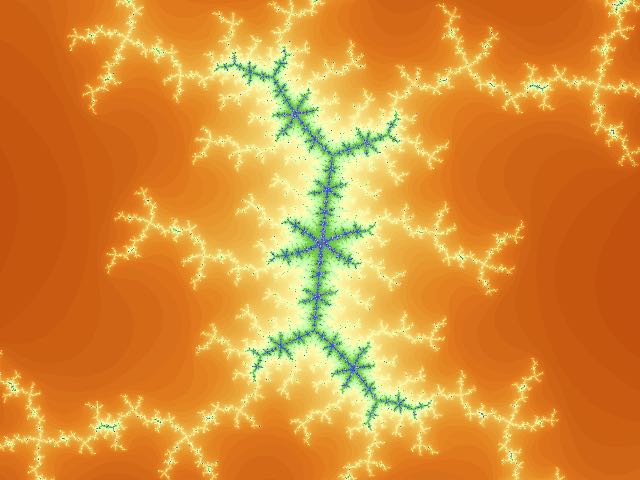} \hskip 5mm

\hskip -40mm
{\small (7)}
\hskip 43mm
{\small (8)}
\hskip 43mm
{\small (9)}

\includegraphics[scale=0.19, bb = 0 0 640 480]{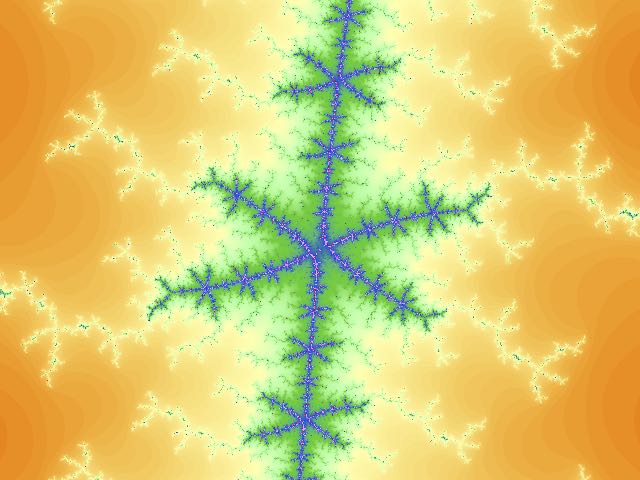} \hskip 5mm
\includegraphics[scale=0.19, bb = 0 0 640 480]{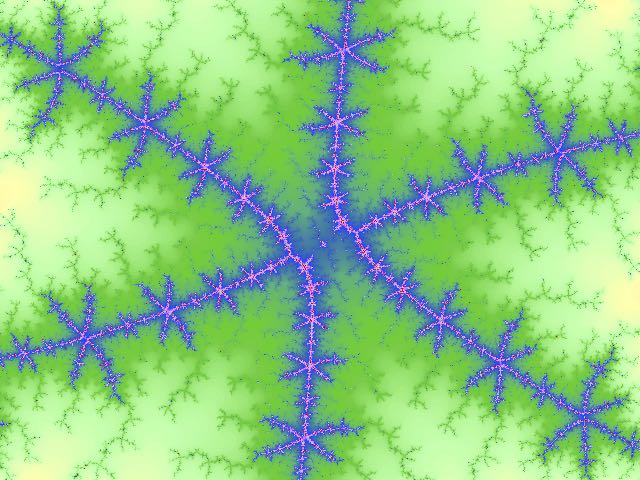} \hskip 5mm
\includegraphics[scale=0.19, bb = 0 0 640 480]{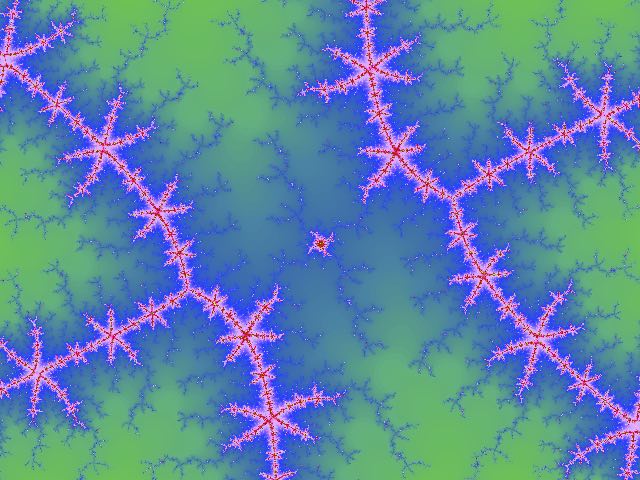} \hskip 5mm

\hskip -38mm
{\small (10)}
\hskip 41mm
{\small (11)}
\hskip 41mm
{\small (12)}

\includegraphics[scale=0.19, bb = 0 0 640 480]{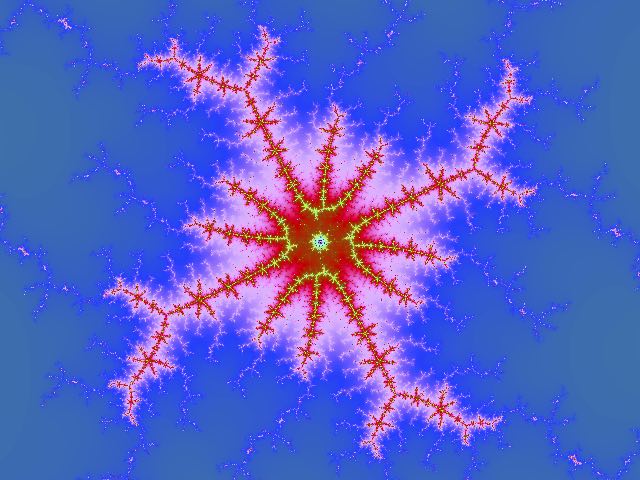} \hskip 5mm
\includegraphics[scale=0.19, bb = 0 0 640 480]{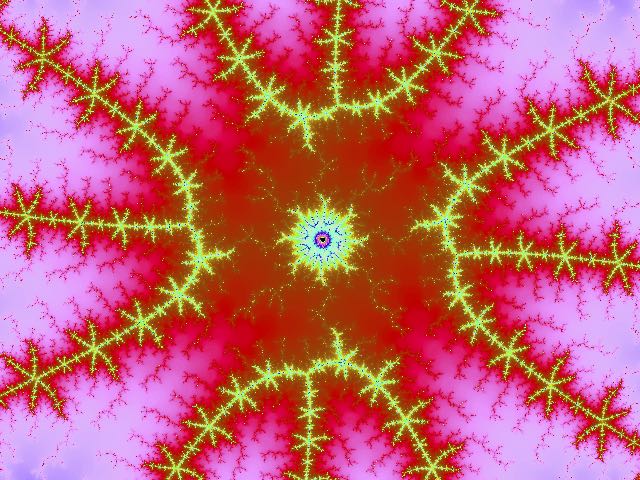} \hskip 5mm
\includegraphics[scale=0.19, bb = 0 0 640 480]{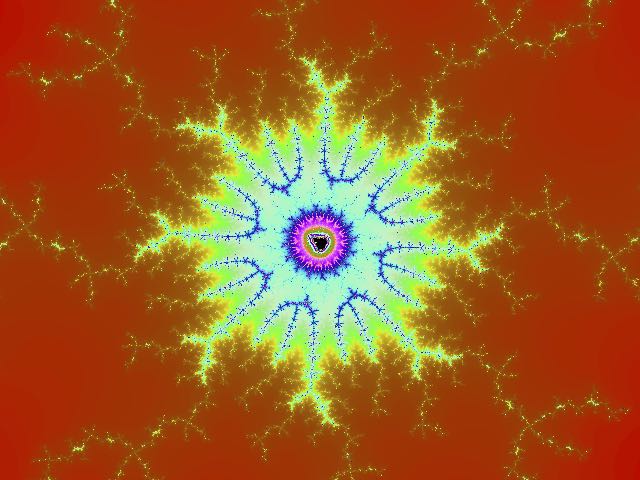} \hskip 5mm

\hskip -38mm
{\small (13)}
\hskip 41mm
{\small (14)}
\hskip 41mm
{\small (15)}

\includegraphics[scale=0.19, bb = 0 0 640 480]{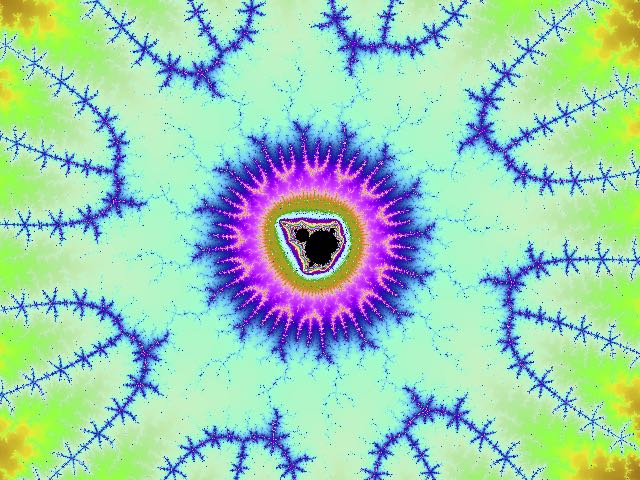} \hskip 5mm
\includegraphics[scale=0.19, bb = 0 0 640 480]{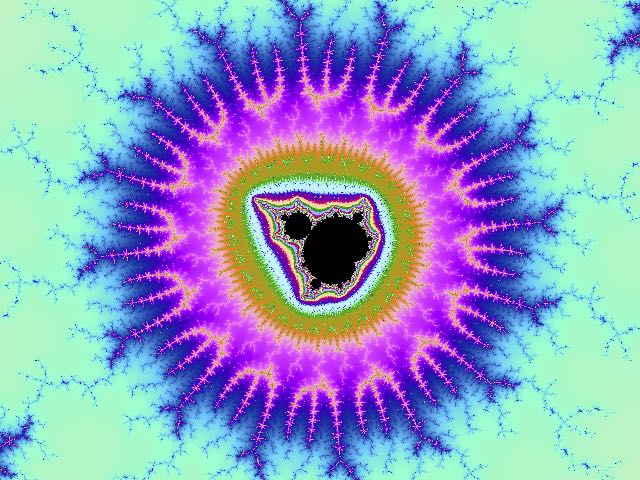} \hskip 5mm
\includegraphics[scale=0.19, bb = 0 0 640 480]{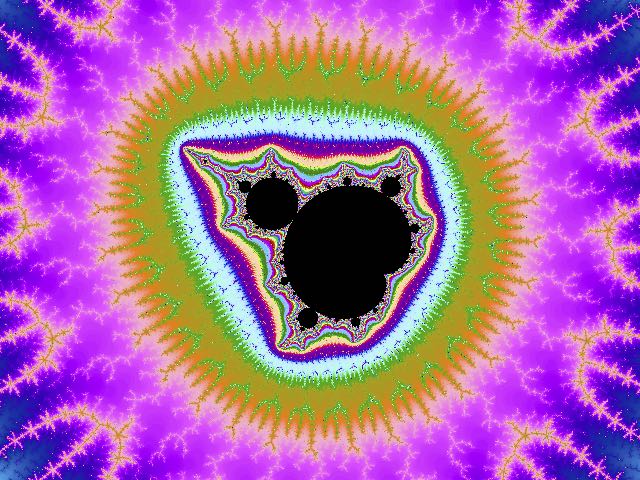} \hskip 5mm

\caption{\small Zooms around a Misiurewicz point 
$c_1 = s_0 \perp c_0$ 
in a primitive small Mandelbrot set $M_{s_0}$, where $c_0$ is a 
Misiurewicz parameter satisfying $P_{c_0}(P_{c_0}^4(0)) = P_{c_0}^4(0)$.
After a sequence of nested structures, another smaller Mandelbrot
set $M_{s_1}$ appears in (15). Here, 
$s_0 \approx 0.3591071125276155 + 0.6423830938166145i$, \
$c_0 \approx -0.1010963638456221 + 0.9562865108091415i$, \
$c_1 \approx 0.3626697754647427 + 0.6450273437137847i$ and
$s_1 \approx 0.3626684938191616 + 0.6450238859863952i$. 
The widths of the figures (1) and (15) are about $10^{-1.5}$
and $10^{-11.9}$, respectively.}
\label{figures of a primitive-Misiurewicz case}
\end{figure}

In this paper, firstly we generalize part of their results (Theorem A).
Actually this kind of phenomena can be observed not only in a small 
neighborhood of the cusp of a \lq\lq primitive small Mandelbrot set", that 
is, the point corresponding to a parabolic parameter $1/4 \in \partial M$, 
but also in every neighborhood of a point corresponding to any Misiurewicz 
or parabolic parameters $c_0$ in a small Mandelbrot set. (For example, 
$c_0=1/4 \in \partial M$ can be replaced by a Misiurewicz parameter 
$c_0=i\in \partial M$ or a parabolic parameter $c_0=-3/4 \in \partial M$ etc.) 
More precisely, we show the following: Take any small Mandelbrot set 
$M_{s_0}$ (Figure 1-(1))) and zoom in the neighborhood of 
$c_1 = s_0 \perp c_0 \in \partial M_{s_0}$ corresponding to $c_0 \in \pt M$ 
which is a Misiurewicz or a parabolic parameter (Figure 1-(2) to (6))). 
(Note that $c_1$ itself is also a Misiurewicz or a parabolic parameter.) 
Then we can find a subset $J' \subset \partial M$ which looks very similar to
$J(P_{c_0})$ (Figure \ref{figures of a primitive-Misiurewicz case}--(6)). 
Zoom in further, then this $J'$ turns out to be similar to $J(P_{c_0+\eta})$
rather than $J(P_{c_0})$, 
where $|\eta|$ is very small and $c_0+\eta \notin M$, 
because $J'$ looks disconnected 
(Figure \ref{figures of a primitive-Misiurewicz case}--(8), (9)). 
Furthermore, as we further zoom in the middle part of $J'$, we can see a 
nested structure which is very similar to the iterated preimages of 
$J(P_{c_0+\eta})$ by $z \mapsto z^2$ (we call these a {\it decoration})
(Figure \ref{figures of a primitive-Misiurewicz case}--(10), (12), (14))
and finally another smaller Mandelbrot set $M_{s_1}$ appears
(Figure \ref{figures of a primitive-Misiurewicz case}--(15)). 
Indeed, we can replace $c_0$ above by {\it any} boundary point of $M$
since both the set of Misiurewicz parameters and 
the set of parabolic parameters form dense subsets of $\partial M$ (Theorem A').

\begin{figure}[htbp] \small
\hskip -40mm
(1)
\hskip 43mm
(2)
\hskip 43mm
(3)

\includegraphics[scale=0.19, bb = 0 0 640 480]{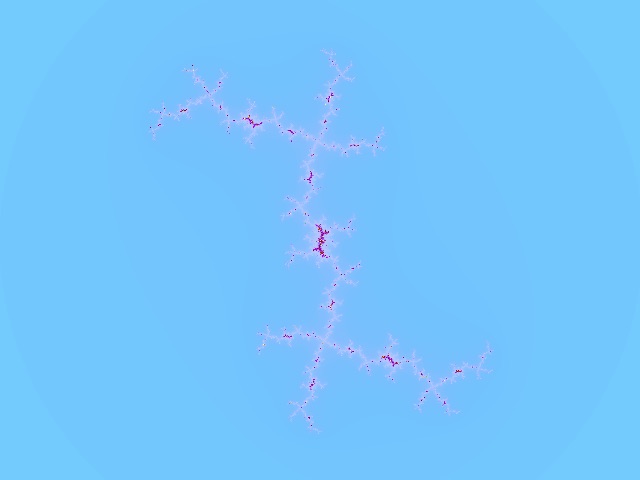} \hskip 5mm
\includegraphics[scale=0.19, bb = 0 0 640 480]{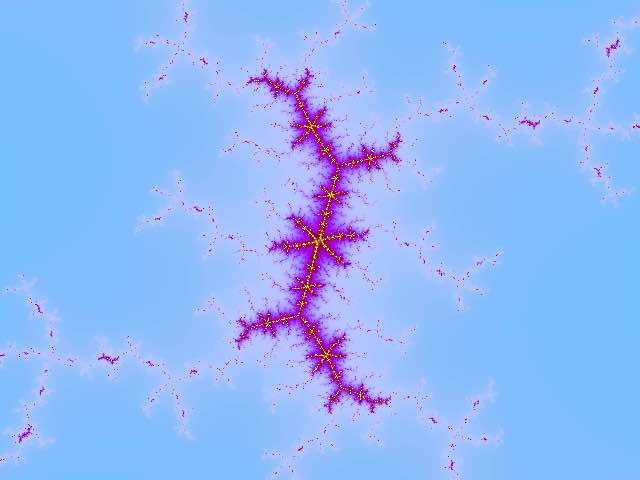} \hskip 5mm
\includegraphics[scale=0.19, bb = 0 0 640 480]{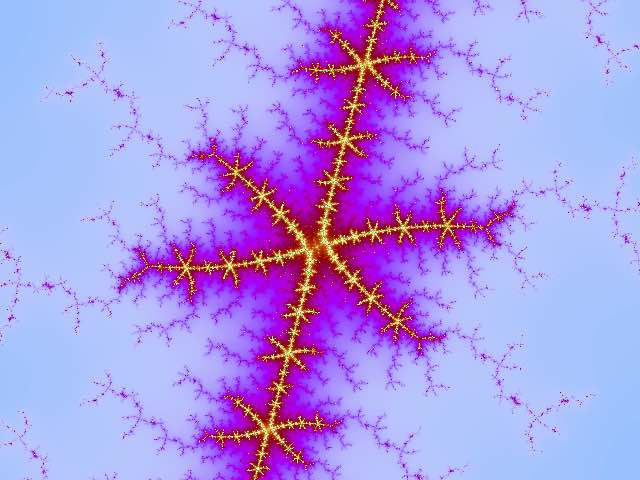}

\hskip -40mm
(4)
\hskip 43mm
(5)
\hskip 43mm
(6)

\includegraphics[scale=0.19, bb = 0 0 640 480]{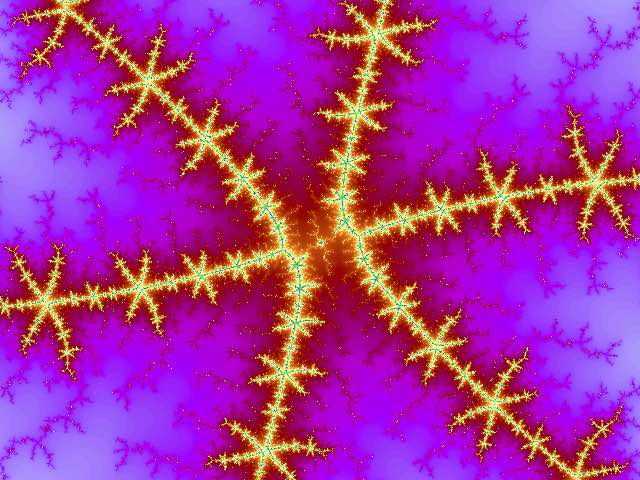} \hskip 5mm
\includegraphics[scale=0.19, bb = 0 0 640 480]{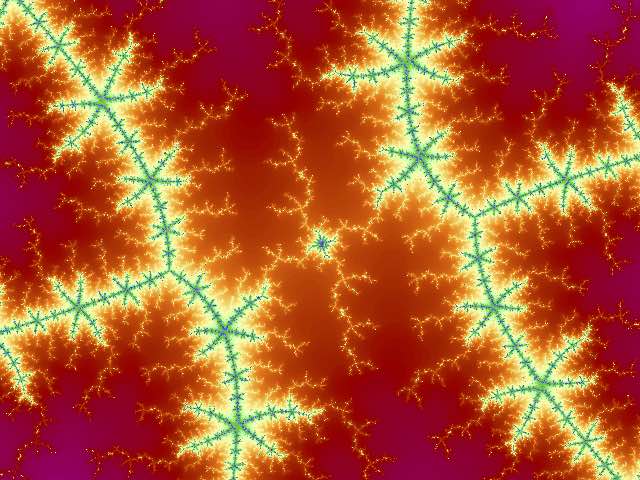} \hskip 5mm
\includegraphics[scale=0.19, bb = 0 0 640 480]{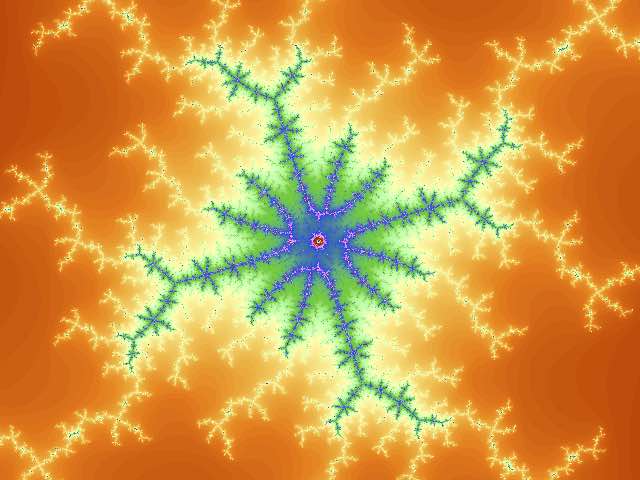}

\hskip -40mm
(7)
\hskip 43mm
(8)
\hskip 43mm
(9)

\includegraphics[scale=0.19, bb = 0 0 640 480]{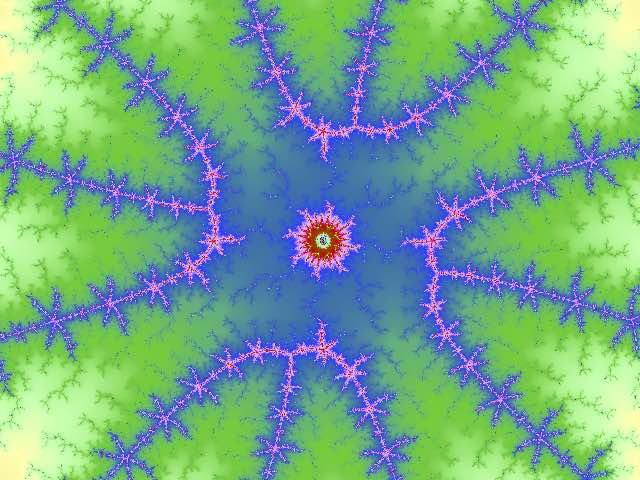} \hskip 5mm
\includegraphics[scale=0.19, bb = 0 0 640 480]{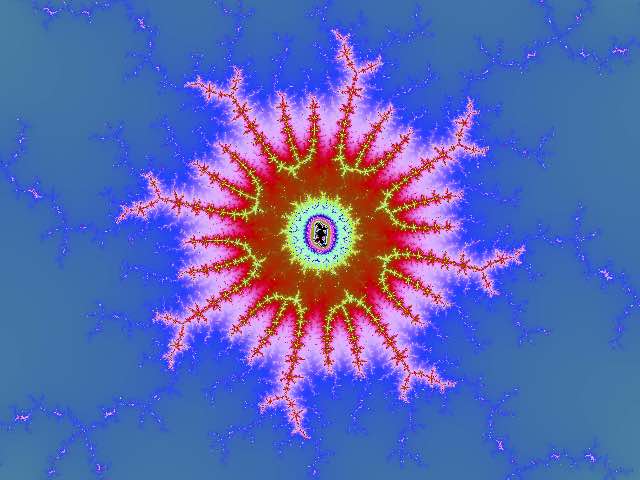} \hskip 5mm
\includegraphics[scale=0.19, bb = 0 0 640 480]{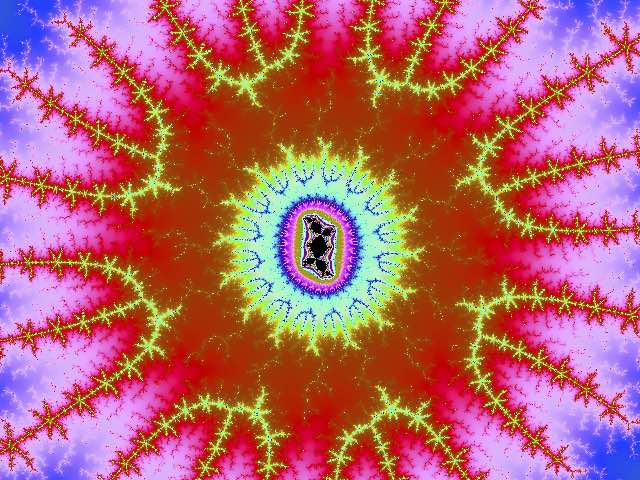}

\hskip -38mm
(10)
\hskip 41mm
(11)
\hskip 41mm
(12)

\includegraphics[scale=0.19, bb = 0 0 640 480]{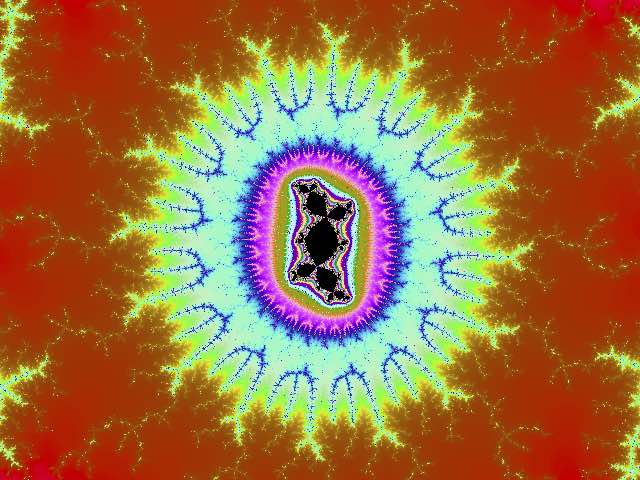} \hskip 5mm
\includegraphics[scale=0.19, bb = 0 0 640 480]{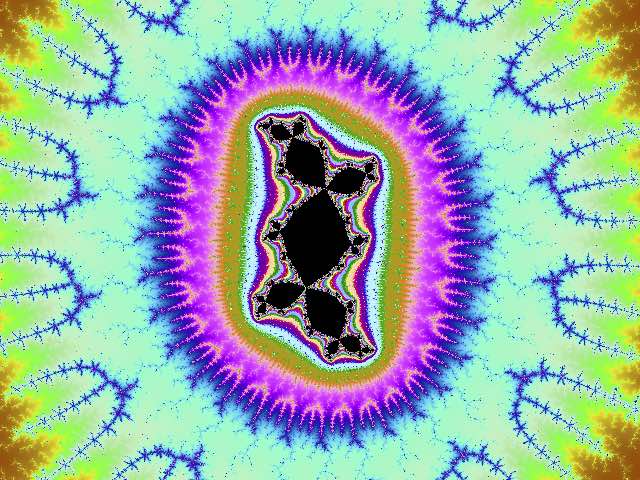} \hskip 5mm
\includegraphics[scale=0.19, bb = 0 0 640 480]{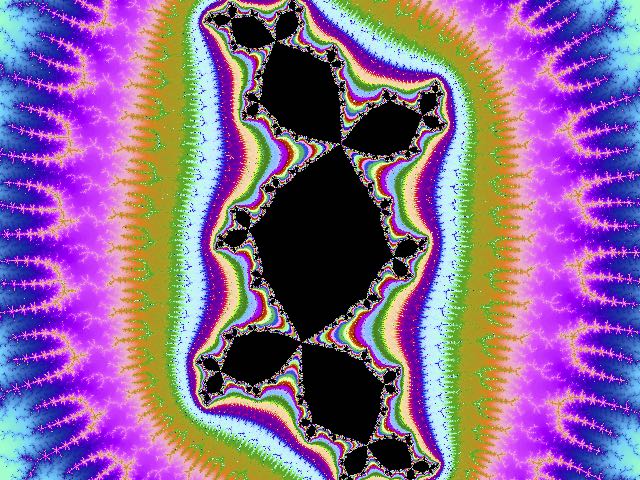} 
\caption{\small Zooms around the critical point $0$ in $K(P_{s_1 \perp c})$
for $s_1 \perp c$ in $M_{s_1}$, 
which is the smaller Mandelbrot set in 
Figure \ref{figures of a primitive-Misiurewicz case}--(15) and $c \in M$
is the parameter for the Douady rabbit. 
$s_1 \approx 0.3626684938191616+0.6450238859863952i$, 
$c \approx -0.12256+0.74486i$ and 
$s_1 \perp c \approx 0.3626684938192285 + 0.6450238859865394i$.
}
\label{nested structure for a filled Julia set}
\end{figure}

Secondly we show the following result for filled Julia sets (Theorem B): 
Take a parameter $s_1 \perp c \ (c \in M)$ from the above smaller Mandelbrot
set $M_{s_1}$ and look at the filled Julia set 
$K(P_{s_1 \perp c})$ and its 
zooms around the neighborhood of $0 \in K(P_{s_1 \perp c})$. Then we can observe
a very similar nested structure to what we saw as zooming in the middle
part of the set $J' \subset \partial M$ 
(see Figure \ref{nested structure for a filled Julia set}).

Thirdly we show that some of the smaller Mandelbrot sets $M_{s_1}$ and their 
decorations are images of  certain model sets by quasiconformal maps whose 
dilatations are arbitrarily close to $1$ (Theorem C). This answers the first
part of the \lq\lq Final remarks" in \cite[p.35]{Douady 2000}.

Finally we show that all 
the parameters belonging to the decorations are semihyperbolic and also the 
set of semihyperbolic but non-Misiurewicz and non-hyperbolic parameters are 
dense in the boundary
of the Mandelbrot set (Corollary D). This together with Theorem C leads to a 
direct and intuitive explanation for the fact that the Hausdorff dimension 
of $\partial M$ is equal to 2, which is a famous result by 
Shishikura (\cite{Shishikura 1998}).

According to Wolf Jung, a structure in the Mandelbrot set which resembles a
whole Julia set in appearance was observed in
computer experiments decades ago by Robert Munafo and Jonathan Leavitt. 
He also claims that he described a general explanation in his
website (\cite{Jung 2015}). We believe some other people  have already 
observed these phenomena so far. For example, we note that Morosawa, 
Nishimura, Taniguchi and Ueda observed this kind of \lq\lq similarity" 
in their book in 1995 (\cite[p.19]{MTU 1995}, \cite[p.26]{MNTU 2000}). 
Further, earlier than this observation, Peitgen observed a kind of local 
similarity between Mandelbrot set and a Julia set by computer experiment 
in 1988 (\cite[Figure 4.23]{Peitgen-Saupe 1988}).


There are different kinds of known results so far which show that
some parts of the Mandelbrot set are similar to some (part of) Julia sets. 
The first famous result for this kind of phenomena is the one by Tan Lei 
(\cite{Tan Lei 1990}). She showed that as we zoom in the neighborhood of
any Misiurewicz parameter $c \in \pt M$, it looks like very much the same as
the magnification of $J(P_c)$ in the neighborhood of $c \in J(P_c)$. Later this
result was generalized to the case where $c$ is a semihyperbolic parameter
by Rivera-Letelier (\cite{Rivera-Letelier 2001}) and its alternative proof
is given by the first author (\cite{Kawahira 2014}). 
On the other hand, some connected Julia sets of quadratic polynomial can
appear quasiconformally in a certain parameter space of a family of cubic 
polynomials. Buff and Henriksen showed that the bifurcation locus of the family 
$\{ f_b(z) = \lambda z + bz^2 + z^3 \}_{b \in \C}$, where  $\lambda \in \C$ with
$|\lambda|=1$ contains quasiconformal copies of $J(\lambda z + z^2)$ 
(\cite{Buff-Henriksen 2001}. See also \cite{Cornell-Rojas-Yampolsky 2017} 
for (non-)computability of the bifurcation locus of such a family for some 
$\lambda$.).

The organization of this paper is as follows: 
In section 2, we construct models for the nested structures mentioned 
above, define the small Mandelbrot set, and show the precise statements 
of the main results (Theorems A, A', B, C and Corollary D). In section 3 we
recall the definitions and basic facts on quadratic-like maps
and Mandelbrot-like families. We prove Theorem A for the Misiurewicz case in 
section 4 and for the parabolic case in section 5. We prove Theorem B in 
section 6.
In section 7 we establish a general formulation
of quadratic-like families that generate \lq\lq fine" copies of the Mandelbrot 
set and we prove Theorem C based on this formulation in section 8. We prove 
Corollary D in section 9
and finally we end this paper with some concluding remarks in section 10. 

\noin
{\bf Acknowledgment: } 
We thank Arnaud Ch\'eritat for informing us a work by Wolf Jung, 
and Wolf Jung for the information of the web pages (\cite{Jung 2015}). 
We also thank the referee for helpful comments. 
The authors were partly supported by JSPS KAKENHI Grants 
16K05193, 17K05296, and 19K03535.

\section{The Model Sets and the Statements of the Results}

\noin
{\bf Notation.} 
We use the following notation for disks and annuli:
\begin{eqnarray*}
& &  D(R) := \{ z \in \C \ | \ |z| < R \}, \quad
  D(\alpha, R) := \{ z \in \C \ | \ |z-\alpha| < R \}, \\
& & A(r, R) := \{ z \in \C \ | \ r < |z| < R \} \quad (0 < r < R).
\end{eqnarray*}
We mostly follow Douady's notations in \cite{Douady 2000} in the 
following. 

\noin
{\bf Models.}
Let $c' \notin M$. Then $J(P_{c'})$ is a Cantor set which does not contain $0$. 
Now take two positive numbers $\rho'$ and $\rho$ such that 
$$
J(P_{c'}) 
\subset 
A(\rho', \rho) \quad (\rho' < \rho). 
$$
We define 
the {\it rescaled Julia set}
$\Gamma_0(c')=\Gamma_0(c')_{\rho',\rho}$
by
$$
\Gamma_0(c') 
:= 
J(P_{c'}) \times \frac{\rho}{(\rho')^2}
=
\braces{ \frac{\rho}{(\rho')^2} \,z \ \big| \ z \in J(P_{c'}) }
$$ 
such that 
$\Gamma_0(c')$ 
is contained in the annulus 
$A(R, R^2)$ with $R:=\rho/\rho'$. 
(In \cite{Douady 2000}, Douady used the radii of the form 
$\rho'=R^{-1/2}$ and $\rho=R^{1/2}$ for some $R>1$ 
such that $\Gamma(c')=J(P_{c'}) \times R^{3/2}$
is contained in $A(R, R^2)$. 
In this paper, however, we need more flexibility
when we are concerned with the dilatation.)

Let $\Gamma_m(c') \ (m \in \N)$ be the inverse image of $\Gamma_0(c')$ by 
$z \mapsto z^{2^m}$. Then $\Gamma_m(c') \ (m=0, \ 1, \ 2, \ \cdots)$ 
are mutually disjoint, because we have
$$
\Gamma_0(c') \subset A(R, R^2), \
\Gamma_1(c') \subset A(R^{1/2}, R), \
\Gamma_2(c') \subset A(R^{1/4}, R^{1/2}), \cdots.
$$

\begin{figure}[htbp]
\includegraphics[width=0.95\textwidth]{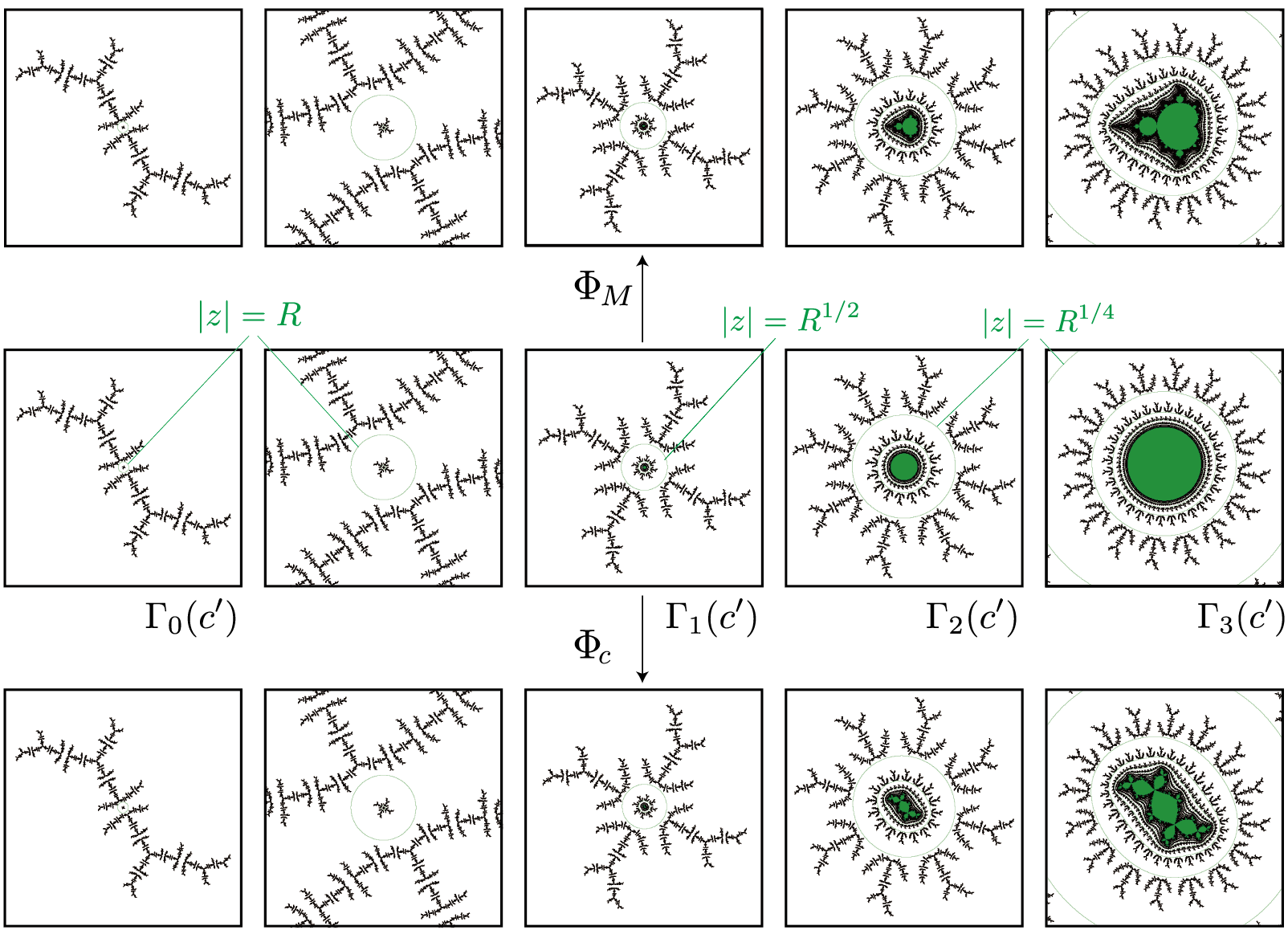}
\caption{\small 
The first row depicts the decorated Mandelbrot set
$\mathcal{M}(c')$ for $c'= -0.10 + 0.97i$
(close to the Misiurewicz parameter $c_0 \approx -0.1011+0.9563i$,
the landing point of the external ray of angle $11/56$)
and $R=220$.
The second row depicts the set $\bigcup_{m \ge 0} \Gamma_m(c')$.
The third row depicts the decorated filled Julia set
$\mathcal{K}_c(c')$ for $c \approx -0.123+0.745$ (the rabbit).
}
\label{figure of models_MK}
\end{figure}

For another parameter $c \in M$, let 
$\Phi_c : \C \smallsetminus K(P_c) \to \C \smallsetminus \overline{\D}$ be the
B\"ottcher coordinate (i.e., $\Phi_c$ is a conformal isomorphism with 
$\Phi_c(P_c(z)) = (\Phi_c(z))^2$). 
Let $\Phi_M : \C \smallsetminus M \to \C \smallsetminus \overline{\D}$ 
be the conformal isomorphism with 
$\Phi_M(c)/c \to 1 \ \text{as} \ |c| \to \infty$. (It is known that
$\Phi_M(c) := \Phi_c(c)$. See \cite{DH Orsay}.) Now define the {\it model sets} 
$\mathcal{M}(c')$ and $\mathcal{K}_c(c')$
as follows (see Figure \ref{figure of models_MK}):
$$
  \mathcal{M}(c') 
:= M \cup \Phi_M^{-1}\Big( \bigcup_{m=0}^\infty \Gamma_m(c') \Big), \quad
  \mathcal{K}_c(c') 
:= K(P_c) \cup \Phi_c^{-1}\Big( \bigcup_{m=0}^\infty \Gamma_m(c') \Big).
$$
We especially call 
${\mathcal M}(c')$ a {\it decorated Mandelbrot set},
${\mathcal M}(c') \smallsetminus M 
= \Phi_M^{-1}\Big( \bigcup_{m=0}^\infty \Gamma_m(c') \Big)$ its 
{\it decoration} and $M \subset {\mathcal M}(c')$
the {\it main Mandelbrot set} of ${\mathcal M}(c')$. 
Also we call
${\mathcal K}_c(c')$ a {\it decorated filled Julia set} and
${\mathcal K}_c(c') \smallsetminus K(P_c)
= \Phi_c^{-1}\Big( \bigcup_{m=0}^\infty \Gamma_m(c') \Big)$ its 
{\it decoration}. We will apply the same terminologies to the images of 
${\mathcal M}(c')$ or
${\mathcal K}_c(c')$ by 
quasiconformal maps.
Note that the sets 
$\Gamma_m(c')~(m \ge 0)$,
$\cM(c')$ and 
$\cK_{c}(c')$
depend on the choice of $\rho'$ and $\rho$.
When we emphasize the dependence,
we denote them by 
$\Gamma_m(c')_{\rho', \rho}$, 
$\cM(c')_{\rho', \rho}$ and
$\cK_{c}(c')_{\rho', \rho}$, 
respectively.

\begin{figure}[htbp]
\includegraphics[width=0.95\textwidth]{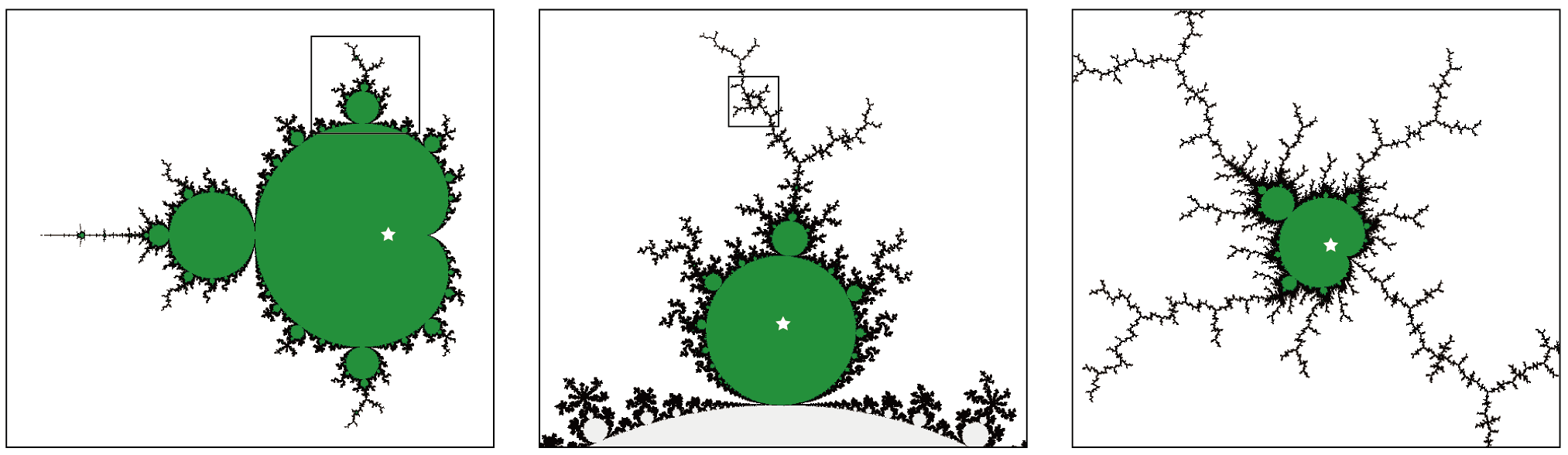}
\caption{\small
The original Mandelbrot set (left), a ``satellite" small Mandelbrot
set (middle),
and a ``primitive" small Mandelbrot set (right).
The stars indicate the central superattracting parameters.
}
\label{primitive and satellite small M-set}
\end{figure}

\noin
{\bf Small Mandelbrot sets.} 
When we zoom in the boundary of $M$, a lot of \lq\lq small Mandelbrot sets" 
appear and it is known that these sets are obtained as follows: (This is 
the result by Douady and Hubbard and its proof can be found in 
\cite[Th\'eor\`eme 1 du Modulation]{Haissinsky 2000}. See also 
\cite{Milnor 2000}.) 
Let $s_0 \ne 0$ be a {\it superattracting parameter}, that is, 
$P_{s_0}(z) = z^2 + s_0$ has a superattracting periodic point, and denote
its period by $p \geq 2$. Then there exists a unique small Mandelbrot set 
$M_{s_0}$ containing $s_0$ and a canonical homeomorphism 
$\chi : M_{s_0} \to M$ with $\chi(s_0) = 0$. Following Douady and Hubbard
we use the notation $s_0 \perp M$ and $s_0 \perp c_0$  to 
denote $M_{s_0} = \chi^{-1}(M)$ and $\chi^{-1}(c_0) \ (c_0 \in M)$, 
respectively. The set $M_{s_0}$ is called the {\it small Mandelbrot set 
with center $s_0$} (see Figure \ref{primitive and satellite small M-set}). 
If $c_1 := s_0 \perp c_0 \ (c_0 \in M)$, then $c_1$ is a parameter in 
$M_{s_0}$ which corresponds to $c_0 \in M$ and it is known that $P_{c_1}$ 
is renormalizable with period $p$ and $P_{c_1}^p$ is hybrid equivalent 
(see section 3) to $P_{c_0}$. 
We say $M_{s_0}$ is {\it primitive} if $K(P_{s_0 \perp (1/4)})$ has a parabolic
periodic point with a single petal. Otherwise we say $M_{s_0}$ is 
{\it satellite}, in which case $K(P_{s_0 \perp (1/4)})$ has a parabolic
periodic point with more than one petal. 
It is known that $M_{s_0}$ is primitive if the hyperbolic component
containing $s_0$ has a cusp on the boundary curve, and that $M_{s_0}$
is satellite if the hyperbolic component containing $s_0$ is attached
to another hyperbolic component at $s_0 \perp (1/4)$
(see Figure \ref{primitive and satellite small M-set}).

\fboxsep=0pt
\fboxrule=1pt
\begin{figure}[htbp]
\begin{center}
(i)\\[.5em]
\fbox{\includegraphics[width=.18\textwidth, bb = 0 0 1000 1002]{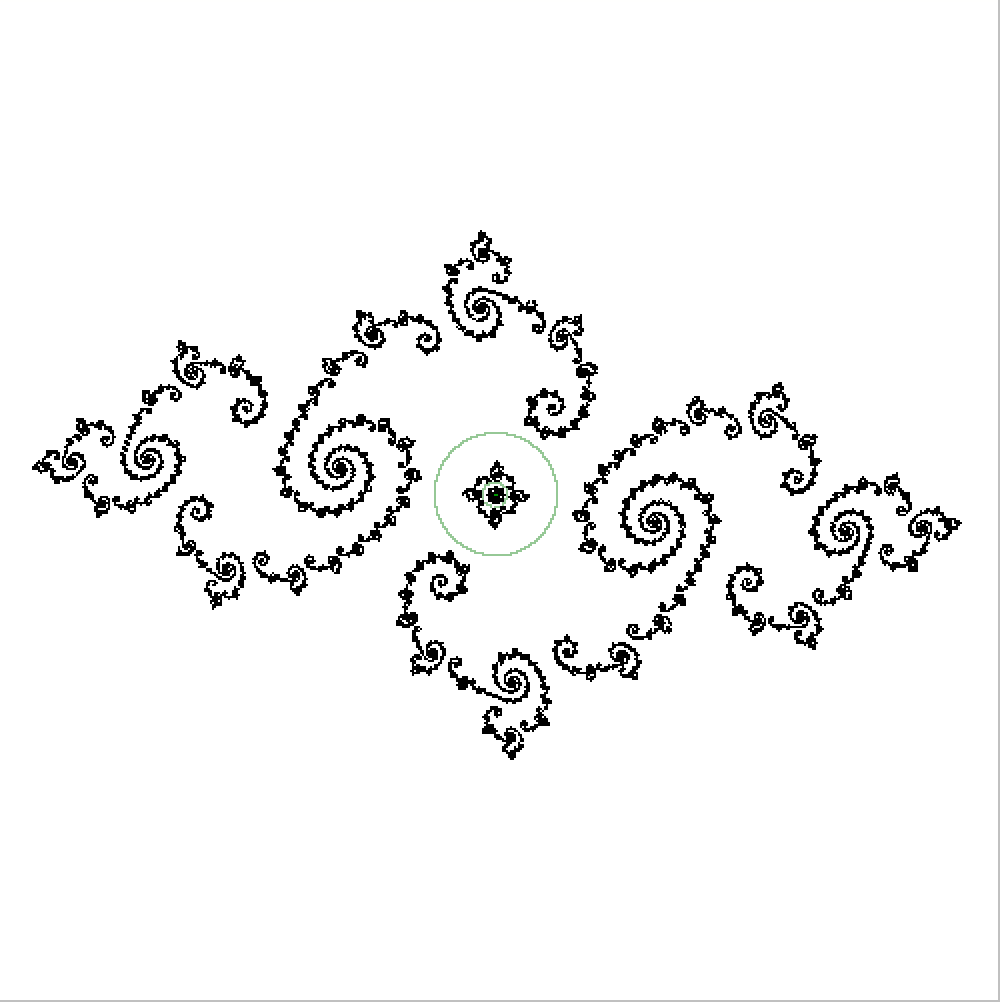}}
\fbox{\includegraphics[width=.18\textwidth, bb = 0 0 1002 1000]{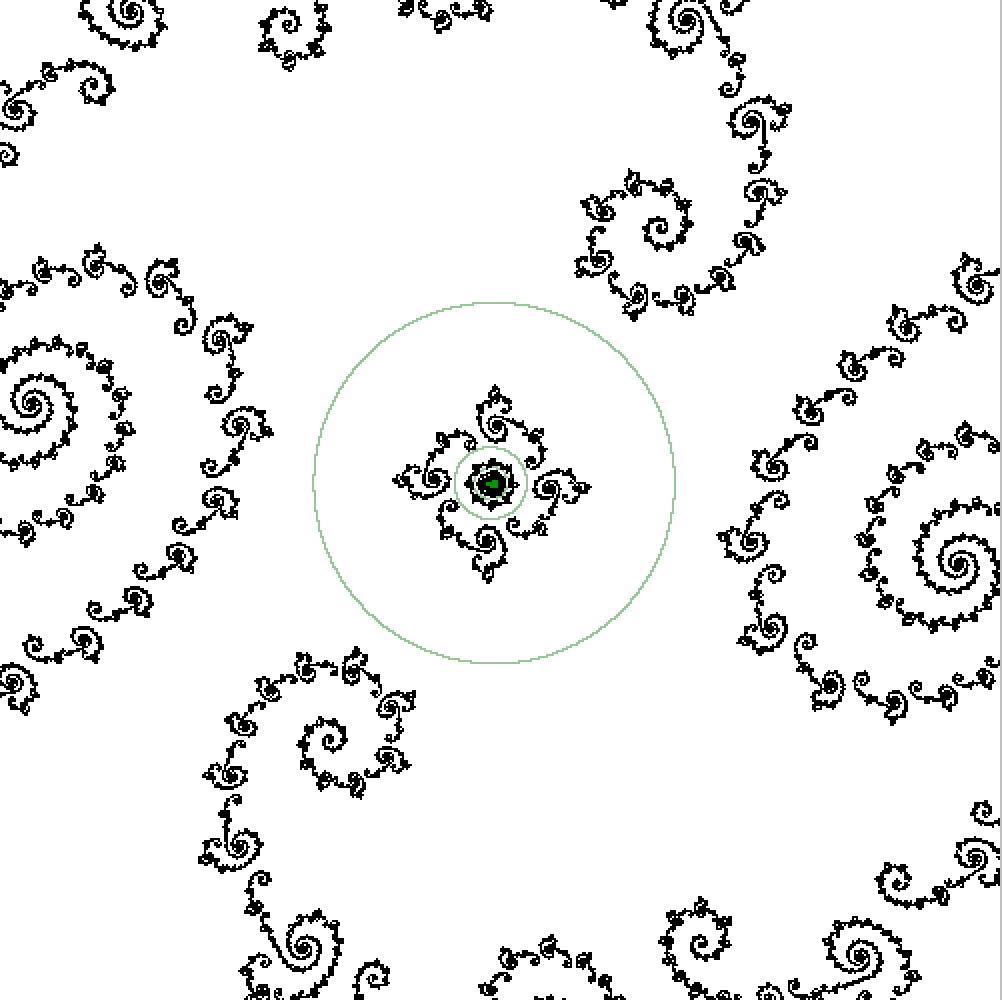}}
\fbox{\includegraphics[width=.18\textwidth, bb = 0 0 1002 998]{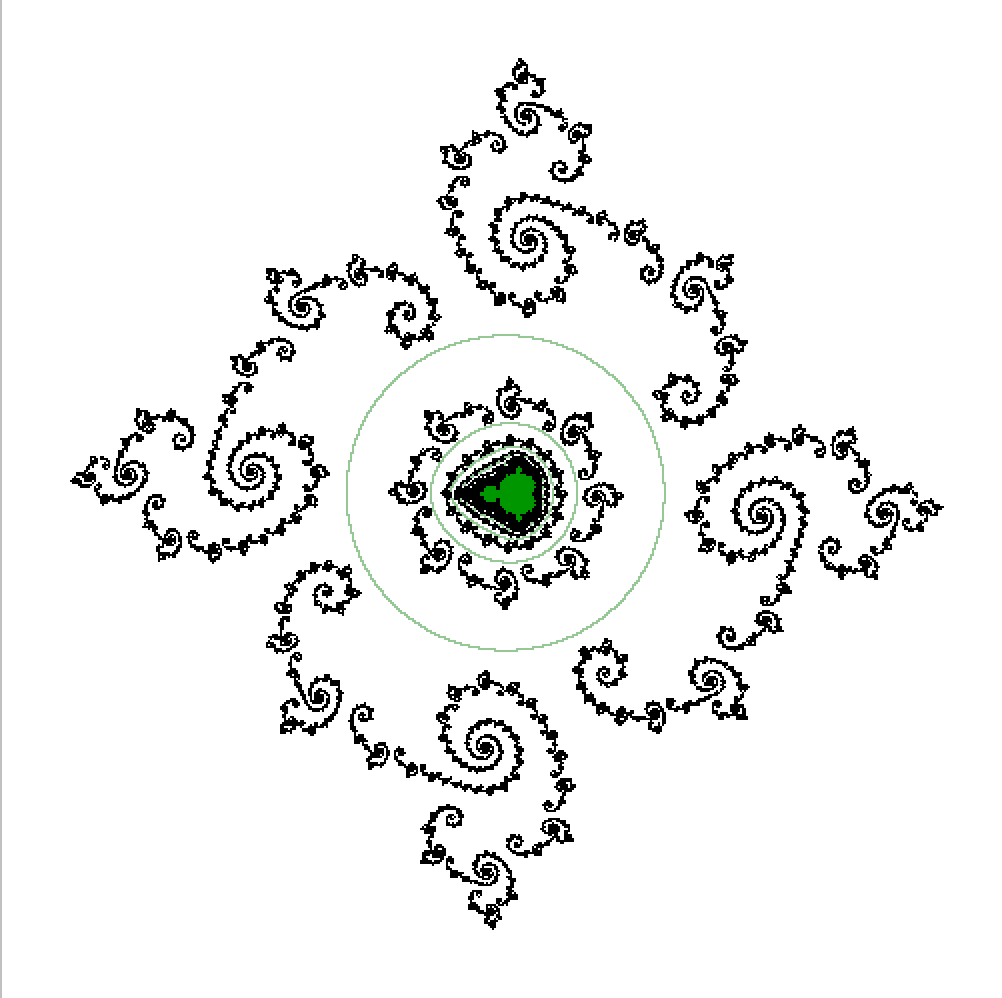}}
\fbox{\includegraphics[width=.18\textwidth, bb = 0 0 1002 1002]{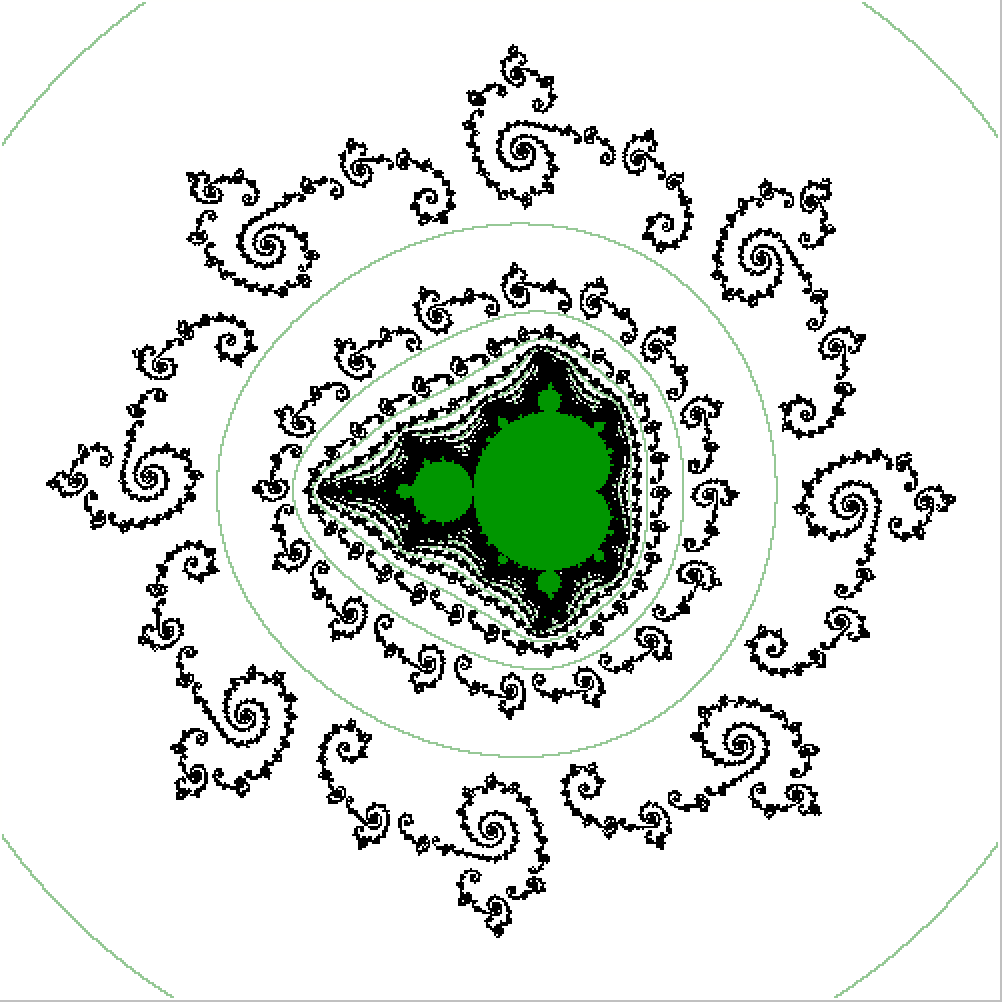}}
\fbox{\includegraphics[width=.18\textwidth, bb = 0 0 1002 1002]{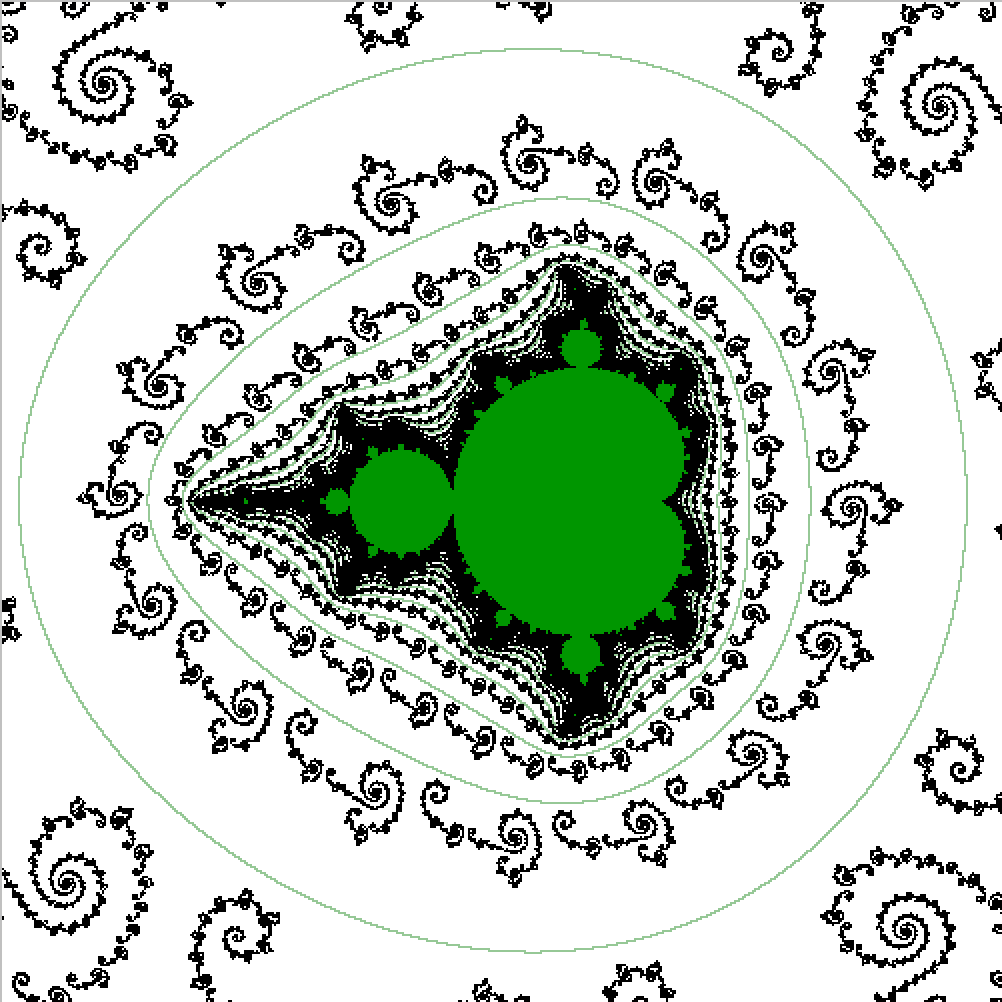}}
\\[.5em]
(ii) \\[.5em]
\fbox{\includegraphics[width=.18\textwidth, bb = 0 0 500 500]{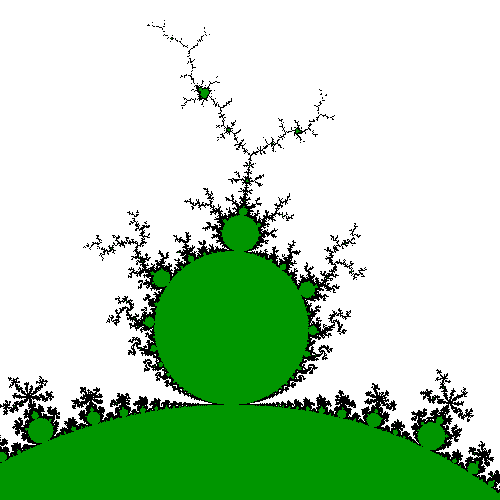}}
\fbox{\includegraphics[width=.18\textwidth, bb = 0 0 500 500]{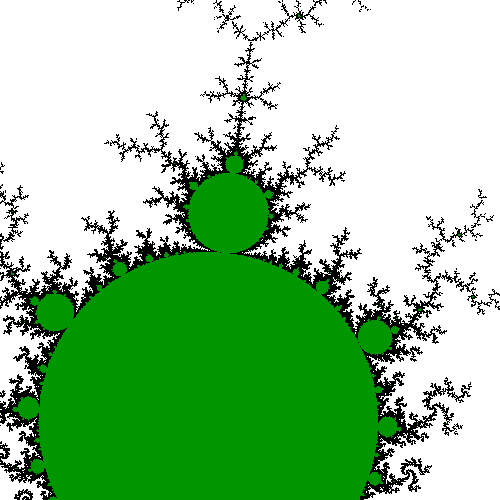}}
\fbox{\includegraphics[width=.18\textwidth, bb = 0 0 500 500]{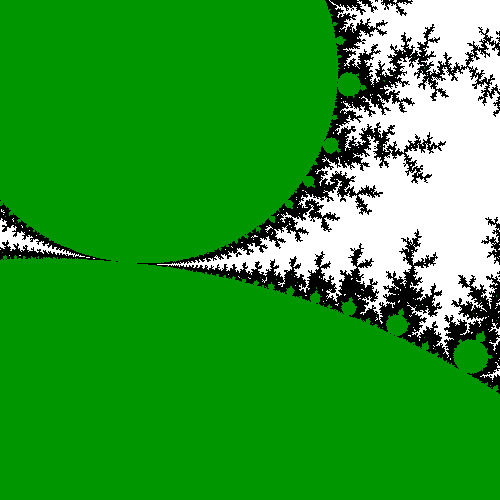}}
\fbox{\includegraphics[width=.18\textwidth, bb = 0 0 500 500]{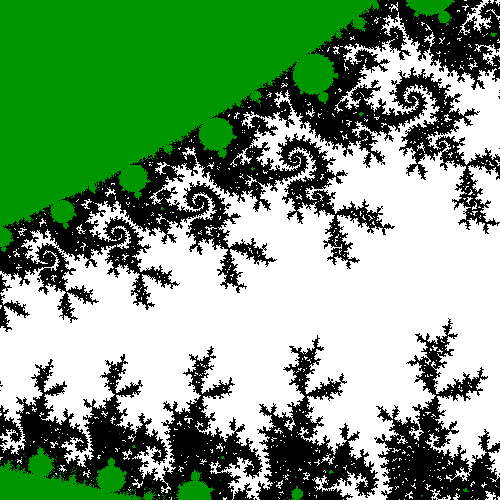}}
\fbox{\includegraphics[width=.18\textwidth, bb = 0 0 500 500]{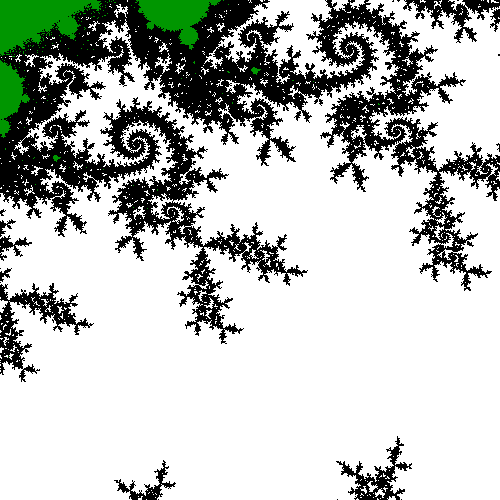}}
\\[.5em]
\fbox{\includegraphics[width=.18\textwidth, bb = 0 0 500 500]{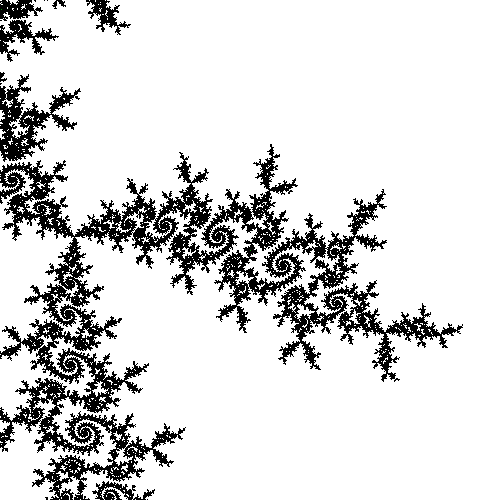}}
\fbox{\includegraphics[width=.18\textwidth, bb = 0 0 500 500]{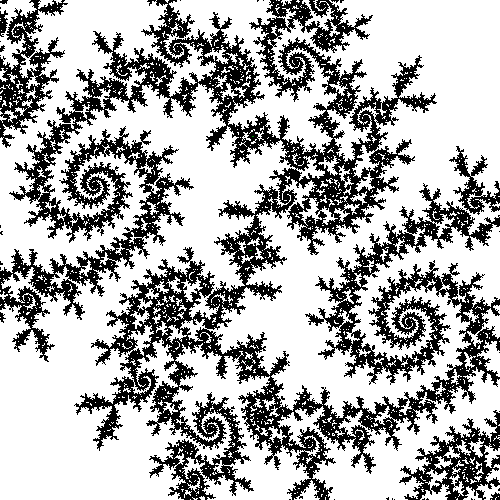}}
\fbox{\includegraphics[width=.18\textwidth, bb = 0 0 500 500]{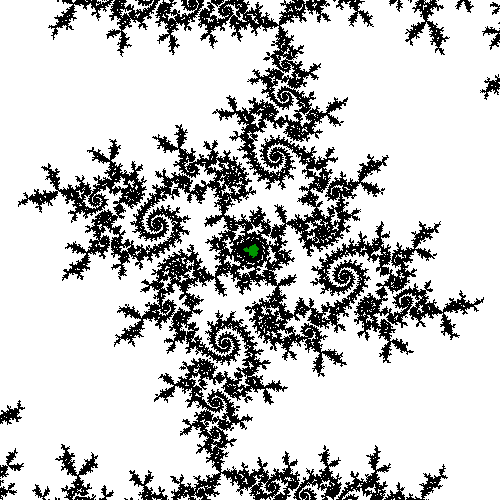}}
\fbox{\includegraphics[width=.18\textwidth, bb = 0 0 500 500]{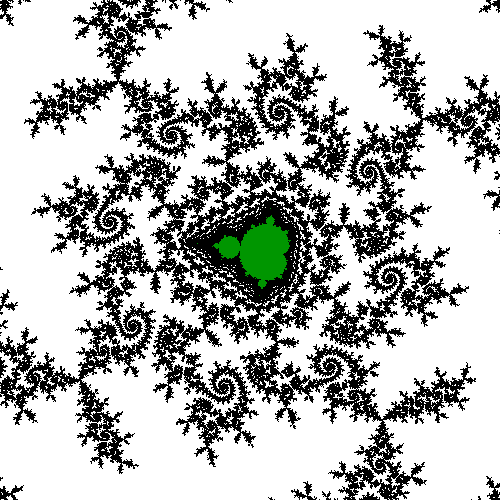}}
\fbox{\includegraphics[width=.18\textwidth, bb = 0 0 500 500]{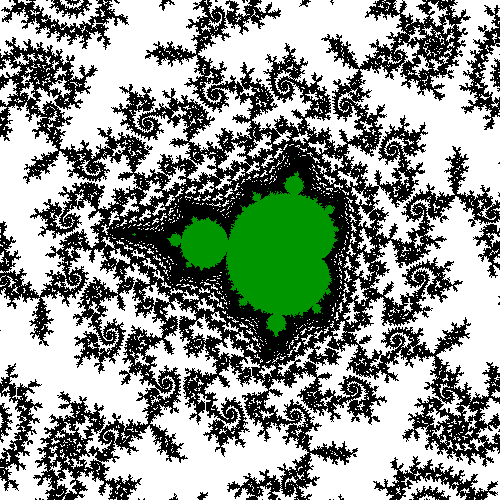}}
\\[.5em]
(iii) \\[.5em]
\fbox{\includegraphics[width=.18\textwidth, bb = 0 0 500 500]{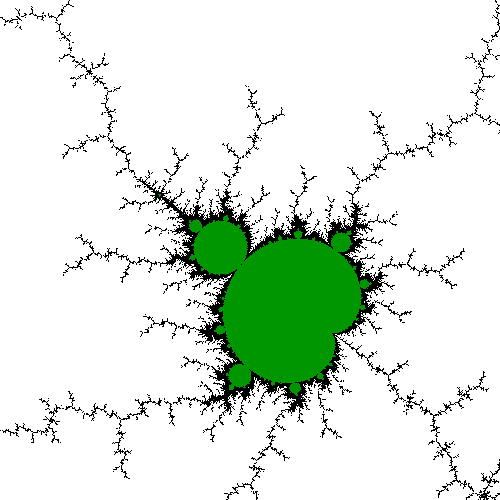}}
\fbox{\includegraphics[width=.18\textwidth, bb = 0 0 500 500]{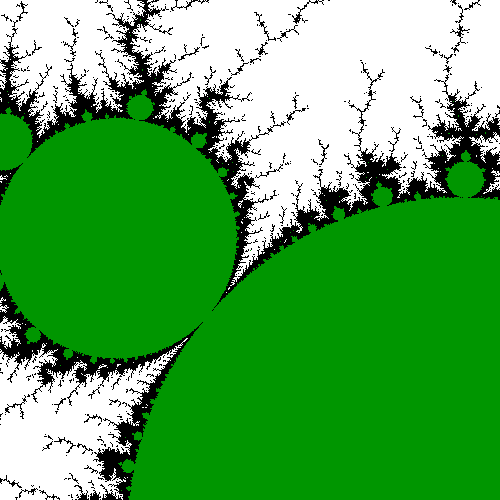}}
\fbox{\includegraphics[width=.18\textwidth, bb = 0 0 500 500]{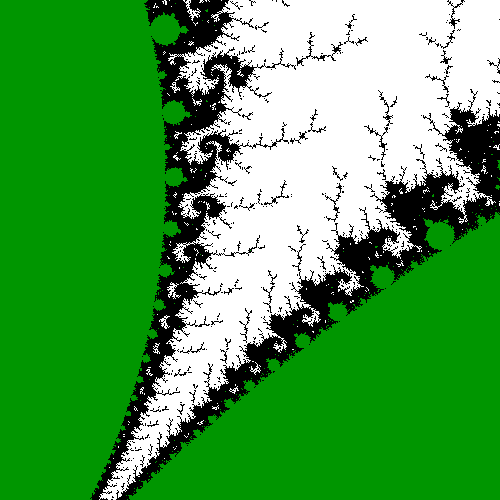}}
\fbox{\includegraphics[width=.18\textwidth, bb = 0 0 500 500]{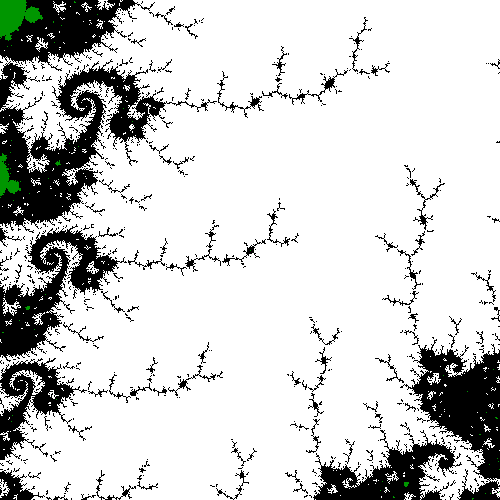}}
\fbox{\includegraphics[width=.18\textwidth, bb = 0 0 500 500]{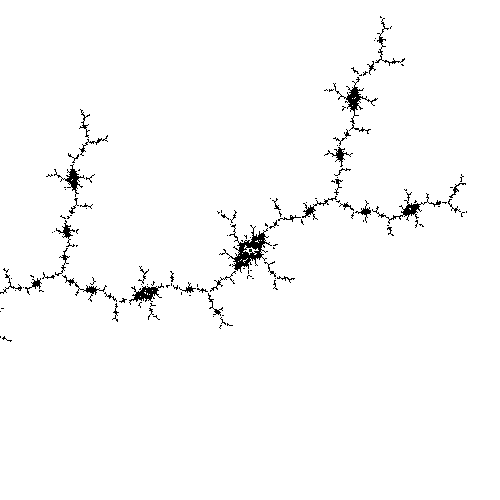}}
\\[.5em]
\fbox{\includegraphics[width=.18\textwidth, bb = 0 0 500 500]{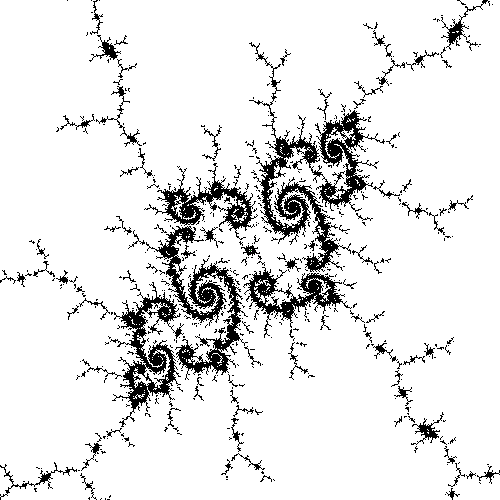}}
\fbox{\includegraphics[width=.18\textwidth, bb = 0 0 500 500]{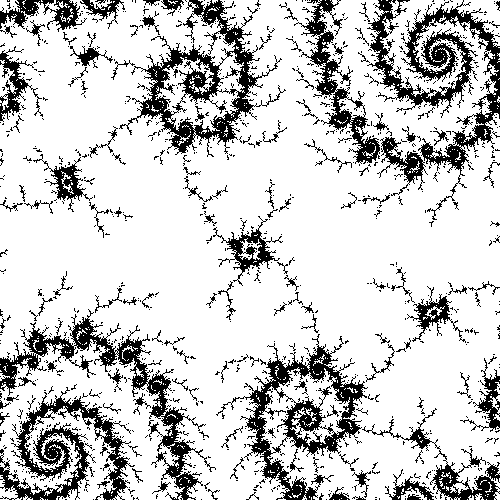}}
\fbox{\includegraphics[width=.18\textwidth, bb = 0 0 500 500]{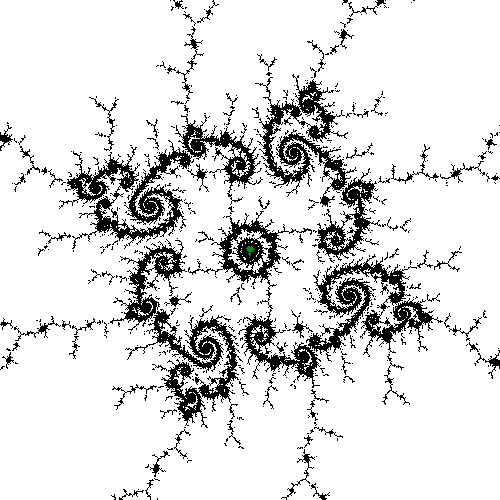}}
\fbox{\includegraphics[width=.18\textwidth, bb = 0 0 500 500]{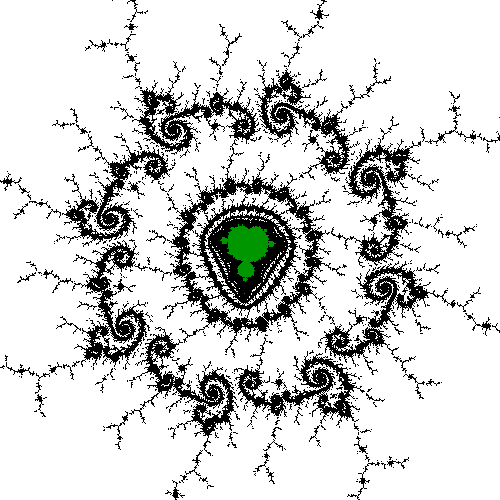}}
\fbox{\includegraphics[width=.18\textwidth, bb = 0 0 500 500]{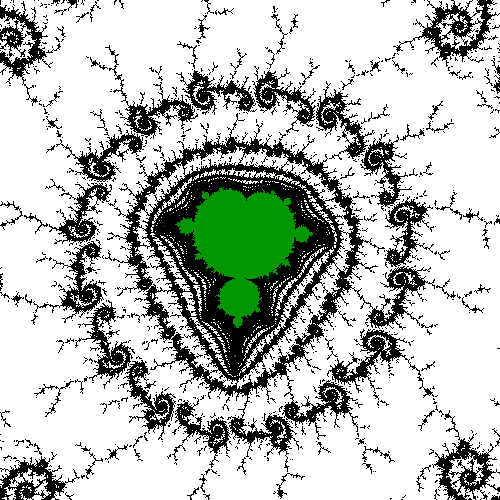}}
\end{center}
\caption{\small
(i): The decorated Mandelbrot set $\mathcal{M}(c')$ for $c'=-0.77+0.18 i$
(close to the parabolic parameter $c_0=-0.75$. 
(ii) and (iii): Embedded quasiconformal copies of $\mathcal{M}(c')$ 
above near the satellite/primitive small Mandelbrot sets 
in Figure \ref{primitive and satellite small M-set}.
}
\label{fig_satellite_primitive}
\end{figure}

\begin{defn}
Let $X$ and $Y$ be non-empty compact sets in $\C$. We say {\it $X$ 
appears ($K$-)quasiconformally in $Y$} or {\it $Y$ contains a 
($K$-)quasiconformal copy of $X$} if there is a ($K$-)quasiconformal map 
$\chi$ on a neighborhood of $X$ such that
$\chi(X) \subset Y$ and $\chi(\partial X) \subset \partial Y$. Note that
the condition $\chi(\partial X) \subset \partial Y$ is to exclude the case
$\chi(X) \subset \text{int}(Y)$.
\end{defn}

Now our results are as follows:

\begin{thmA*}[\bf Julia sets appear quasiconformally in $M$]
Let $M_{s_0}$ be any small Mandelbrot set, where $s_0 \ne 0$ is a
superattracting parameter and $c_0 \in \partial M$ 
any Misiurewicz or parabolic parameter. Then
for every small $\vep > 0$ and $\vep' > 0$, there exists an $\eta \in \C$ 
with $|\eta| < \vep$ and
$c_0+\eta \notin M$ such that
${\mathcal M}(c_0+\eta)$ appears quasiconformally in $M$ in the neighborhood
$D(s_0 \perp c_0, \vep')$ of $s_0 \perp c_0$. In particular, the Cantor
Julia set $J(P_{c_0+\eta})$ appears quasiconformally in $M$.
\end{thmA*}

\noin
Theorem A shows the following: Take any small Mandelbrot set $M_{s_0}$
and zoom in any small neighborhood of $s_0 \perp c_0 \in M_{s_0}$, then we 
can find a quasiconformal image of ${\mathcal M}(c_0+\eta)$. That is, as we 
zoom in, first we observe a quasiconformal image of $J(P_{c_0+\eta})$, which 
corresponds to the $\Phi_M^{-1}$ image of the rescaled Cantor Julia set 
$\Gamma_0(c_0+\eta)$ in ${\mathcal M}(c_0+\eta)$ and its iterated preimages
(decoration) by $z \mapsto z^2$ and finally the main Mandelbrot set of the
quasiconformal image of ${\mathcal M}(c_0+\eta)$, say $M_{s_1}$, appears.

Figure \ref{figures of a primitive-Misiurewicz case} shows zooms around
a Misiurewicz parameter $c_1 = s_0 \perp c_0$ in a primitive small 
Mandelbrot set $M_{s_0}$ 
(Figure \ref{figures of a primitive-Misiurewicz case}--(1)). 
After a sequence of nested structures, a smaller 
\lq\lq small Mandelbrot set" $M_{s_1}$ appears 
(Figure \ref{figures of a primitive-Misiurewicz case}--(15)). 
Here $M_{s_0}$ is the relatively big \lq\lq small 
Mandelbrot set" which is located in the upper right part of $M$. The map 
$P_{s_0}$ has a superattracting periodic point of period 4 and $c_0$ is the 
Misiurewicz parameter which satisfies $P_{c_0}(P_{c_0}^4(0)) = P_{c_0}^4(0)$
and corresponds to the \lq\lq junction of three roads" as shown in
Figure \ref{primitive and satellite small M-set} in the middle.

Since both the set of Misiurewicz parameters and the set of parabolic parameters form dense subsets of $\partial M$,
we can reformulate Theorem A as follows: 
\begin{thmA'*}
Let $M_{s_0}$ be any small Mandelbrot set, where $s_0 \ne 0$ is a
superattracting parameter and $c_0 \in \partial M$ 
any parameter. Then
for every small $\vep > 0$ and $\vep' > 0$, there exists an $\eta \in \C$
with $|\eta| < \vep$ and $c_0+\eta \notin  M$ such that
${\mathcal M}(c_0+\eta)$ appears quasiconformally in $M$ in the 
neighborhood $D(s_0 \perp c_0, \vep')$. In particular, the Cantor
Julia set $J(P_{c_0+\eta})$ appears quasiconformally in $M$. 
\end{thmA'*}

Next we show that the same decoration of $\cM(c_0+\eta)$ in Theorem A
appears quasiconformally also in some filled Julia sets.

\begin{thmB*}[\bf Decoration in filled Julia sets]
Let $M_{s_1}$ denote the main Mandelbrot set of the quasiconformal image of
${\mathcal M}(c_0+\eta)$ in Theorem A. Then for every 
$c \in M$, ${\mathcal K}_{c}(c_0+\eta)$ appears 
quasiconformally in $K(P_{s_1 \perp c})$, where $s_1 \perp c \in M_{s_1}$.
\end{thmB*}

\noin
Theorem B shows the following: Choose any parameter from the main Mandelbrot
set $M_{s_1}$ in the quasiconformal image of ${\mathcal M}(c_0+\eta)$,
that is, choose any $c \in M$ and consider $s_1 \perp c$ $\in M_{s_1}$
and zoom in the neighborhood of $0 \in K(P_{s_1 \perp c})$. Then we can find
a quasiconformal image of ${\mathcal K}_{c}(c_0+\eta)$, whose decoration is
conformally the same as that of ${\mathcal M}(c_0+\eta)$.

Next we show that there are smaller Mandelbrot sets and their 
decorations which are images of model sets by quasiconformal maps
whose dilatations are arbitrarily close to 1.

\begin{thmC*}[\bf Almost conformal copies]
Let $c_0$ be any Misiurewicz or parabolic parameter
and $B$ any small closed disk whose interior intersects with $\partial M$.
Then for any small $\vep > 0$ and $\kappa>0$,
there exist an $\eta \in \C$ with $|\eta| < \vep$ and two positive numbers
$\rho$ and $\rho'$ with $\rho' < \rho$ 
such that $c_0 +\eta \notin M$ and $\cM(c_0+\eta)_{\rho', \rho}$ 
appears $(1+\kappa)$-quasiconformally 
in $B \cap M$.
In particular, $B \cap \partial M$ contains a $(1+\kappa)$-quasiconformal
copy of the Cantor Julia set $J(P_{c_0+\eta})$.
\end{thmC*}
\noin
One can formulate this theorem as in Theorem A' but we skip the details.

\begin{defn}[\bf Semihyperbolicity]
A quadratic polynomial $P_c(z) = z^2 + c$ (or the parameter $c$) is 
called {\it semihyperbolic} if  
\begin{itemize}
\item[(1)]  the critical point $0$ is non-recurrent, that is, 
$0 \notin \omega(0)$, where $\omega(0)$ is the $\omega$-limit set of 
the critical point $0$ and

\item[(2)] $P_c$ has no parabolic periodic points.
\end{itemize}
\end{defn}

\noin
It is easy to see that if $P_c$ is hyperbolic then it is semihyperbolic.
If $P_c$ is semihyperbolic, then it is known that it has no Siegel
disks and Cremer points 
(\cite{Mane 1993}, \cite{Carleson-Jones-Yoccoz 1994}). Also it is not 
difficult to see that $J(P_c)$ is measure 0 from the result by Lyubich 
(\cite{Lyubich 1991}) and Shishikura. (This also follows from 
\cite[p.2, Theorem 1.1]{Carleson-Jones-Yoccoz 1994}.) Thus the 
semihyperbolic dynamics is relatively understandable. A typical 
semihyperbolic but non-hyperbolic parameter $c$ is a Misiurewicz 
parameter. But there seems less concrete examples of semihyperbolic
parameter $c$ which is neither hyperbolic nor Misiurewicz. Next 
corollary shows that we can visually identify these parameters 
everywhere in $\partial M$.

\begin{corD*}[\bf Abundance of semihyperbolicity]
For every parameter $c$ belonging to the quasiconformal image
of the decoration of ${\mathcal M}(c_0+\eta)$ in Theorem A, 
$P_c$ is semihyperbolic. Also the set of semihyperbolic parameters which 
are not Misiurewicz and non-hyperbolic is dense in $\partial M$. 
\end{corD*}

\noindent
Corollary D together with Theorem C explains the following famous result
by Shishikura:

\begin{thm*}[Shishikura, 1998]
Let
$$
  SH := \{ c \in \partial M \ | \ P_c \text{ is semihyperbolic} \},
$$
then the Hausdorff dimension of $SH$ is 2. In particular, the Hausdorff
dimension of the boundary of $M$ is 2.
\end{thm*}

\noin
{\bf Explanation.} \ 
Since there exist quadratic Cantor Julia sets with Hausdorff dimension
arbitrarily close to 2 and we can find such parameters in every 
neighborhood of a point in $\partial M$ 
(\cite[p.231, {\it proof of Theorem B} and p.232, {\it Remark} 1.1 (iii)]{Shishikura 1998}), we can find an 
$\eta$ such that $\dim_H(J(P_{c_0+\eta}))$ is
arbitrarily close to 2. Then by Theorem C and Corollary D it follows that
we can find a subset of $\partial M$ with Hausdorff dimension
arbitrarily close to 2 and consisting of semihyperbolic parameters 
as a quasiconformal image of the decoration of ${\mathcal M}(c_0+\eta)$.
This implies that $\dim_H(SH) = 2$. 
\QED


\begin{figure}[htbp]
\hskip -40mm
{\small (1)}
\hskip 43mm
{\small (2)}
\hskip 43mm
{\small (3)}

\includegraphics[scale=0.19, bb = 0 0 640 480]{fig_M-4dendrite-015.jpg} \hskip 5mm
\includegraphics[scale=0.19, bb = 0 0 640 480]{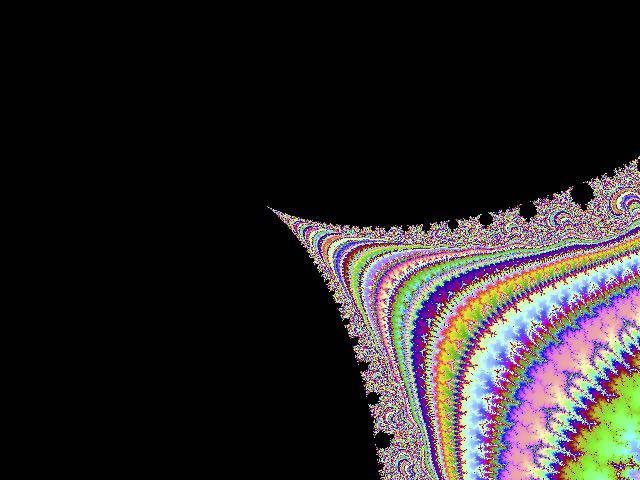} \hskip 5mm
\includegraphics[scale=0.19, bb = 0 0 640 480]{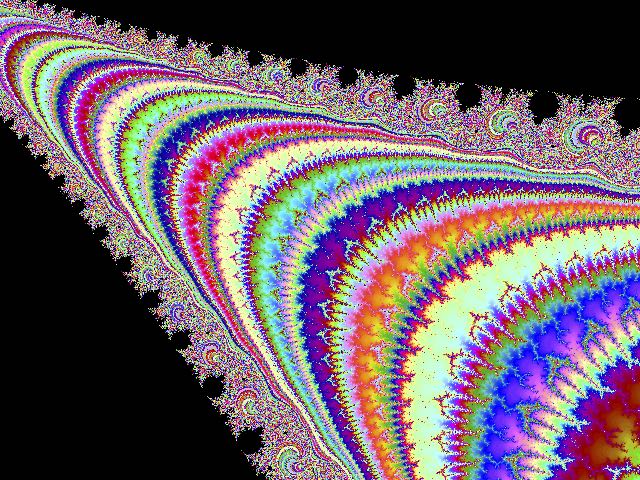} \hskip 5mm

\hskip -40mm
{\small (4)}
\hskip 43mm
{\small (5)}
\hskip 43mm
{\small (6)}

\includegraphics[scale=0.19, bb = 0 0 640 480]{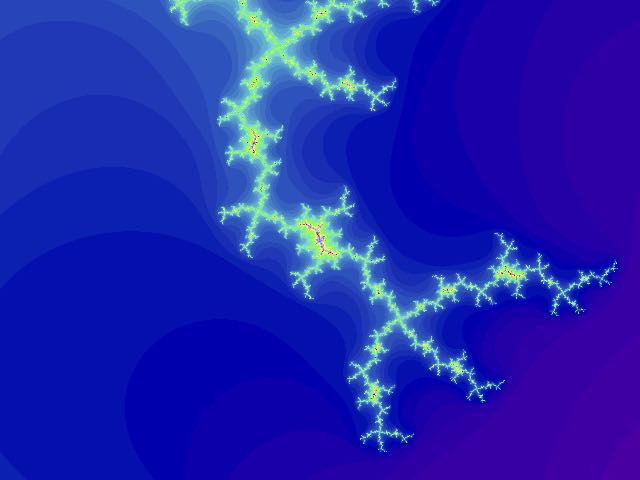} \hskip 5mm
\includegraphics[scale=0.19, bb = 0 0 640 480]{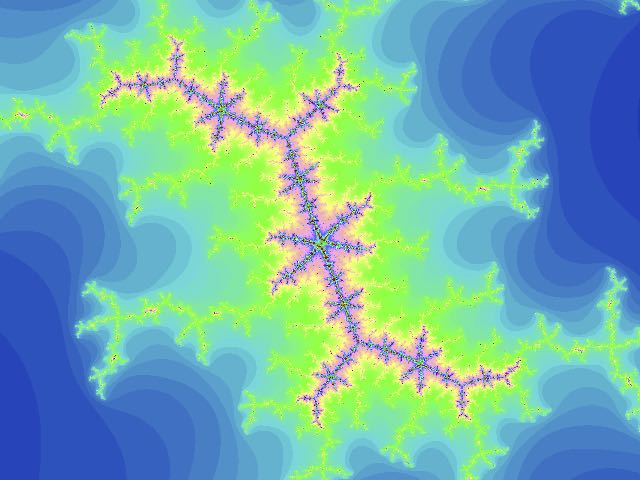} \hskip 5mm
\includegraphics[scale=0.19, bb = 0 0 640 480]{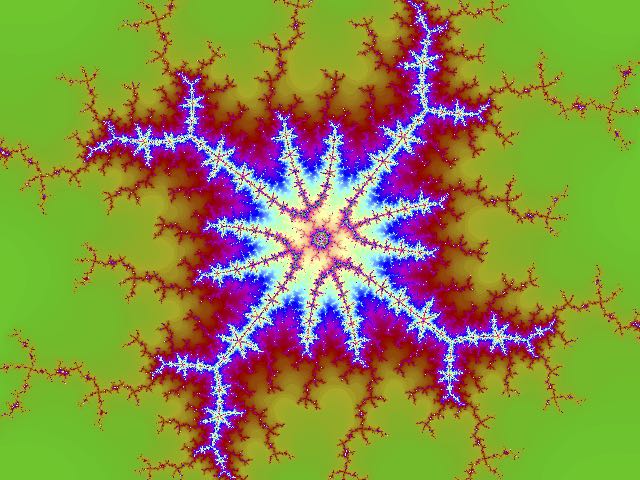} \hskip 5mm

\hskip -40mm
{\small (7)}
\hskip 43mm
{\small (8)}
\hskip 43mm
{\small (9)}

\includegraphics[scale=0.19, bb = 0 0 640 480]{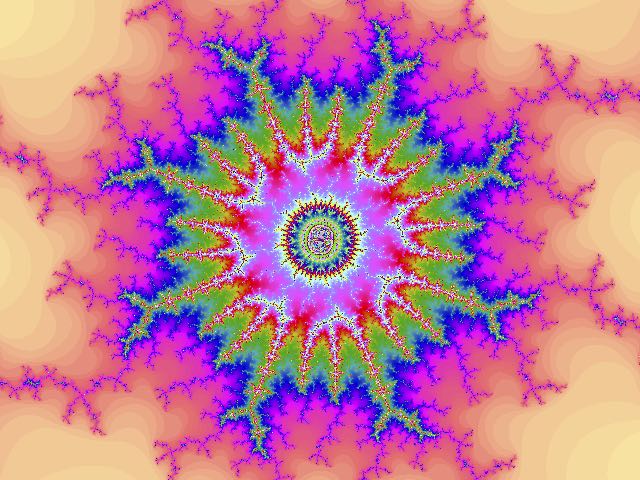} \hskip 5mm
\includegraphics[scale=0.19, bb = 0 0 640 480]{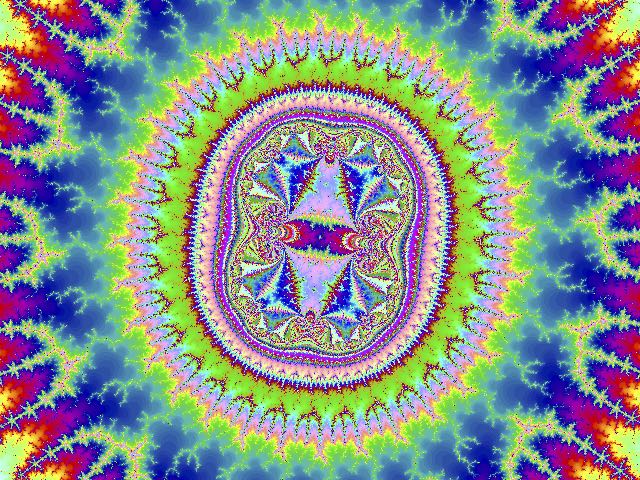} \hskip 5mm
\includegraphics[scale=0.19, bb = 0 0 640 480]{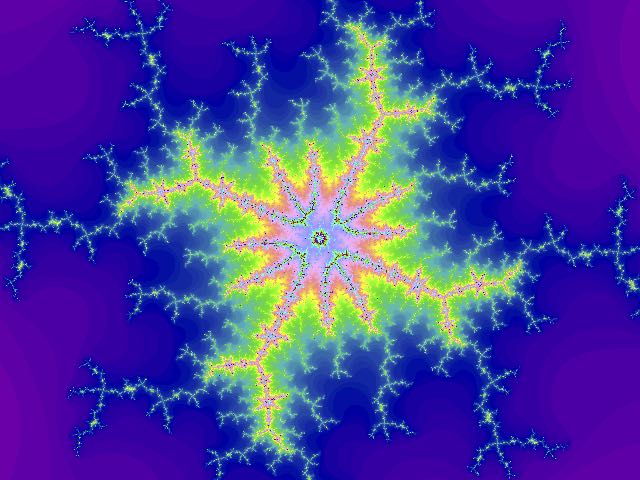} \hskip 5mm

\hskip -38mm
{\small (10)}
\hskip 41mm
{\small (11)}
\hskip 41mm
{\small (12)}

\includegraphics[scale=0.19, bb = 0 0 640 480]{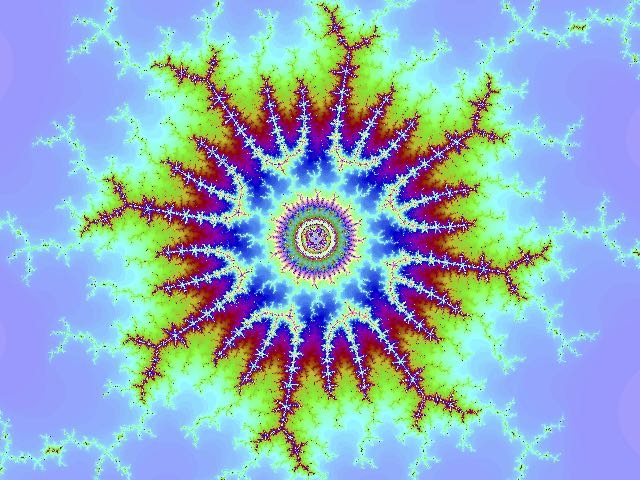} \hskip 5mm
\includegraphics[scale=0.19, bb = 0 0 640 480]{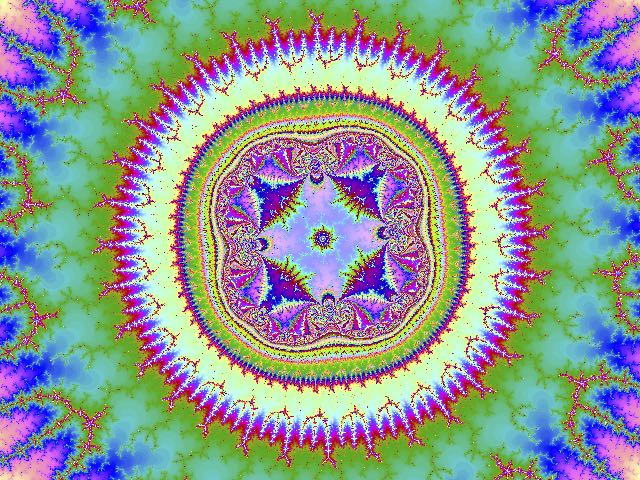} \hskip 5mm
\includegraphics[scale=0.19, bb = 0 0 640 480]{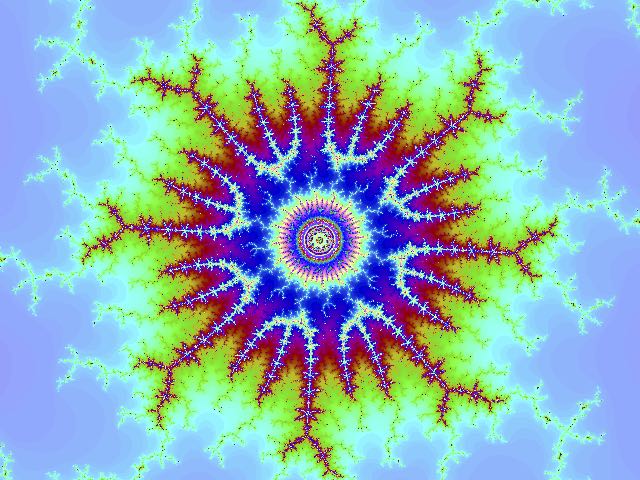} \hskip 5mm

\hskip -38mm
{\small (13)}
\hskip 41mm
{\small (14)}
\hskip 41mm
{\small (15)}

\includegraphics[scale=0.19, bb = 0 0 640 480]{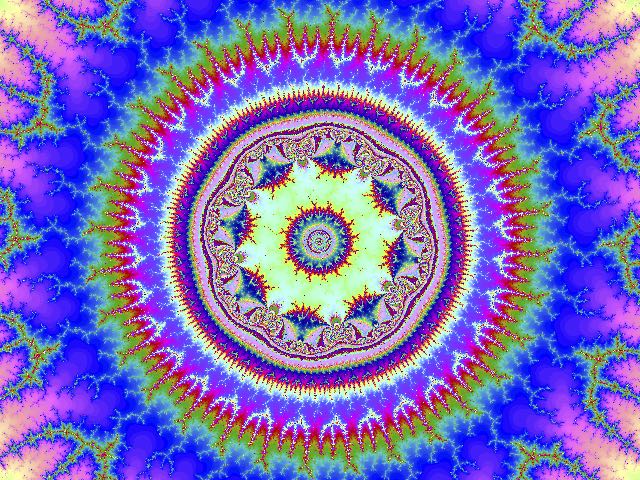} \hskip 5mm
\includegraphics[scale=0.19, bb = 0 0 640 480]{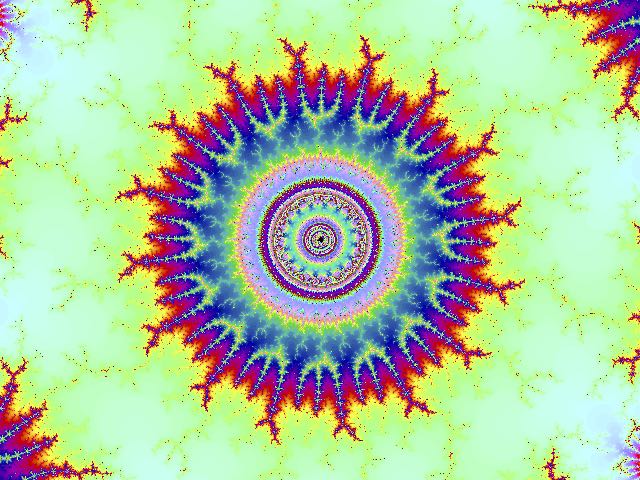} \hskip 5mm
\includegraphics[scale=0.19, bb = 0 0 640 480]{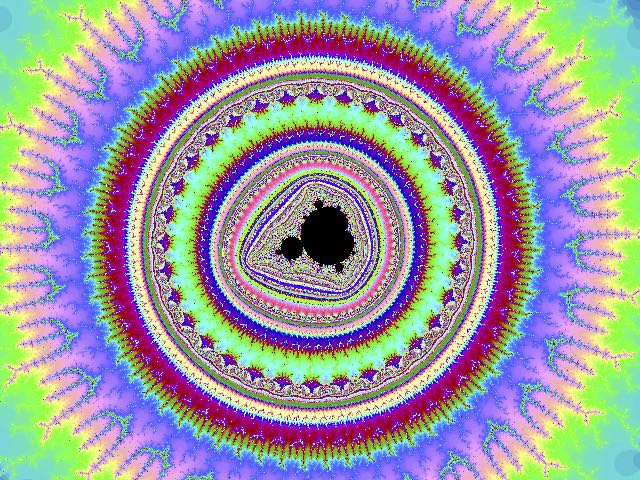} \hskip 5mm

\caption{\small Zooms around a parabolic point $s_1 \perp c_1$ 
in a primitive small Mandelbrot set $M_{s_1}$. After a sequence of 
complicated nested structures, another smaller Mandelbrot set 
$M_{s_2}$ appears ((15)).}
\label{figures of a complicated structure}
\end{figure}

\begin{rem*}
(1) \ 
A similar result to Theorem B still holds even when $c \in \C \smallsetminus M$
is sufficiently close to $M$. Actually when $M_{s_1}$ is primitive (resp. satellite),
the homeomorphism $\chi : M_{s_1} \to M$ can be extended to a homeomorphism between
some neighborhoods of $M_{s_1}$ and $M$ 
(resp. $M_{s_1} \smallsetminus \{ s_1 \perp (1/4) \}$ and 
$M \smallsetminus \{ 1/4 \}$)
and so $s_1 \perp c = \chi^{-1}(c)$ can be still defined for such a 
$c \in \C \smallsetminus M$. Then in this case, by modifying the definition 
of ${\mathcal K}_{c}(c_0+\eta)$ for this $c \in \C \smallsetminus M$ we can 
prove that a \lq\lq ${\mathcal K}_{c}(c_0+\eta)$" appears quasiconformally
in $K(P_{s_1 \perp c})$, where $s_1 \perp c \in \C \smallsetminus M_{s_1}$ is a
point which is sufficiently close to $M_{s_1}$. We omit the details.

\noin
(2) \ 
Corollary D and the above Theorem (Shishikura, 1998) means that relatively
understandable dynamics is abundant in $\partial M$ (provided that Lebesgue
measure of $\partial M$ is 0). 
In \cite[p.225, THEOREM A]{Shishikura 1998}, Shishikura actually proved
$\dim_H(SH) = 2$, which immediately implies $\dim_H(\partial M) = 2$.
A new point of our \lq\lq explanation" is that we constructed a decoration 
in $M$ which contains a quasiconformal image of a {\it whole} Cantor Julia set
and consists of semihyperbolic parameters. 
So now we can say that \lq\lq$\dim_H(\partial M) = 2$ holds, because
we can see a lot of almost conformal images of Cantor quadratic Julia sets
whose Hausdorff dimension are arbitrarily close to 2".

\noin
(3) \ 
Take a small Mandelbrot set $M_{s_1}$ (e.g. 
Figure \ref{figures of a primitive-Misiurewicz case}--(15) 
= Figure \ref{figures of a complicated structure}--(1)) 
and another Misiurewicz or parabolic parameter $c_*$ (e.g. $c_* = 1/4$ in
Figure 6) and zoom in the neighborhood of $s_1 \perp c_*$. Then we
see much more complicated structure than we expected as follows:
According to Theorem A, by replacing $s_0$ with $s_1$ and $c_0$ with $c_*$, 
it says that ${\mathcal M}(c_* + \eta)$
appears quasiconformally in $D(s_1 \perp c_*, \vep')$. This means that as 
we zoom in, we first see a quasiconformal image of $J(P_{c_*+\eta_*})$, say 
$\wt{J}_{c_*+\eta_*}$ (e.g. \lq\lq broken cauliflower", when $c_* = 1/4$). 
But in reality as we zoom in, what we first see is a $\wt{J}_{c_0+\eta_0}$ 
(e.g. \lq\lq broken dendrite". 
See Figure \ref{figures of a complicated structure}--(5)).
This seems to contradict with Theorem A, but actually it does not. As 
we zoom in further in the middle part of $\wt{J}_{c_0+\eta_0}$,  we
see iterated preimages of $\wt{J}_{c_0+\eta_0}$ by $z \mapsto z^2$ 
(Figure \ref{figures of a complicated structure}--(6), (7)) 
and then 
$\wt{J}_{c_*+\eta_*}$ appears 
(Figure \ref{figures of a complicated structure}--(8)). After that we 
see again iterated preimages of $\wt{J}_{c_0+\eta_0}$ by $z \mapsto z^2$ 
(Figure \ref{figures of a complicated structure}--(9), (10)) 
and then a once iterated preimage of $\wt{J}_{c_*+\eta_*}$ appears
(Figure \ref{figures of a complicated structure}--(11)). This complicated 
structure continues and finally, we see a smaller Mandelbrot set,
say $M_{s_2}$ (Figure \ref{figures of a complicated structure}--(15)). 
We can explain this complicated phenomena as follows: What we see
in the series of magnifications above is a quasiconformal image of 
$\cM (\cK_{c_*+\eta_*}(c_0+\eta_0))$, where
\begin{eqnarray*}
\cM (\cK_{c_*+\eta_*}(c_0+\eta_0))
& := &
M \cup 
\Phi_M^{-1}\Big( 
       \bigcup_{m=0}^\infty \Gamma_m(\cK_{c_*+\eta_*}(c_0+\eta_0)) \Big), \\
\Gamma_m(\cK_{c_*+\eta_*}(c_0+\eta_0)) 
& := &
\text{the inverse image of } \Gamma_0(\cK_{c_*+\eta_*}(c_0+\eta_0)) \ \text{by} \
z \mapsto z^{2^m}. 
\end{eqnarray*}
Here $\cM (\cK_{c_*+\eta_*}(c_0+\eta_0))$ is obtained just by replacing
$J(P_{c'})$ with $\cK_{c_*+\eta_*}(c_0+\eta_0)$ in the definition of 
$\cM(c')$. Although $c_*+\eta_* \notin M$, 
$\cK_{c_*+\eta_*}(c_0+\eta_0)$ can be defined in the similar manner. 
See the Remark (1) above. So what we first 
see as we zoom in the neighborhood of $s_1 \perp c_*$ is a quasiconformal image of
$\cK_{c_*+\eta_*}(c_0+\eta_0)$, whose outer most part is 
$\wt{J}_{c_0+\eta_0}$ ($=$ broken dendrite) and inner most part is 
$\wt{J}_{c_*+\eta_*}$ ($=$ broken cauliflower). As we zoom in further, 
we see quasiconformal image of the preimage of $\cK_{c_*+\eta_*}(c_0+\eta_0)$ by 
$z \mapsto z^2$, whose inner most part is a once iterated preimage of 
$\wt{J}_{c_0+\eta_0}$. After we see successive preimages of 
$\cK_{c_*+\eta_*}(c_0+\eta_0)$ by $z \mapsto z^2$, a much more smaller Mandelbrot
set $M_{s_2}$ finally appears. Since $\cK_{c_*+\eta_*}(c_0+\eta_0)$
itself has a nested structure, the total picture has this very 
complicated structure. The proof is completely the same as for the 
Theorem A. 

\par
\noin
(4) \
Theorem A is an extension of a Douady's result
but his original formulation is a little more general
(see THEOREM 1 and THEOREM 2 of \cite[pp.22-23]{Douady 2000}). 
It is possible to state our results in his formulation.
See the remark in Section 5.
 
\end{rem*}

\section{The quadratic-like maps and the Mandelbrot-like family}

In this section, we briefly recall the definitions of the quadratic-like map and
the Mandelbrot-like family and explain the key Proposition \ref{D-BDS Proposition} 
which is crucial for the proof of Theorem A.

A map $h : U' \to U$ is called {\it a polynomial-like map} if
$U', \ U \subset \C$ are topological disks with $U' \Subset U$ (which means 
$\overline{U'} \subset U$) 
and $h$ is holomorphic and proper map of degree $d$ with respect to $z$. 
It is called a {\it quadratic-like map} when $d=2$.
The {\it filled Julia set} $K(h)$ and the {\it Julia set} $J(h)$ of a 
polynomial-like map $h$ are defined by
\begin{eqnarray*}
  K(h) 
& :=  &
\{ z \in U' \ | \ h^n(z) \text{ are defined for every } n \in \N \}
= \bigcap_{n=0}^\infty h^{-n}(U'), \\
  J(h)
& := &
\partial K(h).
\end{eqnarray*}
The famous straightening theorem by Douady and Hubbard 
(\cite[p.296, THEOREM 1]{Douady-Hubbard 1985}) says that
every polynomial-like map $h : U' \to U$ of degree $d$ is quasiconformally 
conjugate to a polynomial $P$ of degree $d$. More precisely $h$ is 
{\it hybrid equivalent} to $P$, that is, there exists a quasiconformal map
$\phi$ sending a neighborhood of $K(h)$ to a neighborhood of $K(P)$ such
that $\phi \circ h = P \circ \phi$ and $\overline{\partial} \phi = 0 \ \text{a.e.}$
on $K(h)$.
Also if $K(h)$ is connected, then $P$ is unique up to conjugacy by an affine map.

A family of holomorphic maps $\boldsymbol{h} = \{ h_\lambda \}_{\lambda \in W}$ 
is called a {\it Mandelbrot-like family} if the 
following (1)--(8) hold:

\begin{itemize}
\setlength{\itemsep}{5pt}
\item[(1)]
$W \subset \C$ is a Jordan domain with $C^1$ boundary $\partial W$.

\item[(2)]
There exists a family of maps 
$\Theta = \{ \Theta_\lambda \}_{\lambda \in W}$ such that
for every $\lambda \in W$, $\Theta_\lambda : \overline{A(R, R^2)} \to \C$ is a 
quasiconformal embedding and that $\Theta_\lambda(Z)$ is holomorphic in 
$\lambda$ for every $Z \in \overline{A(R, R^2)}$.

\begin{figure}[htbp]
\hskip 1cm
\includegraphics[width=.65\textwidth]{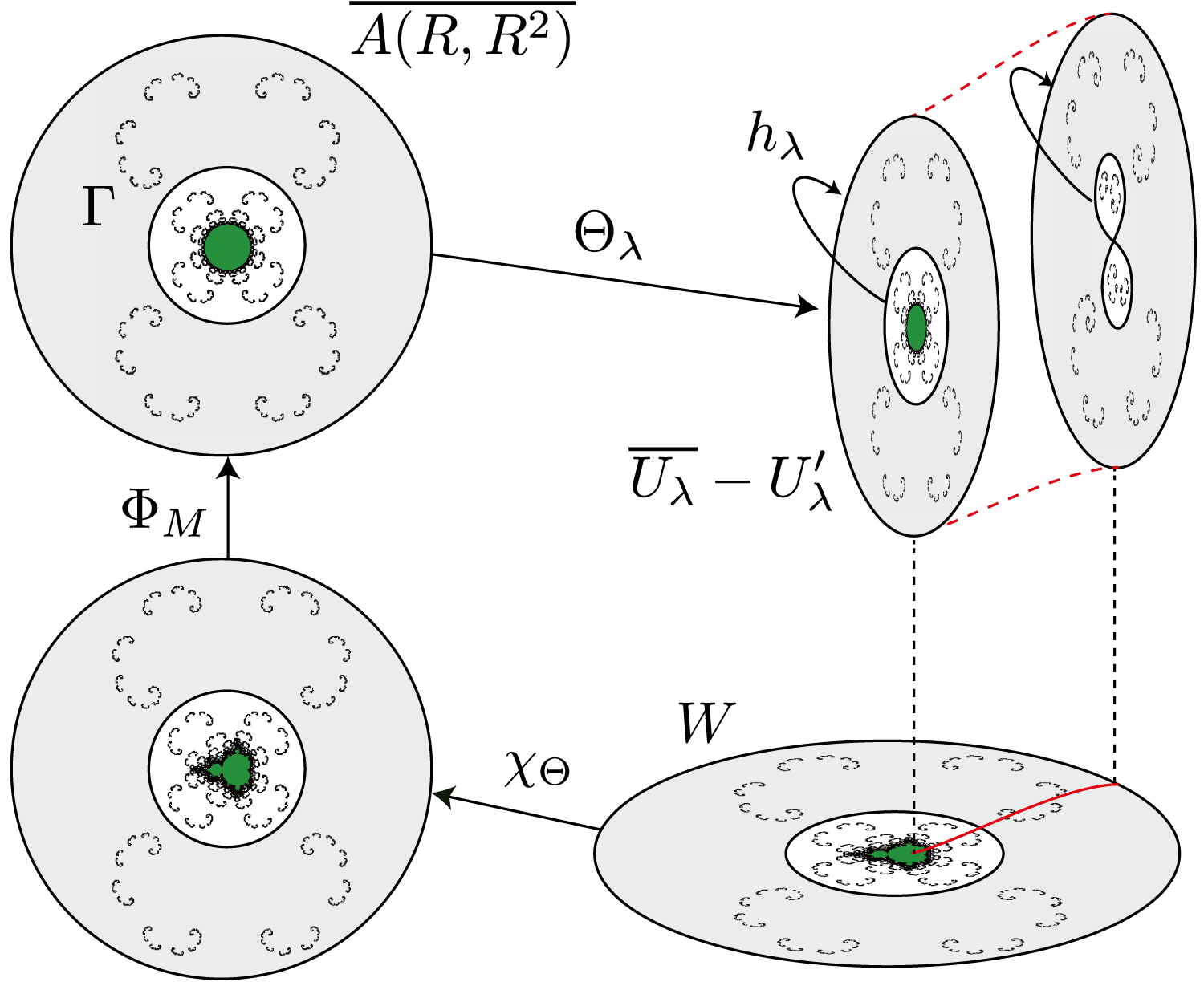}
\caption{\small Tubing 
$\Theta = \{ \Theta_\lambda \}_{\lambda \in W}$.}
\end{figure}

\item[(3)]
Define $C_\lambda := \Theta_\lambda(\partial D(R^2)), \
C'_\lambda := \Theta_\lambda(\partial D(R))$
and let 
$U_\lambda \ (\text{resp.} \ U'_\lambda)$ be the Jordan domain bounded by 
$C_\lambda \ (\text{resp.} \ C'_\lambda)$. Then
$h_\lambda : U'_\lambda \to U_\lambda$ is a quadratic-like map with
a critical point $\omega_\lambda$. Also let
$$
  {\mathcal U} := \{ (\lambda, z) \ | \ \lambda \in W, \ z \in U_\lambda \}, \quad
  {\mathcal U'} := \{ (\lambda, z) \ | \ \lambda \in W, \ z \in U'_\lambda \}
$$
then
$\boldsymbol{h} : {\mathcal U} \to {\mathcal U'}, \ 
(\lambda, z) \mapsto (\lambda, h_\lambda(z))$ is analytic and proper.

\item[(4)]
$\Theta$ is equivariant on the boundary, i.e.,
$\Theta_\lambda(Z^2) = h_\lambda(\Theta_\lambda(Z))$ for $Z \in \partial D(R)$.

\vskip 2mm

\noin
The family of maps 
$\Theta = \{ \Theta_\lambda \}_{\lambda \in W}$ satisfying the above conditions
(1)--(4) is called a {\it tubing}. 

\vskip 2mm

\item[(5)]
$\boldsymbol{h}$ extends continuously to a map 
$\overline{{\mathcal U}'} \to \overline{\mathcal U}$ and 
$\Theta_\lambda : (\lambda ,z) \mapsto (\lambda, \Theta_\lambda(Z))$ extends 
continuously to a map 
$\overline{W} \times \overline{A(R, R^2)} \to \overline{\mathcal U}$
such that $\Theta_\lambda$ is injective on $A(R, R^2)$ for $\lambda \in \partial W$.

\item[(6)]
The map $\lambda \mapsto \omega_\lambda$ extends continuously to $\overline{W}$.

\item[(7)]
$h_\lambda(\omega_\lambda) \in C_\lambda$ for $\lambda \in \partial W$.

\item[(8)]
{\it The one turn condition}: 
When $\lambda$ ranges over $\partial W$ making one turn, then the vector
$h_\lambda(\omega_\lambda) - \omega_\lambda$ makes one turn around 0.
\end{itemize}

Now let $M_{\boldsymbol{h}}$ be the {\it connectedness locus} of the family
$\boldsymbol{h} = \{ h_\lambda \}_{\lambda \in W}$:
$$
M_{\boldsymbol{h}} := \{ \lambda \in W \ | \ K(h_\lambda)  \text{ is connected} \}
= \{ \lambda \in W \ | \ \omega_\lambda \in K(h_\lambda) \}.
$$
Douady and Hubbard (\cite[Chapter IV]{Douady-Hubbard 1985}) showed that there
exists a homeomorphism 
$$
  \chi : M_{\boldsymbol{h}} \to M.
$$ 
This is just a correspondence by the Straightening Theorem, that is, for every
$\lambda \in M_{\boldsymbol{h}}$ there exist a unique
$c = \chi(\lambda) \in M$ such that $h_\lambda$ is hybrid equivalent	
to $P_c(z) = z^2 + c$. Furthermore they
showed that this
$\chi$ can be extended to a homeomorphism $\chi_\Theta : W \to W_M$ by using
$\Theta = \{ \Theta_\lambda \}_{\lambda \in W}$, where
$$
  W_M := \{ c \in \C \ | \ {\mathcal G}_M(c) < 2 \log R \}, 
  \quad {{\mathcal G}_M := }  \text{ the Green function of } M
$$
is a neighborhood of $M$. Also Lyubich showed that $\chi_\Theta$ is 
quasiconformal on any $W'$ with $W' \Subset W$ 
(\cite[p.366, THEOREM 5.5 (The QC Theorem)]{Lyubich 1999}).

Then Douady et al. showed the following:

\begin{prop}\cite[p.29, PROPOSITION 3]{Douady 2000}
For any $\Gamma \subset A(R, R^2)$, 
let $\Gamma_m$ be the preimage of $\Gamma$ by $z \mapsto z^{2^m}$. Then 
$$
  \chi_\Theta^{-1}(\Phi_M^{-1}(\Gamma_m)) 
=\{ \lambda \in W \ | \ h_\lambda^{m+1}(\omega_\lambda) \in \Theta_\lambda(\Gamma) \}
$$
and therefore
$$
  M_{\boldsymbol h}
\cup \{ \lambda \ | \ 
   h_\lambda^k(\omega_\lambda) \in \Theta_\lambda(\Gamma) \ 
   \text{for some} \ k \in \N \}
=
\chi_\Theta^{-1}
\bigg( M \cup \Big( \bigcup_{m=0}^\infty \Phi_M^{-1}(\Gamma_m) \Big) \bigg).
$$
\label{D-BDS Proposition}
\end{prop}

\noin
In what follows, we shall apply this proposition to the rescaled Julia set 
$\Gamma = \Gamma_0(c') (:= J(P_{c'}) \times (\rho/(\rho')^2))$ contained
in $A(R, R^2)$, where $c' \notin M$, $J(P_{c'}) \subset A(\rho', \rho)$ 
and $R := \rho/\rho'$. Then we have
$$
\chi_\Theta^{-1}
\bigg( M \cup \Big( \bigcup_{m=0}^\infty \Phi_M^{-1}(\Gamma_m) \Big) \bigg) 
= 
\chi_\Theta^{-1}
\bigg( M \cup \Phi_M^{-1} \Big( \bigcup_{m=0}^\infty \Gamma_m(c') \Big) \bigg) 
= 
\chi_\Theta^{-1}({\mathcal M}(c')).
$$

\section{Proof of Theorem A for the Misiurewicz case}

Let $M_{s_0}$ be any small Mandelbrot set, where $s_0 \ne 0$ is a 
superattracting parameter (i.e., the critical point $0$ is a periodic 
point of period $p \geq 2$ for $P_{s_0}$).
By the tuning theorem by Douady and Hubbard 
\cite[p.42, Th\'eor\`eme 1 du Modulation]{Haissinsky 2000},
there exists a simply connected domain $\Lambda=\Lambda_{s_0}$ in 
the parameter plane  with the following properties:
\begin{itemize}
\item
If $M_{s_0}$ is a primitive small Mandelbrot set, then 
$M_{s_0}\subset \Lambda$.
\item
If $M_{s_0}$ is a satellite small Mandelbrot set,   
then $M_{s_0}\sminus\{s_0 \perp (1/4)\}\subset \Lambda$.
\item
For any $c \in \Lambda$, $P_c$ is renormalizable with period $p$. 
More precisely, there exist two Jordan domains $\widetilde{U}_c'$ and $\widetilde{U}_c$ with piecewise analytic boundaries such that
$$
  f_c := P_c^p|_{\widetilde{U}_c'} : \widetilde{U}_c' \to \widetilde{U}_c
$$
is a quadratic-like map with a critical point $0 \in \widetilde{U}_c'$. 
In particular, the boundaries of $\widetilde{U}_c'$ and 
$\widetilde{U}_c$ move holomorphically with respect to $c$ over $\Lambda$.
\end{itemize}
In this section let $c_0 \in \partial M$ be any Misiurewicz parameter and $c_1 := s_0 \perp c_0 \in M_{s_0}$. 
The proof of Theorem A for the Misiurewicz case 
breaks into four steps, (M1) to (M4).

\vskip 2mm

\noin
{\bf Step (M1): Definitions of $U_c$, $U_c'$ and $V_c$.} 

\indent
In this step, we shall first construct a 
family $\{f_c:U_c' \to U_c\}_{c\, \in \,S}$ 
of quadratic-like maps over a neighborhood $S$ of $c_1$
by slightly shrinking the original domains of $f_c$,
together with a family of isomorphisms 
$\{g_c:V_c \to U_c\}_{c\, \in \,S}$ 
with a specific property.

We start with the dynamics of $P_{c_1}$ for the parameter $c_1=s_0 \perp c_0$.
Since the quadratic-like map $f_{c_1}= P_{c_1}^p|_{\widetilde{U}_{c_1}'}$ is hybrid equivalent to Misiurewicz $P_{c_0}$,
the ``small" Julia set $J(f_{c_1})$ of $f_{c_1}$ 
is a connected subset of the ``global" Julia set $J(P_{c_1})$ of $P_{c_1}$.
Note that the parameter $c_1$ is also Misiurewicz and thus $J(P_{c_1})$ itself is connected. 
Moreover,
$$
  q_{c_1} := f_{c_1}^l(0) \quad \text{for some} \ l \in \N
$$
is a repelling periodic point of period $k$ with 
multiplier 
$$
  \mu_{c_1} := (f_{c_1}^k)'(q_{c_1}).
$$
Let $\phi_{c_1}: \Omega_{c_1} \to \C$ 
be a linearizing coordinate of $q_{c_1}$
defined on a neighborhood $\Omega_{c_1}$ of $q_{c_1}$
such that $\phi_{c_1}(q_{c_1}) = 0$ and 
$\phi_{c_1}(f_{c_1}^k(z)) = \mu_{c_1} \cdot \phi_{c_1}(z)$.

\begin{lem}
There exist Jordan domains $U_{c_1}$, $U_{c_1}'$ and $V_{c_1}$ with $C^1$
boundaries and integers $N, \ j \in \N$ which satisfy 
the following:
\begin{itemize}
\item[\rm (1)]
$0 \in U_{c_1}' \subset \widetilde{U}_{c_1}'$ and 
$f_{c_1}: U_{c_1}' \to U_{c_1}$ is a 
quadratic-like map.

\item[\rm (2)]
$g_{c_1} := P_{c_1}^N|_{V_{c_1}} : V_{c_1} \to U_{c_1}$ is an isomorphism and
$\overline{f_{c_1}^j(V_{c_1})} 
\subset U_{c_1} \smallsetminus \overline{U_{c_1}'}$.

\item[\rm (3)]
$\overline{V_{c_1}} \subset \Omega_{c_1}$. Also we can take $V_{c_1}$ arbitrarily
close to $q_{c_1} \in \Omega_{c_1}$.
\end{itemize}
\label{def of U_{c_1} etc}
\end{lem}

\begin{figure}[htbp]
\begin{center}
\includegraphics[width=.60\textwidth]{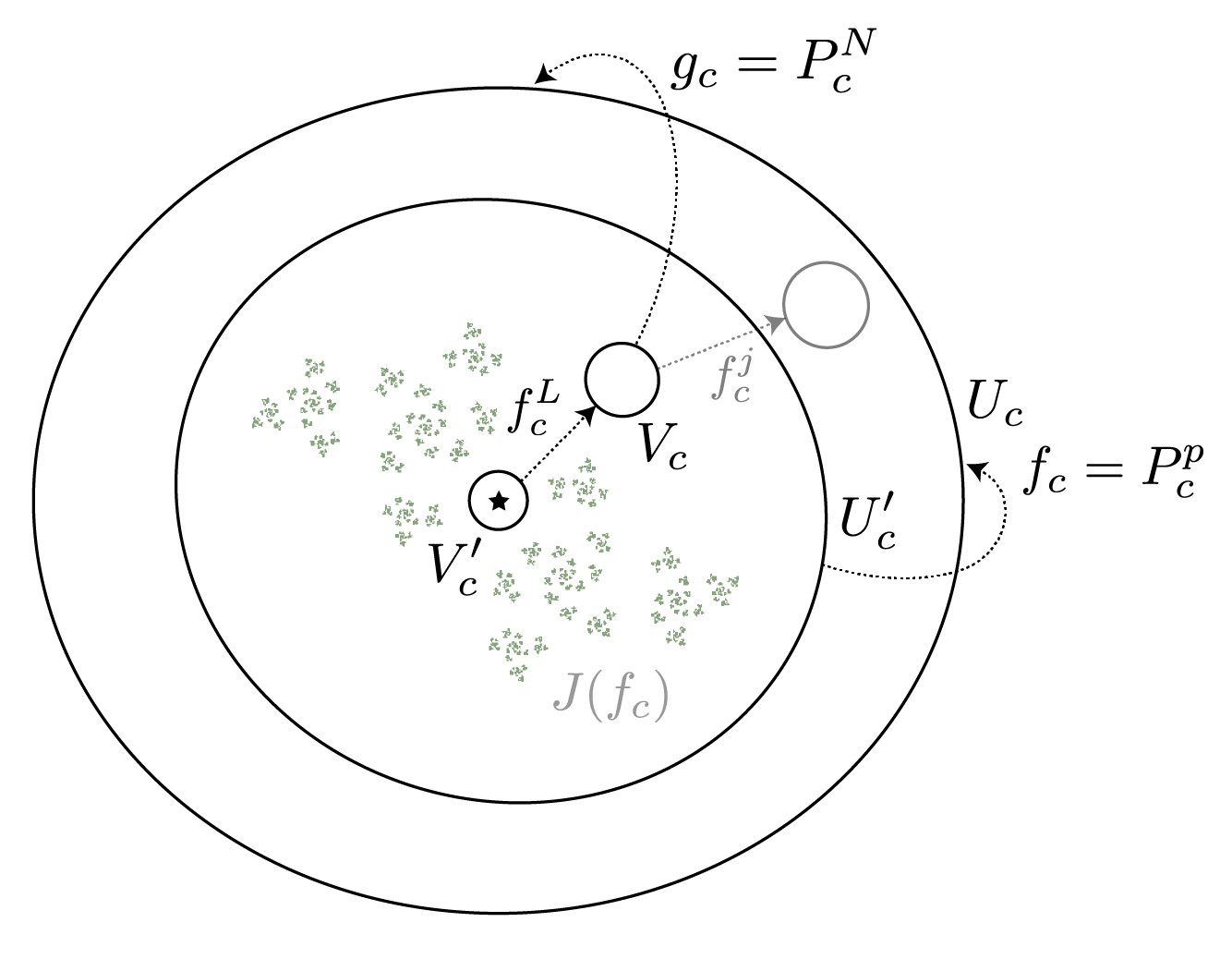}
\caption{\small
The Jordan domains $U_c,\,U_c',\, V_c$, and $V_c'$. 
}\label{UV}
\end{center}
\end{figure}

\paragraph{\bf Proof.}
By shrinking $\widetilde{U}_{c_1}$ and $\widetilde{U}_{c_1}'$ slightly,
we can take Jordan domains $U_{c_1}$ and $U_{c_1}'$ 
with $C^1$ boundaries which are neighborhoods of $J(f_{c_1})$ 
and $f_{c_1} : U_{c_1}' \to U_{c_1}$ is a quadratic-like map. 

For $j \in \N$ let
$$
  A_j := f_{c_1}^{-j}(U_{c_1} \smallsetminus \overline{U_{c_1}'}).
$$
Since $f_{c_1}$ is conjugate to $z^2$ on 
$U_{c_1}' \smallsetminus J(f_{c_1})$,
the annulus $A_j$ is uniformly close to $J(f_{c_1})$ and  
thus $A_j \cap \Omega_{c_1} \ne \emptyset$ for every sufficiently large $j$.
Also since the ``global" Julia set $J(P_{c_1})$ of $P_{c_1}$ is a connected set containing the ``small" Julia set $J(f_{c_1})$, 
the annulus $U_{c_1} \smallsetminus \overline{U_{c_1}'}$
intersects with $J(P_{c_1})$ 
and so does $A_j$.
In particular, for every sufficiently large $j$, 
$A_j \cap \Omega_{c_1}$ contains 
a point $z_0 \in J(P_{c_1})$ 
arbitrarily close to the repelling periodic point 
$q_{c_1} \in J(f_{c_1})$.
Let $B$ be any closed disk in $\Omega_{c_1} \cap A_j$ 
centered at this $z_0$ with an arbitrarily small radius. 

Note that the postcritical set of 
the map $P_{c_1}$ in $\mathbb C$ 
is contained in 
$U_{c_1}' \cup P_{c_1}(U_{c_1}') \cup \cdots \cup P_{c_1}^{p-1}(U_{c_1}')$.
There are two disjoint connected components 
$X:=P_{c_1}^{p-1}(U_{c_1}')$
and 
$-X:=\{-x \in \mathbb{C} ~:~ x \in X\}$
of $P_{c_1}^{-1}(U_{c_1})$, 
where $-X$ does not intersect with neither 
the postcritical set of $P_{c_1}$ nor the critical point $0$.
Hence for any $n \in \N$ and any connected component $V$ of
$P_{c_1}^{-n}(-X)$,
$P_{c_1}^{n+1}:V \to U_{c_1}$ is an isomorphism.

Since the inverse images of $-X$ in the dynamics of $P_{c_1}$
accumulate on any point in the Julia set $J(P_{c_1})$ of $P_{c_1}$ 
(by Montel's theorem), 
the shrinking lemma (\cite[p.86]{Lyubich-Minsky 1997} or 
\cite[Lem.2.9]{Cui-Tan 2018})
implies that we can find a component $V_{c_1}$ of 
$P_{c_1}^{-N+1}(-X)$ contained in the closed disk $B$ 
for some $N \in \N$. 
This gives a desired isomorphism $g_{c_1}:=P_{c_1}^N:V_{c_1} \to U_{c_1}$.

\hfill
\QED ~{\small (Lemma \ref{def of U_{c_1} etc})}
\\

\begin{rem*}
This proof indicates that there are infinitely many different choices of $V_{c_1}$ and one can choose $V_{c_1}$ with arbitrarily small diameter. 
Indeed, each choice of $V_{c_1}$ will give a different \lq\lq decorated small Mandelbrot set".
\end{rem*}

For $ \vep'>0$ given in the statement of Theorem A,
we take a small neighborhood $S$ of $c_1$ contained in 
$D(c_1, \vep') \cap \Lambda$ and 
let $U_c := U_{c_1}$ for each $c \in S$.
By taking a smaller $U_{c_1}$ 
in the previous lemma
and a smaller $S$ if necessary,
we may assume that 
$U_c \,(\equiv U_{c_1})$ is contained in 
$\widetilde{U}_{c}$ 
for each $c \in S$ and 
$U_c':=f_c^{-1}(U_c) \subset \widetilde{U}_{c}'$ 
gives a quadratic-like map $f_c:U_c' \to U_c$
that is a restriction of the original $f_c: \widetilde{U}_{c}' \to \widetilde{U}_{c}$.
Again by taking a smaller $S$ if necessary, 
we may assume in addition that 
for $c \in S$ there exists a component $V_c$
of $P_c^{-N}(U_c)$ which is close to $V_{c_1}$ such that
$$
  \overline{f_c^j(V_c)} \subset U_c \smallsetminus \overline{U_c'}
$$
and $g_c : V_c \to U_c \,(\equiv U_{c_1})$ is an isomorphism. 
See 
Figure \ref{UV}. Let $b_c := g_c^{-1}(0) \in V_c$ and call it
a {\it pre-critical point}. 

\noindent
{\bf Linearizing coordinates.}
Since the fixed point $q_{c_1}$ of $f_{c_1}^k$ is repelling, 
there exists a repelling fixed point $q_c$ of $f_c^k$
that depends holomorphically on $c$ near $c_1$
by the implicit function theorem.\footnote{
While $q_{c_1} = f_{c_1}^l(0)$, 
we have $q_c \ne f_c^l(0)$ for any $c \ne c_1$ 
in a sufficiently small neighborhood of $c_1$. 
Otherwise $q_c = f_c^l(0)$ for any $c$ close to $c_1$ by the identity
theorem.}
To take advantage of the linearizing coordinates
in the next step, 
let us replace the neighborhood $S$ of $c_1$
by an even smaller one such that for each $c \in S$ there exists a unique linearizing coordinate 
$\phi_c : \Omega_c \to \C$ satisfying the following conditions (\cite[\S 8]{Milnor 2006}):
\begin{itemize}
\item
The domain $\Omega_c$ is a neighborhood of $q_c$
and $\phi_c(q_c) = 0$.
\item
Let $\mu(c) := (f_c^k)'(q_c)$. Then
$
  \phi_c(f_c^k(z)) = \mu(c) \cdot \phi_c(z)
$
if both $z$ and $f_c^{k}(z)$ are contained in $\Omega_c$.
\item
(Holomorphic dependence) 
Every compact set $E$ in $\Omega_{c_1}$ 
is contained in $\Omega_c$ for $c \in S$ sufficiently close to $c_1$,
and $\phi_c(z)$ depends holomorphically on $c$ near $c_1$ for each $z \in \Omega_{c_1}$.
\item
$\overline{V_c} \subset \Omega_c$ for any $c \in S$. 
\item
(Normalization) 
$\phi_c(b_c)=1$, where $b_c=g_c^{-1}(0) \in V_c$ is the pre-critical point.
\end{itemize}

\vskip 2mm

\noin
{\bf Step (M2): Construction of the Mandelbrot-like family 
$\boldsymbol{G} = \boldsymbol{G}_n$.} \ 

We shall construct a Mandelbrot-like family 
$$
\boldsymbol{G} = \boldsymbol{G}_n 
:= 
\{ G_c=G_{c, n} : V_c'=V'_{c,n} \to U_c 
   \}_{c \in W=W_n} 
$$
such that $V_c'=V'_{c,n} \subset U_c'$ and $W = W_n \subset S$ for every
sufficiently large $n \in \N$. Note that $G_c$ and $V_c'$ depend on $n$ 
but $U_c'$ and $U_c$ do not. 
Define
$$
  a_c := f_c^l(0)
$$
then as we mentioned above, note that $a_{c_1} = q_{c_1}$ and $a_c \ne q_c$
when $c \ne c_1$ is close to $c_1$.  Hence for such a parameter $c$,  
$a_c$ is repelled from $q_c$ by the dynamics of $f_c^k$. 
By taking a sufficiently small $S$, we may assume that 
$a_c = f_c^l(0) \in \Omega_c$
for each $c \in S$. Now we define
$$
W = W_n 
:= 
\{ c \in S \ | \ f_c^{ki}(a_c) \in \Omega_c \ \text{ for } \ i=1, \dots, n-1
\ \text{ and } \ 
   f_c^{kn}(a_c) (= f_c^{l+kn}(0)) \in V_c \}, 
$$
that is, we consider the parameter $c$ such that the orbit of $a_c$ by $f_c^k$
hits $V_c$.

\begin{lem}
By shrinking $U_c \equiv U_{c_1}$ slightly, the set $W=W_n$ is a non-empty
Jordan domain with $C^1$ boundary for every 
sufficiently large $n$. Moreover there exists an $s_n \in W_n$ such that 
$f_{s_n}^{kn}(a_{s_n}) = b_{s_n}$, which implies
$g_{s_n} \circ f_{s_n}^{l+kn}(0) = P_{s_n}^{(l+kn)p+N}(0) = 0$ and hence
$P_{s_n}$ has a superattracting periodic point. 
\label{W_n is not empty}
\end{lem}

The existence of $W_n$ (and $s_n$, with no explicit description of $\partial W_n$)
 is originally shown by Douady and Hubbard \cite[Chapter V]{Douady-Hubbard 1985},
and independently by Eckmann and Epstein \cite{EE 1985}. (See also McMullen \cite[Theorem 3.1]{McMullen 2000} for a reformulation in a more general context.) 
Here we present a proof that directly show that 
the boundary of $W_n$ inherits regularity (smoothness) from that of $U_c$,
and it is indeed an analytic Jordan curve.

\paragraph{\bf Proof.}
We work with the original $U_c \equiv U_{c_1}$ for the moment and then
shrink (i.e., change the definition of) $U_c \equiv U_{c_1}$ slightly
later to get the result.
In order to do this, 
we observe the dynamics near $q_c$ through the linearizing coordinate 
$\phi_c : \Omega_c \to \C$ of $q_c$. Let
$$
  \tau(c) := \phi_c(a_c), \quad \wt{V}_c := \phi_c(V_c)
$$
then $c \in W_n$ if and only if 
$$
  \mu(c)^n \tau(c) \in \wt{V}_c.
$$
Next recall that
$$
  g_c : V_c \to U_c \equiv U_{c_1}
$$
is an isomorphism by Lemma~\ref{def of U_{c_1} etc} (2) and let
$$
  u : U_c \equiv U_{c_1} \to \D, \quad u(0) = 0
$$
be a Riemann map of $U_c \equiv U_{c_1}$. Then
$$
  u \circ g_c : V_c \to \D, \quad u \circ g_c(b_c)=0
$$
is a Riemann map of $V_c$ and hence
$$
  u \circ g_c \circ \phi_c^{-1} : \wt{V}_c \to \D, \quad 
  u \circ g_c \circ \phi_c^{-1}(1) = 0
$$
is a Riemann map of $\wt{V}_c$. Take the inverse of this map and define
$$
  v(c, \zeta) := \phi_c \circ (u \circ g_c)^{-1}(\zeta), \quad 
  \zeta \in \D \ \text{ with } \ v(c, 0)=1.
$$
Now we solve the equation with respect to the variable $c$
\begin{equation}
  \mu(c)^n \tau(c) = v(c,\zeta)
\label{eqn for v(c,zeta)}
\end{equation}
for each fixed $\zeta \in \D$. 
Since $\mu(c)$ and $\tau(c)$ depend holomorphically on $c$, 
there exist $\alpha \in \N, \ M_0 \ne 0$ and $K_0 \ne 0$ such that
\begin{eqnarray*}
  \mu(c) 
& = &
(f_c^k)'(q_c) = \mu_{c_1} + M_0(c-c_1)^\alpha + O((c-c_1)^{\alpha+1}),  \\
\tau(c) 
& = &
K_0(c-c_1) + O((c-c_1)^2), \quad (c \to c_1).
\end{eqnarray*}
The fact that $K_0 \ne 0$ for the expansion of $\tau(c)$ follows from the 
result by Douady and Hubbard
(\cite[p.333, Lemma 1]{Douady-Hubbard 1985}. See also 
\cite[p.609, Lemma 5.4]{Tan Lei 1990}). 
Now we have
\begin{eqnarray*}
\mu(c)^n \tau(c)
& = &
\big( \mu_{c_1} + M_0(c-c_1)^\alpha + O((c-c_1)^{\alpha+1}) \big)^n 
\cdot \big( K_0(c-c_1) + O((c-c_1)^2) \big) \\
& = & 
\bigg\{ 
  \mu_{c_1}\Big(1 + \frac{M_0}{\mu_{c_1}}(c-c_1)^\alpha + O((c-c_1)^{\alpha+1}) \Big)
\bigg\}^n
\cdot K_0(c-c_1)(1+ O(c-c_1)) \\
& = & 
\mu_{c_1}^nK_0(c-c_1) 
\bigg( 1 + \frac{nM_0}{\mu_{c_1}}(c-c_1)^\alpha + O((c-c_1)^{\alpha+1}) \bigg)
\big( 1+ O(c-c_1) \big) \\
& = & 
\mu_{c_1}^nK_0(c-c_1) + h(c), 
\end{eqnarray*}
where 
$h(c)=O(n\mu_{c_1}^n(c-c_1)^2)$ when $\alpha=1$ and
$h(c)=O(\mu_{c_1}^n(c-c_1)^2)$ when $\alpha \geq 2$. So 
$$
  h(c)=O(n\mu_{c_1}^n(c-c_1)^2)
$$
is true in any case. 
Then the equation (\ref{eqn for v(c,zeta)}) can be rewritten as
\begin{equation}
F(c,\zeta) + G(c,\zeta) = 0,
\label{eqn F+G=0}
\end{equation}
where
$$
  F(c,\zeta) := \mu_{c_1}^nK_0(c-c_1) - v(c_1,\zeta), \quad
  G(c,\zeta) := h(c) - \big( v(c,\zeta)- v(c_1,\zeta) \big).
$$
The equation $F(c,\zeta)=0$ has a unique solution
$$
  c = c_n(\zeta) := c_1 + \frac{v(c_1,\zeta)}{\mu_{c_1}^nK_0}.
$$
Let
$$
  r_n(\zeta) := \bigg| \frac{v(c_1,\zeta)}{\mu_{c_1}^nK_0} \bigg| 
= O(\mu_{c_1}^{-n})
$$
and $\beta := 1/2$. Consider (\ref{eqn F+G=0}) in the disk 
$D(c_n(\zeta), r_n(\zeta)^{1+\beta})$. Since it is easy to see that
$$
  |F(c,\zeta)| = O(r_n(\zeta)^\beta) = O(\mu_{c_1}^{-\beta n}), \quad 
  |G(c,\zeta)| = O(n\mu_{c_1}^{-n})
$$
on the boundary $C := \partial D(c_n(\zeta), r_n(\zeta)^{1+\beta})$ of this disk, 
we have $|F(c,\zeta)| > |G(c,\zeta)|$ on $C$ for sufficiently large $n$. 
By Rouch\'e's theorem (\ref{eqn F+G=0}) has a unique solution 
$c = \check{c}_n(\zeta)$ in $D(c_n(\zeta), r_n(\zeta)^{1+\beta})$, so it 
satisfies
$$
  \check{c}_n(\zeta) 
= c_1 + \frac{v(c_1,\zeta)}{\mu_{c_1}^nK_0}
\big( 1 + O(\mu_{c_1}^{-\beta n}) \big).
$$
By using this solution, we can write
$$
  W_n = \{ \check{c}_n(\zeta) \in \C \ | \ \zeta \in \D \}.
$$

\vskip 2mm

\paragraph{\bf Claim.} \ 
{\it 
{\rm (1)} The map $\check{c}_n : \D \to W_n$ is holomorphic.

\noindent
{\rm (2)} For every $r \in (0,1)$, $\check{c}_n$ is univalent on 
$\overline{\D(r)}$ for every sufficiently large $n$.
}

\vskip 2mm

\paragraph{\bf Proof.}
(1) By the argument principle, for each $\zeta \in \D$ we have
$$
  \check{c}_n(\zeta)
= \frac{1}{2 \pi i}
  \int_C  H(c,\zeta) c \cdot dc, 
$$
where
$$
H(c,\zeta) 
:= 
\frac{\frac{\partial}{\partial c} \big(F(c,\zeta) + G(c,\zeta)\big)}
{F(c,\zeta)+G(c,\zeta)}, \quad
C = \{ z \ | \ |c-c_n(\zeta)|=r_n(t)^{1+\beta} \}
$$
Hence if $|\Delta \zeta| \ll 1$ and $c_n(\zeta+\Delta \zeta) \in \text{int}(C)$,
we have
$$
  \check{c}_n(\zeta+\Delta\zeta)
= \frac{1}{2 \pi i}
  \int_C  H(c,\zeta+\Delta\zeta) c \cdot dc.
$$
Then it follows that $H$ is holomorphic with respect to $\zeta$ and hence
$\check{c}_n(\zeta)$ is holomorphic in a neighborhood of $\zeta$. Thus
$\check{c}_n : \D \to W_n$ is holomorphic.

\noindent
(2) Let $v(\zeta) := v(c_1, \zeta)$ and
$$
  v_n(\zeta) := \mu_{c_1}^nK_0(\check{c}_n(\zeta)-c_1)
= v(c_1,\zeta)\big( 1 + O(\mu_{c_1}^{-\beta n}) \big)
= v(\zeta)\big( 1 + O(\mu_{c_1}^{-\beta n}) \big).
$$
Then for every $r \in (0,1)$, we have
$$
  v_n(\zeta) \to v(c_1,\zeta) = v(\zeta) \quad (n \to \infty) 
$$
uniformly on $\overline{\D(r)}$. In order to show the assertion, it is 
enough to show that $v_n$ is injective on $\overline{\D(r)}$ for every
sufficiently large $n$. Suppose that
$$
  v_n(\zeta_n) = v_n(\zeta_n')
$$
for some $\zeta_n, \ \zeta_n' \in \overline{\D(r)}, \ \zeta_n \ne \zeta_n'$,
where $n$ ranges over a subsequence $\{ n_k \}_{k=1}^\infty$. By taking a
further subsequence, we may assume that
$$
  \zeta_n \to \hat{\zeta}, \quad \zeta_n' \to \hat{\zeta}' \quad \text{ for }
  \ n=n_k, \ k \to \infty.
$$

\noindent
(a) When $\hat{\zeta} = \hat{\zeta}'$ : 
Let $A_0 := v'(\hat{\zeta}) \ne 0$, then there exists a $\delta > 0$
such that
$$
  |v'(\zeta)-A_0| \leq \frac{|A_0|}{4} 
\quad \text{on} \ \D(\hat{\zeta}, \delta).
$$
Since $v_n \to v$ uniformly, we have
$$
  |v_n'(\zeta)-A_0| \leq \frac{|A_0|}{2} 
\quad \text{on} \ \D(\hat{\zeta}, \delta) \ \text{for} \ n \gg 0.
$$
Hence for $\zeta, \ \zeta' \in \D(\hat{\zeta}, \delta)$ we have
\begin{eqnarray*}
\big| \{ v_n(\zeta) - v_n(\zeta') \} - A_0(\zeta-\zeta') \big|
& = &
\Bigg| \int_\zeta^{\zeta'} (v_n'(\zeta) - A_0) d\zeta \Bigg| \\
& \leq &
\frac{|A_0|}{2} |\zeta-\zeta'|
\end{eqnarray*}
It follows that
$$
   \frac{|A_0|}{2} |\zeta-\zeta'|
\leq |v_n(\zeta) - v_n(\zeta')| 
\leq 
\frac 32 |A_0| |\zeta-\zeta'|.
$$
In particular, $v_n$ is injective on $\D(\hat{\zeta}, \delta)$ for
$n \gg 0$. However, $\zeta_n, \ \zeta_n' \in \D(\hat{\zeta}, \delta)$ 
for $n \gg 0$ and this is a contradiction.

\noindent
(b) When $\hat{\zeta} \ne \hat{\zeta}'$ : 
We have
$$
  |v(\hat{\zeta}) - v(\hat{\zeta}')| 
\leq |v(\hat{\zeta}) - v(\zeta_n)| 
+|v(\zeta_n) - v(\zeta_n')| 
+|v(\zeta_n') - v(\hat{\zeta}')|.
$$
As $n=n_k \to \infty$, the first and the third terms of the right hand side 
of this inequality tend to $0$ by the continuity of $v$. Also by
using $v_n(\zeta_n) = v_n(\zeta_n')$,  we have
\begin{eqnarray*}
|v(\zeta_n) - v(\zeta_n')| 
& \leq &
|v(\zeta_n) - v_n(\zeta_n)| + |v_n(\zeta_n') - v(\zeta_n')| \\
& \leq &
2 \sup_{\zeta \in \overline{\D(r)}} |v(\zeta) - v_n(\zeta)| \to 0 \quad
(n=n_k \to \infty), 
\end{eqnarray*}
since $v_n \to v$ uniformly on $\overline{\D(r)}$. This implies 
$v(\hat{\zeta}) = v(\hat{\zeta}')$, but this contradicts the univalence of
$v$.
\QED (Claim)

\medskip

\noindent
By shrinking $U_c \equiv U_{c_1}$ slightly and using the Riemann map $u$ 
of the original $U_c \equiv U_{c_1}$, the boundary of the new 
$U_c \equiv U_{c_1}$ is parametrized as $u^{-1}(\gamma(t))$, where 
$\gamma(t) = re^{2\pi it} \subset \D \ (t \in [0,1])$ and $r \in (0,1)$ 
is close to $1$. Then 
$\partial \wt{V}_c$ is parameterized as $v(c, \gamma(t))$ and hence 
$\partial W_n$ (for the new $W_n$) is parameterized as 
$\check{c}_n(\gamma(t))$ by using the solution $\check{c}_n(\zeta)$ for 
the equation (\ref{eqn for v(c,zeta)}). Clearly this is a $C^1$ Jordan curve
and $W_n$ is the image of $\D(r)$ by $\check{c}_n(\zeta)$. This shows that 
$W_n$ is a non-empty Jordan domain with $C^1$ boundary. 
In particular, let $s_n := \check{c}_n(0)$ then this satisfies
$\mu(s_n)^n\tau(s_n) = 1$. This means that
$f_{s_n}^{kn}(a_{s_n}) = b_{s_n}$,
which implies $g_{s_n} \circ f_{s_n}^{l+kn}(0) = P_{s_n}^{(l+kn)p+N}(0)= 0$. 
Hence $P_{s_n}$ has a superattracting periodic point. 
This completes the proof of Lemma \ref{W_n is not empty}.
\QED ~{\small (Lemma \ref{W_n is not empty})}
\\

We call 
$s_n \in W_n$ the {\it center} of $W_n$. Now let $L = L_n := l+kn$ and
$V_c'=V'_{c,n}$ be the Jordan domain bounded by the component of $f_c^{-L}(V_c)$ 
containing $0$ and define
$$
  G_c = G_{c, n} := g_c \circ f_c^L : V_c' \to U_c
  \ \  \text{and} \ \ 
 \boldsymbol{G} = \boldsymbol{G}_n 
 := 
\{ G_c \}_{c \in W_n}, 
$$
where $W_n = \{ c \in S \ | \ f_c^L(0) \in V_c \}$.
See Figure \ref{UV}.

\medskip

\noin
{\bf Step (M3): Proof for $\boldsymbol{G} = \boldsymbol{G}_n$ 
being a Mandelbrot-like family.}

The map
$f_c^L : V_c' \to V_c$ is a branched covering of degree 2 and
$g_c : V_c \to U_c$ is a holomorphic isomorphism. Hence 
$G_c := g_c \circ f_c^L : V_c' \to U_c $ is a quadratic-like map.

Next we construct a tubing 
$\Theta = \Theta_n = \{ \Theta_c \}_{c \in W_n}$ for $\boldsymbol{G}_n$
as follows: For $s_n \in W_n$,  since $f_{s_n}^L(0) \in V_{s_n}$ and 
$f_{s_n}^j(V_{s_n}) \subset U_{s_n} \smallsetminus \overline{U_{s_n}'}$, 
from Lemma \ref{def of U_{c_1} etc}, we have
$f_{s_n}^{L+j}(0) \notin U_{s_n}'$. 
It follows that $J(f_{s_n})$ is a Cantor set,
 which is quasiconformally homeomorphic to a quadratic Cantor
Julia set $J(P_{c_0+\eta_n})$ for some $\eta=\eta_n$ with 
$c_0+\eta_n \notin M$ by the Straightening Theorem.
By continuity of the straightening of $f_c$ for 
$c \in \Lambda$, we have $|\eta_n| < \vep$ for sufficiently large $n$.
Let $\Psi_{s_n}$ be the quasiconformal straightening map
that conjugates $f_{s_n}$ and $P_{c_0+\eta_n}$ 
defined on a neighborhood of $K(f_{s_n})$.
Then the image of $J(f_{s_n})$ by $\Psi_{s_n}$ is $J(P_{c_0+\eta_n})$.
Take an $R > 1$ and let $\rho' := R^{-1/2}$ and $\rho := R^{1/2}$ such
that $J(P_{c_0+\eta_n}) \subset A(\rho', \rho)$. 
Define the rescaled Julia set
$$
  \Gamma := \Gamma_0(c_0+\eta_n) 
= \Gamma_0(c_0+\eta_n)_{\rho', \rho}
:= J(P_{c_0+\eta_n}) \times R^{3/2}
\subset A(R, R^2). 
$$

\begin{lem}\label{lem_Theta_n^0}
There exists a quasiconformal homeomorphism
$$
\Theta_n^0 : \overline{A(R, R^2)} \to 
              \overline{U}_{s_n} \smallsetminus V_{s_n}'
$$ 
for $s_n$ such that 
\begin{itemize}
\setlength{\itemsep}{0.5mm} 
\item
$\Theta_n^0$ is quasiconformal,

\item
$\Theta_n^0$ is equivariant on the boundary, i.e.,
$\Theta_n^0(Z^2) = G_{s_n}(\Theta_n^0(Z))$ for $|Z| = R$, 

\item
$\Theta_n^0(Z) = \Psi_{s_n}^{-1}(R^{- 3/2}Z)$ 
for $Z \in \Gamma_0(c_0+\eta_n)$,

\item
$\Theta_n^0(\Gamma_0(c_0+\eta_n)) = J(f_{s_n})$. 
\end{itemize}
\end{lem}
Note that we have no dilatation control of the quasiconformal map $\Theta_n^0$ for this lemma, which may depend on the shapes of the image $\overline{U}_{s_n} \smallsetminus V_{s_n}'$ and $J(f_{s_n})$. 
However, in the proof of Theorem C (Section 8, Claim 6), 
we will show that such a $\Theta_n^0$ can be almost conformal 
by choosing an appropriate $s_n$.

\paragraph{\bf Proof.}
Since the boundary components of 
the closed annuli $\overline{A(R, R^2)}$ and
$\overline{U}_{s_n} \smallsetminus V_{s_n}'$ are smooth,
we can take a smooth homeomorphism $\psi_0$ between 
$\partial D(R^2)$ and $\partial {U}_{s_n}$.
By letting $\psi_0(Z)$ be an appropriate branch of 
$G_{s_n}^{-1}(\psi_0(Z^2))$ for $Z \in \partial D(R)$, 
we have a smooth, equivariant homeomorphism $\psi_0$ between the
boundaries of the closed annuli.

Next we consider $\Gamma=J(P_{c_0+\eta_n}) \times R^{3/2}$.
Recall that the quasiconformal (straightening) map 
$\Psi_{s_n}$ 
sends a neighborhood of $J(f_{s_n})$
to that of $J(P_{c_0+\eta_n})$.
There exists a neighborhood $D^\ast$ of $\Gamma$
such that the quasiconformal map 
$\psi_1:=\Psi_{s_n}^{-1}(R^{-3/2} Z)$
defined on $D^\ast$ 
sends $\Gamma$ to $J(f_{s_n})$.
Since $\Gamma$ is a Cantor set,
we may choose $D^\ast$ such that $D^\ast$ 
is a finite union of smooth Jordan domains
satisfying $D^\ast \Subset A(R,R^2)$ 
and 
$\psi_1(D^\ast) \Subset 
U_{s_n} \smallsetminus\overline{V_{s_n}'}$. 
Now the sets 
$A(R,R^2)  \smallsetminus \overline{D^\ast}$ 
and 
${U}_{s_n} \smallsetminus \overline{V_{s_n}' \cup \psi_1(D^\ast)}$
are multiply connected domains with the same connectivity.
By a standard argument in complex analysis
(see \cite[Chapter 6, Theorem 10]{Ahlfors 1978} for example),
they are conformally equivalent to round annuli with 
concentric circular slits, and there is a 
quasiconformal homeomorphism $\psi_2$ between these domains.
Since each component of $\psi_1(D^\ast)$ is a quasidisk,
we can modify $\psi_2$ such that the boundary correspondence 
agrees with $\psi_0$ and $\psi_1$. 
Hence we obtain a desired quasiconformal homeomorphism 
$\Theta_n^0$
by gluing $\psi_0$, $\psi_1$, and this modified $\psi_2$. 
\QED{\small (Lemma \ref{lem_Theta_n^0})}\\

The Julia set
$J(f_{c}) \subset U_c' \smallsetminus \overline{V_c'}$ is a Cantor set
for every $c \in W_n$ for the same reason for $J(f_{s_n})$ and this, 
as well as $\partial U_c$ and $\partial V_c'$ undergo holomorphic motion 
(see \cite[p.229]{Shishikura 1998}).
By S{\l}odkowski's theorem (\cite{Slodkowski 1991}) there exists a 
holomorphic motion $\iota_c$ on $\C$
which induces these motions. Finally 
define $\Theta_c := \iota_c \circ \Theta_n^0$ , then 
$\Theta = \Theta_n := \{ \Theta_c \}_{c \in W_n}$ is a
tubing for $\boldsymbol{G}_n$.

Now we have to check that
$\boldsymbol{G}_n$ with $\Theta_n$
satisfies the conditions (1)--(8) for a Mandelbrot-like family.
The condition (1) is already shown in Lemma~\ref{W_n is not empty}. 
It is easy to check the conditions (2)--(7). 
Finally the one turn condition (8) is proved as follows: 
Note that $\check{c}_n(\gamma(t))$ satisfies
$$
\mu (\check{c}_n(\gamma(t)))^n 
\tau(\check{c}_n(\gamma(t)))
= v(\check{c}_n(\gamma(t)),\gamma(t)).
$$
When $c$ ranges over $\partial W_n$ making one turn, the variable $t$ for
both sides varies from $t=0$ to $t=1$. 
Since $v(\check{c}_n(\gamma(t)),\gamma(t))$ 
is very close to $v(c_1, \gamma(t))$, which is a parameterization of
$\partial \wt{V}_{c_1}$ for sufficiently large $n$, 
$v(\check{c}_n(\gamma(t)),\gamma(t))$ and hence
$\mu(\check{c}_n(\gamma(t)))^n \tau(\check{c}_n(\gamma(t)))$ 
makes one turn in a very thin tubular neighborhood of 
$\partial \wt{V}_{c_1}$ as $t$ moves from 0 to 1. 
This implies that $f_c^{kn}(a_c) = f_c^{l+kn}(0) = f_c^L(0)$ 
makes one turn in a very thin tubular neighborhood of $\partial V_{c_1}$. 
Hence $G_c(0)-0 = G(0) = g_c \circ f_c^{l+kn}(0) = g_c \circ f_c^L(0)$ makes
one turn in a very thin tubular neighborhood of $\partial U_{c_1}$. 
In particular this shows that $G_c(0)-0$ makes one turn around $0 \in U_{c_1}$.

\vskip 2mm

\noin
{\bf Step (M4): End of the proof of Theorem A for the Misiurewicz case.}

For every $\vep > 0$ and $\vep' > 0$, take a sufficiently large $n \in \N$ such
that $c_0 + \eta = c_0 + \eta_n \in D(c_0, \vep) \smallsetminus M$. We conclude that
the model ${\mathcal M}(c_0+\eta)$ appears quasiconformally in $M$ in the neighborhood
$D(c_1, \vep') = D(s_0 \perp c_0, \vep')$ of $c_1 = s_0 \perp c_0$ by applying 
Proposition \ref{D-BDS Proposition} to the Mandelbrot-like family 
$\boldsymbol{G}=\boldsymbol{G_n}$ with $\Theta = \Theta_n$. Indeed from
Proposition \ref{D-BDS Proposition}, the set
$$
{\mathcal N}
:=
M_{\boldsymbol G}
\cup 
\{ c \ | \ 
     G_c^k(0) \in \Theta_c(\Gamma_0(c_0+\eta)) \quad \text{for some} \ k \in \N \}
$$
is the image of ${\mathcal M}(c_0+\eta)$ by the quasiconformal map 
$\chi_{\Theta}^{-1} = \chi_{\Theta_n}^{-1}$, where $M_{\boldsymbol G}$ is the
connectedness locus of ${\boldsymbol G}$. On the other hand, for 
$c \in M_{\boldsymbol G}$, the orbit of the critical point $0$ by 
$G_c = g_c \circ f_c^{l+kn} = P_c^{Lp+N}$ is bounded, which
implies that the orbit of $0$ by $P_c$ is also bounded and hence $c \in M$.
If $G_c^k(0) \in \Theta_c(\Gamma_0(c_0+\eta))$ for some $k \in \N$, then
$c \in M$ as well. So the set $\mathcal N$ is a subset of $M$.
In particular, since a conformal image $\Phi_M^{-1}(J(P_{c_0+\eta}) \times R^{3/2})$
of $J(P_{c_0+\eta})$ is a subset of ${\mathcal M}(c_0+\eta)$, we conclude that 
$J(P_{c_0+\eta})$ appears 
quasiconformally in $M$. This completes the proof of Theorem A for the 
Misiurewicz case.
\QED 
\medskip

\begin{rem*}
In \cite{Douady 2000}, there is no proof for $\partial W_n$ being a $C^1$ 
Jordan curve and also the proof for the one turn condition (8) is intuitive. 
\end{rem*}

\section{Proof of Theorem A for the parabolic case}
As in the previous section, 
let $M_{s_0}$ be the small Mandelbrot set with center $s_0 \neq 0$
such that $0$ is a periodic point of period $p \ge 2$,
and let $\Lambda=\Lambda_{s_0}$ be the simply connected domain 
where the family 
$
\{
f_c := P_c^p|_{\widetilde{U}_c'} : \widetilde{U}_c' \to \widetilde{U}_c
\}_{c \, \in \, \Lambda}
$
of quadratic-like maps is defined.
In this section let $c_0 \in \partial M$ be any parabolic parameter 
and $c_1:=s_0 \perp c_0 \in M_{s_0}$.

A simple way to show Theorem A for the parabolic case is the following:
since the Misiurewicz parameters are dense in the boundary of the Mandelbrot set,
we can find a Misiurewicz parameter $c_0'$ that is arbitrarily close to the parabolic parameter $c_0$. By continuity of the tuning map 
$c \mapsto s_0 \perp c$, we may apply Theorem A for the Misiurewicz case.

There is another proof that is 
independent of the Misiurewicz case,
based on Douady's original proof for the cauliflower. 
(Hence by the same logic the parabolic case implies the Misiurewicz case.)
Details will be given in a forthcoming paper 
\cite{Kawahira-Kisaka 2023}.

\begin{rem*}
Theorem A is a kind of generalization of the Douady's result but
the statements of the results of ours and his are not quite parallel.
Actually Douady considered not only the case of the quadratic family
but also more general situation and proved a theorem 
(\cite{Douady 2000}, p.23, THEOREM 2) and then showed the theorem for 
the Mandelbrot set (\cite{Douady 2000}, p.22, THEOREM 1) by using it. 
Douady's result also shows that a sequence of quasiconformal images of 
${\mathcal M}(1/4 + \vep)$ appears in 
$D(s_0 \perp (1/4), \vep')$. It is possible to state our result like 
Douady's. But in order to do this, it is necessary to assume several 
conditions which are almost obvious for the quadratic family case and 
this would make the argument more complicated. So we just concentrated 
on the case of the quadratic family. We avoided stating our result like 
\lq\lq a sequence of quasiconformal images of 
${\mathcal M}(c_0 + \eta)$ appears" for the same reason.

In what follows, we summarize the general situation under which a result similar
to THEOREM 2 in \cite{Douady 2000} (that is, Theorem A'' below) hold and this 
implies our Theorem A. These are the essential assumptions for more
general and abstract settings, which leads to the general result Theorem A''.

\medskip

\noin
$\bullet$ \
$\{ f_c : U_c' \to U_c \}_{c \,\in\, \Lambda}$ is an analytic family of
quadratic-like maps with a critical point $\omega_c$, where $\Lambda \subset \C$
is an open set. The parameter $c_1 \in \Lambda$ is either Misiurewicz or 
parabolic.

\noin
$\bullet$ \
$\{ g_c : V_c \to U_c \}_{c \,\in\, \Lambda}$ is an analytic family of analytic
isomorphism, where $V_c$ satisfies 
$\overline{f_c^j(V_c)} \subset U_c \smallsetminus \overline{U_c'}$ 
for some $j \in \N$.

\noin
$\bullet$ \
The open sets $U_c, \ U_c'$ and $V_c$ are Jordan domains with $C^1$ boundary
and move by a holomorphic motion. Let $z_c(t)$ be a parametrization
of $\partial V_c$. Then $z_c(t)$ is holomorphic in $c$ and $C^1$ in $t$ and
$\frac{\partial^2}{\partial c \partial t}z_c(t)$ exists and continuous.

\noin
$\bullet$ \
(1) \ When $c_1$ is Misiurewicz, for some $l \in \N$, 
$f_{c_1}^l(\omega_{c_1})$ is a repelling periodic point of period $k$ 
and we let $a_c := f_c^l(\omega_c)$.
Let $q_c$ be the repelling periodic point persisting when $c$ is perturbed 
from $c_1$. Then assume that $a_c \ne q_c$ for $c (\ne c_1)$ which is 
sufficiently close to $c_1$.

(2) \ 
When $c_1$ is parabolic, $f_{c_1}$ has a parabolic periodic
point $q_{c_1}$ of period $k$ with multiplier
$$
\mu_{c_1} :=
(f_{c_1}^k)'(q_{c_1})
=e^{2 \pi i \nu'/\nu},
$$
where $\nu'$ and $\nu$ are coprime integers.
Then assume the following:
The parabolic fixed point $q_{c_1}$ of $f_{c_1}^k$
splits into one fixed point $q_c$ and
a cycle of period $\nu$ of $f_c^k$ for each $c \neq c_1$.
There exists a suitable sector $S$ in the parameter space and
holomorphic local coordinate $w=\psi_c(z) \ (c \in S)$ near
$q_{c}$ with $\psi_c(q_{c})=0$ such that
$$
\psi_c \circ f_{c}^{k\nu} \circ \psi_c^{-1}(w)
=\mu_c^\nu w\,(1 + w^{\nu}+O(w^{2\nu})),
$$
where $\mu_c^\nu \to 1$ and $\psi_c \to \psi$
uniformly as $c \in S$ tends to $c_1$.
(Douady gives a sufficient
condition for this condition when $k = 1$ and  $\mu_{c_1}= 1$
in \cite[p.23]{Douady 2000}.)

\medskip

\noin
Note that there exists a $c_0$ such that $f_{c_1}$ is hybrid equivalent to
$P_{c_0}$. Now define the map $F_c : U_c' \cup V_c \to U_c$ so that
$F_c := f_c$ on $U_c'$ and $F_c := g_c$ on $V_c$. Also define
\begin{eqnarray*}
& & K(F_c) := \{ z \ | \ F_c^n(z) \ \text{is defined for all} \ n \ 
              \text{and} \ F_c^n(z) \in U_c' \cup V_c \}, \\
& & M_F := \{ c \in \Lambda \ | \ \omega_c \in K(F_c) \}.
\end{eqnarray*}

\noin
Under the above assumptions, we can show the following theorem which
implies our Theorem A:

\medskip

\noindent
{\bf Theorem A''}
{\it
For every small $\vep > 0$ and $\vep' > 0$, there exists an 
$\eta \in \C$ with $|\eta| < \vep$ and
$c_0+\eta \notin M$ such that the
decorated Mandelbrot set ${\mathcal M}(c_0+\eta)$ appears quasiconformally
in $M_F$. 
}

\end{rem*}

\setcounter{section}{5}
\section{Proof of Theorem B}
Let $M_{s_1}$ be the main Mandelbrot set of the quasiconformal copy
of the decorated Mandelbrot set $\cM(c_0+\eta)$ given in Theorem A. 
Choose any $c \in M$ and set $\sig:=s_1 \perp c \in M_{s_1}$. 
(For example, 
let $c_0$ be the Misiurewicz parameter for which $P_{c_0}^5(0)=P_{c_0}^4(0)$ and $c$ the parameter for Douady's rabbit
as in Figure \ref{nested structure for a filled Julia set}.)

Let $\Phi_{c}:\C \sminus K(P_{c}) \to \C  \sminus \Bar{\D}$ be 
the B\"ottcher coordinate for $P_c$. 
For any $R>1$ with 
$J(P_{c_0 + \eta}) \subset A(R^{-1/2}, R^{1/2})$,
we take the Jordan domains 
$\Omega_1'$ and $\Omega_1$ in $\C$ with $\Omega_1' \Subset \Omega_1$
whose boundaries are the inner and the outer boundaries of
$\Phi_{c}^{-1}(A(R, R^{2}))$.
(That is, we take Douady's radii $\rho'=R^{-1/2}$ and $\rho=R^{1/2}$
in the definition of rescaled Julia set. See section 2.) 
Then $P_{c}:\Omega_1' \to \Omega_1$ is a 
quadratic-like restriction of $P_{c}$, 
and the decorated filled Julia set $\cK_{c}(c_0+\eta)=\cK_{c}(c_0+\eta)_{R^{-1/2},R^{1/2}}$ 
is a compact set in $\Omega_1$.

Now we want to show that for $\sig=s_1 \perp c \in M_{s_1}$
the filled Julia set $K(P_{\sig})$ contains 
a quasiconformal copy $\cK$ of the model set $\cK_{c}(c_0+\eta)$.
Consider the quadratic-like maps
$f_{\sig}: U_{\sig}' \to U_{\sig}$ and $G_{\sig}: V_{\sig}' \to U_{\sig}$ 
given in the proof of Theorem A.
Since we have $J(f_{\sig}) \subset U_{\sig}  \sminus  \Bar{V_{\sig}'}$,
the filled Julia set $K(G_{\sig}) \subset V_{\sig}'$
is surrounded by the set
$$
{\Gamma}:= \bigcup_{m \ge 0} G_{\sig}^{-m}(J(f_{\sig})).
$$
Let $\cK$ be the union $K(G_{\sig})  \cup \Gamma$, 
which is a compact subset of $U_{\sig}$.
Then the boundary $\partial \cK$ is contained in 
$\partial K(P_{\sig})$,
since the set of points that eventually lands on 
a repelling cycle of $f_\sigma$ or $G_\sigma$
is dense in $\partial \cK$.
Hence it is enough to show that there exists a 
quasiconformal map on a domain 
that maps the model set $\cK_{c}(c_0+\eta)$ to $\cK$.

Let $h=h_{\sig}:U_{\sig} \to \C$ be 
a straightening map of $G_{\sig}:V_{\sig}' \to U_{\sig}$.
By setting $\Omega_2:=h(U_{\sig})$ and $\Omega_2':=h(V_{{\sig}}')$,
the map $P_{c}=h \cc G_{\sig} \cc h^{-1}:\Omega_2' \to \Omega_2$
is also a quadratic-like restriction of $P_{c}$
such that 
$h(K(G_{\sig}))=K(P_{c})$ and 
$h(J(f_\sigma)) \subset h(U_\sigma \smallsetminus \overline{V_\sigma'})
=\Omega_2\smallsetminus\overline{\Omega_2'}$.
By slightly shrinking $\Omega_2$, we may assume that 
the boundaries of $\Omega_2$ and $\Omega_2'$ are smooth Jordan curves.
Since $h$ is quasiconformal, 
it suffices to show that there exists a quasiconformal map 
$H:\Omega_1 \to \Omega_2$  
that maps the model set $\cK_{c}(c_0+\eta)$ 
onto $h(\cK)$.

Now we claim:
\begin{lem}
\label{lem_B-0}
There exists a quasiconformal homeomorphism 
$H:\Bar{\Omega_1} \sminus {\Omega_1'} \to \Bar{\Omega_2} \sminus {\Omega_2'}$
such that
\begin{itemize}
\item
$H$ is equivariant. That is, 
$P_{c}(H(z))=H(P_{c}(z))$ for any  $z \in \partial \Omega_1'$.
\item
$H$ maps 
$J^\ast:=\Phi_{c}^{-1}(J(P_{c_0 + \eta}) \times R^{3/2})$ in the model set onto $h(J(f_{\sig}))$.
\end{itemize}
\end{lem}

\paragraph{\bf Proof.}
Since the boundary components of these annuli are smooth,
we can take a smooth homeomorphism between $\partial \Omega_1$
and $\partial \Omega_2$.
By pulling it back by the action of $P_{c}$,
we have a smooth, equivariant homeomorphism $\psi_0$ between the
boundaries of the closed annuli $\Bar{\Omega_1} \sminus {\Omega_1'}$
and $\Bar{\Omega_2} \sminus {\Omega_2'}$.

Next we consider $J^\ast$.
Recall that $\sigma \in M_{s_1} \subset W \subset \Lambda$,
where $W$ and $\Lambda$ are given in the proof of Theorem A.
Hence there exists a straightening map $\hat{h}:U_\sigma \to \C$ 
that quasiconformally conjugates 
the quadratic-like map $f_{\sig}:U_\sigma' \to U_\sigma$ 
to $P_{\sigma'}:\hat{h}(U_\sigma') \to \hat{h}(U_\sigma)$ for some 
$\sigma' \in \C \smallsetminus M$.
Consider a sequence of homeomorphisms
$$
J^\ast=\Phi_{c}^{-1}(J(P_{c_0 + \eta}) \times R^{3/2})
\stackrel{(1)}{\longrightarrow}
J(P_{c_0+\eta})
\stackrel{(2)}{\longrightarrow}
J(P_{\sig'})
\stackrel{(3)}{\longrightarrow}
J(f_{\sig}) 
\stackrel{(4)}{\longrightarrow}
h(J(f_{\sig})), 
$$
where (1) -- (4) are given as follows: 
\begin{enumerate}
\item 
This is just a conformal map $z \mapsto R^{-3/2} \Phi_c(z)$
restricted to $J^\ast$.
\item
Take a simply connected domain $W'$ in $\C \sminus M$ 
containing $c_0+\eta$ and $\sigma'$. 
Then there exists a holomorphic motion of $J(P_{c_0+\eta})$
over $W'$ that gives a quasiconformal map
on the plane that sends $J(P_{c_0+\eta})$ to $J(P_{\sig'})$
by the Bers-Royden theorem (\cite[Thoerem 1]{Bers-Royden 1986}).
\item
This is $(\hat{h}|_{J(f_\sigma)})^{-1}$, which is a 
restriction of a quasiconformal map 
$\hat{h}^{-1}:\hat{h}(U_\sigma) \to U_\sigma$.
\item
This is $h|_{J(f_\sigma)}$,
which is a restriction of the quasiconformal straightening map 
${h}:U_\sigma \to h(U_\sigma)$ of the quadratic-like map 
$G_\sigma:V_\sig' \to U_\sig$.
\end{enumerate}
Hence there exists a neighborhood $D^\ast$ of $J^\ast$
and a quasiconformal map $\psi_1:D^\ast \to \C$ 
that sends $J^\ast$ to $h(J(f_{\sig}))$.

The remaining construction of the map is the same as 
Lemma \ref{lem_Theta_n^0}.
We define a quasiconformal map $\psi_2$ 
between 
$\Omega_1  \sminus\overline{\Omega_1' \cup D^\ast}$ 
and $\Omega_2  \sminus \overline{\Omega_2'\cup \psi_1(D^\ast)}$.
Then we obtain a desired quasiconformal homeomorphism $H$
by gluing $\psi_0$, $\psi_1$, and a modified $\psi_2$. 
\QED {\small (Lemma 6.1)}

\medskip
By pulling back the map $H$ given in Lemma \ref{lem_B-0}
by the dynamics of $P_{c}$,
we have a unique homeomorphic extension
$H:\Bar{\Omega_1} \sminus K(P_{c}) \to \Bar{\Omega_2} \sminus K(P_{c})$
such that $P_{c}(H(z))=H(P_{c}(z))$ for any $\Omega_1' \sminus K(P_{c})$
and that $H$ maps the decoration of $\cK_{c}(c_0+\eta)$ 
onto $h(\Gamma)$.

We employ the following lemmas. 
(For the proofs, see Lyubich's book in preparation \cite{Lyubich Book}.\footnote{It is still being updated. 
The numbers of subsections below are tentative.})

\begin{lem}[{\cite[\S 41.3]{Lyubich Book}}] 
\label{lem_B_1}
Let $f:U' \to U$ be a quadratic-like map with connected Julia set.
Let $W_1 \subset U$ and $W_2 \subset U$ be two open annuli whose
inner boundary is $J(f)$. 
Let $H:W_1 \to W_2$ be an automorphism of $f$,
that is, $f(H(z))=H(f(z))$ on $f^{-1}(W_1)$.
Then $H$ admits a continuous extension to
a map $H:W_1 \cup J(f) \to W_2 \cup J(f)$
identical on the Julia set. 
\end{lem}

\begin{lem}[\bf Bers' Gluing Lemma, {\cite[\S 13.3]{Lyubich Book}}]
\label{lem_Bers}
Let $K$ be a compact set in $\C$ and 
let $\Omega_1$ and $\Omega_2$ 
be neighborhoods of $K$ 
such that there exists two quasiconformal maps
$H_1:\Omega_1  \sminus K \to \C$ and $H_2: \Omega_2 \to \C$
that match on $\partial K$, i.e., the map
$H:\Omega_1 \to \C$ defined by $H(z):=H_1(z)$
for $z \in \Omega_1  \sminus K$ 
and $H(z):=H_2(z)$ for $z \in K$
is continuous. 
Then $H$ is quasiconformal and $\mu_H=\mu_{H_2}$
for almost every $z \in K$.
\end{lem}

Now we apply Lemma \ref{lem_B_1} by regarding $P_{c}$ and $\Omega_j \sminus K(P_{c})$
as $f$ and $W_j$ for each $j=1,\,2$.
It follows that the restriction $H_1:=H|_{\Omega_1 \sminus K(P_{c})}$
of the map 
$H:\Bar{\Omega_1} \sminus K(P_{c}) \to \Bar{\Omega_2} \sminus K(P_{c})$
admits a continuous extension to 
$H_1:\Omega_1 \sminus \mathrm{int}(K(P_{c})) \to \Omega_2 \sminus \mathrm{int}(K(P_{c}))$
that agrees with the identity map $H_2:=\mathrm{id}:\Omega_2 \to \Omega_2$ on $\partial K(P_c)$.
By Lemma \ref{lem_Bers},
we have a quasiconformal map $H:\Omega_1 \to \Omega_2$ 
such that $H(\cK_{c}(c_0+\eta))=h(\cK)$.
 \QED

\section{Almost conformal straightenings}
In this section we establish a general formulation
of quadratic-like families 
that generate ``fine" copies of the Mandelbrot set. 

The notation here (for example, the way we use $f_c, U_c', U_c, \ldots$) is different from that in the other sections. 

\begin{lem}[\bf Almost conformal straightening]
\label{lem_almost_conformal}
Fix two positive constants $r$ and $\delta$.
Suppose that for some 
$R>\max\{2(r+\delta), 6\}$ we have a function $u=u_c(z)=u(z,c)$ 
that is holomorphic with $u_c'(0)=0$ and $|u(z,c)|<\delta$
in both $z \in D(R)$ and $c \in D(r)$.
Let 
$$
f_c(z):=z^2 + c+u(z,c), 
$$
$U_c:=D(R)$ and $U_c':=f_c^{-1}(U_c)$.
Then for any $c \in D(r)$
the map
$f_c:U_c' \to U_c$ is a quadratic-like map
 with a critical point $z=0$.
Moreover, it satisfies the following properties 
for sufficiently large $R$:
\begin{enumerate}[\rm (1)]
\item
There exists a family of smooth $(1+O(R^{-1}))$-quasiconformal maps (a tubing) 
$$
\Theta=\{\Theta_c: \Bar{A(R^{1/2},R)} \to \Bar{U_{c}} \sminus {U_{c}'} \}_{c\,\in\, D(r)}
$$
such that 
\begin{itemize}
\item $\Theta_c$ is identity on $\partial D(R)$ for each $c$; 
\item
$\Theta_c$ is equivariant on the boundary, i.e., 
$\Theta_c(z^2)=f_c(\Theta_c(z))$ on $\partial D(R^{1/2})$; and 
\item
for each $z \in \Bar{A(R^{1/2},R)}$
 the map $c \mapsto \Theta_c(z)$ is holomorphic in $c \in D(r)$. 
\end{itemize}
\item
Each $\Theta_c$ induces a straightening map 
$h_c$ defined on $\C$
that is uniformly $(1+O(R^{-1}))$-quasiconformal for $c \in D(r)$.
\end{enumerate}
\end{lem}
\noindent
Thus we obtain an analytic family of quadratic-like maps
 $\bs{f}=\skakko{f_c:U_c' \to U_c}_{c \,\in \,D(r)}$.
(Note that $\bs{f}$ with tubing $\Theta$ 
is not necessarily a Mandelbrot-like family.) 

\paragraph{\bf Proof.}
One can check that $f_c:U_c' \to U_c$ is a quadratic-like map
as in Example 1 of \cite[p.329]{Douady-Hubbard 1985}: 
Indeed, if $w \in \Bar{U_c}=\Bar{D(R)}$ and 
$R > \max\{2(r+\delta), 6\}$, 
then $R^2/4 > 3R/2 >r+\delta +R$ and 
thus the equation $f_c(z)=w$ 
has two solutions in $D(R/2)$ by Rouch\'e's theorem. 
This implies that $\overline{U_c'} \subset D(R/2)$.
By the maximum principle, 
$f_c:U_c' \to U_c$ is a proper branched covering of degree two.
Since $|f_c(0)| \le r+\delta <R/2$,
we have $f_c(0) \in U_c$ and hence 
the critical point $0$ is contained in $U_c'$.
Thus the Riemann-Hurwitz formula
implies that $U_c'$ is a topological disk contained in $D(R/2) \Subset U_c$.

Next we construct $\Theta$:
Let $w(a):=Re^{ia}~(0 \le a \le 4\pi)$ be a path making two turns 
along $\partial U_c$.
Let $z=g_c(w)$ be the univalent branch of $f_c^{-1}(w)=\sqrt{w-(c+u)}$ 
defined on the disk centered at $w(0)=R$ with radius $R/2$
such that $g_c(R)$ is close to $R^{1/2}$.
(Note that if $|z|=R>6$
we have $|f_c(z)| \ge R^2-(r+\delta) >R^2-R/2>3R/2$.
Thus $D(w(0),R/2) \subset f_c(D(R))$,
and by $f_c(0) \notin D(w(0),R/2)$ 
we obtain such a univalent branch.)
We take the analytic continuation of this branch $g_c$ 
along the path $w(a)~(0 \le a \le 4\pi)$ that is univalent on each $D(w(a), R/2)$,
in such a way that a path $g_c(w(a))~(0 \le a \le 4\pi)$ 
makes one turn along $\partial U_c'$. 
We will construct $\Theta$ by interpolating the paths 
$w(t)=R e^{it} \in \partial U_c$
and $z(t)= g_c(w(2t)) \in \partial U_c'$ for $0 \le t \le 2 \pi$. 

For the analytic continuation $g_c$ on $\Bar{D(w(t), R/2)}$ as above,
we have 
$$
\log z = \log g_c(w)= 
\frac{1}{2}\log w + \frac{1}{2} \log \paren{1-\frac{c+u}{w}},
$$
where $u=u(g_c(w),c)$ by taking appropriate branches of the complex logarithm.
Let
$$
\Upsilon_c(w):=\frac{1}{2} \log \paren{1-\frac{c+u}{w}},
$$
where we choose a branch of the logarithm such that $\log 1 = 0$.
Then we obtain $|\Upsilon_c(w)| =O(R^{-1})$
and hence
$\dfrac{d\Upsilon_c}{dw}(w(2t))$ $= O(R^{-2})$
by the Schwarz lemma.
Let 
$$
v_c(t):=\Upsilon_c(w(2t)) \qquad (0 \le t  \le 2 \pi)
$$
such that 
$t \mapsto w(t)=\exp(\log R+ it)$ 
and 
$t \mapsto \exp((\log R)/2+ it+v_c(t))$
parametrize the boundaries of $U_c$ and $U_c'$.
To give a homeomorphism between the closed annuli
$\Bar{A(R^{1/2},R)}$ and
$\Bar{A_c}:=\Bar{U_c} \sminus {U_c'}$, we take their logarithms:
Set $\ell:=(\log R)/2$, and consider 
the rectangle $E:=\{ s+it \st \ell \le s \le 2\ell,  0\le t \le 2\pi\}$. 
Fix a smooth decreasing function $\eta_0:[0,1] \to [0,1]$
such that: $\eta_0(0)=1$; $\eta_0(1)=0$; 
and the $j$-th derivative 
of $\eta_0$ tends to $0$ as $x  \searrow 0$ and as 
$x \nearrow 1$ for any $j$.
Set 
$$
\eta(s):=\eta_0(s/\ell -1),\quad s \in [\ell, 2\ell]. 
$$ 
Then we have $\dfrac{d\eta}{ds}(s) = O(\ell^{-1})$, 
$v_c(t) =O(R^{-1})$, and 
$\dfrac{dv_c}{dt}(t)=O(R^{-1})$.

\begin{figure}[htbp]
\begin{center}
\includegraphics[width=.65\textwidth]{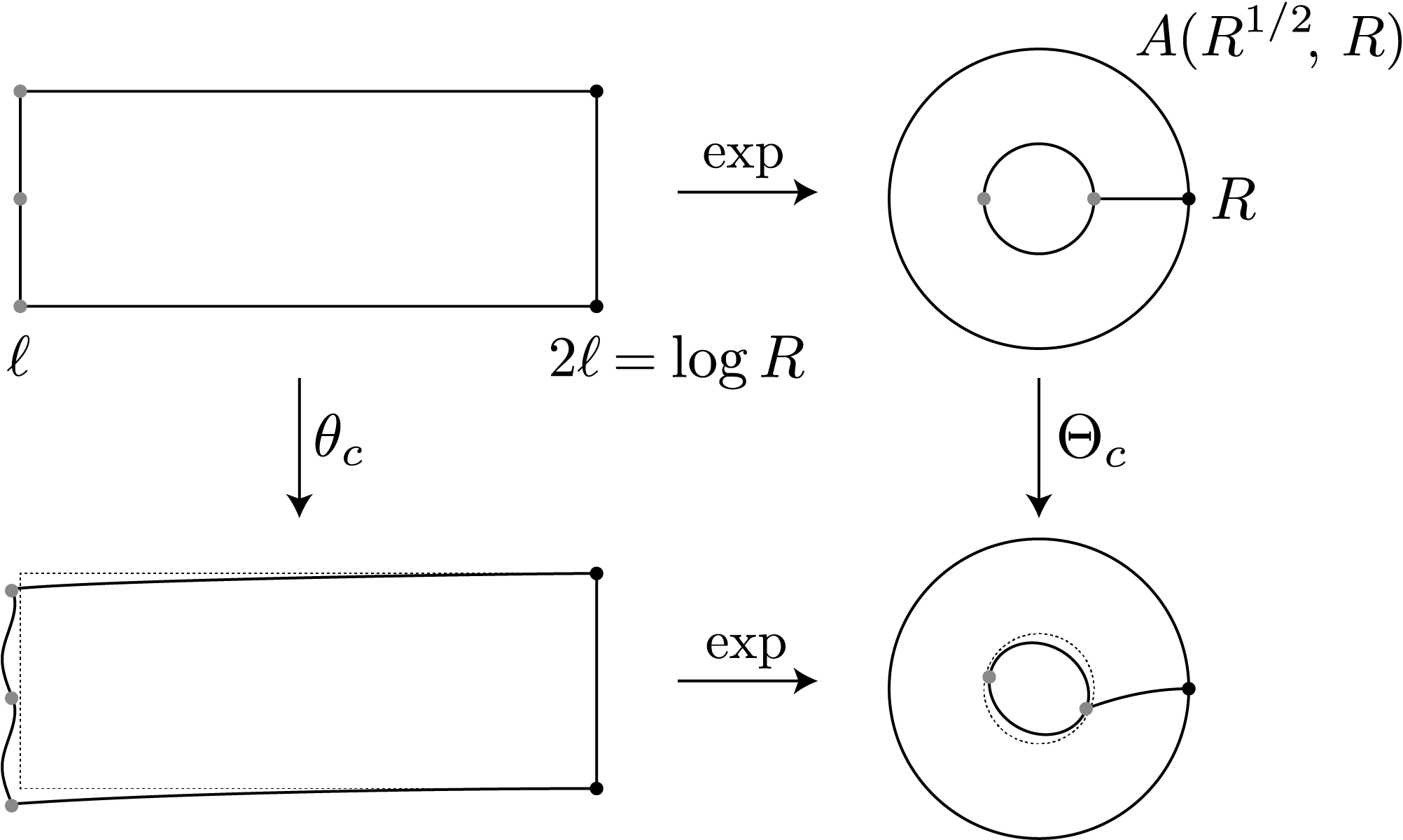}
\end{center}
\caption{\small Construction of the tubing $\Theta$ for $\bs{f}$.}
\label{fig_tubing}
\end{figure}
Now consider the smooth map $\theta_c: E \to \C$ defined by
$$
\theta_c(s+it):= s+it + \eta(s)v_c(t).
$$
The map $\theta_c$ is injective for sufficiently large $R$ since 
\begin{align}
&|\theta_c(s+it)-\theta_c(s'+it')| \notag \\ 
\ge& |(s-s') + i(t-t')| 
- |\eta(s)(v_c(t)-v_c(t'))| 
- |(\eta(s)-\eta(s'))v_c(t')|  \notag \\
\ge &
|(s-s') + i(t-t')| 
- O(R^{-1})|t-t'|
- O(\ell^{-1}R^{-1})|s-s'|. \label{eq_injectivity}
\end{align}
The Beltrami coefficient of $\theta_c$ is given by 
$$
\mu_{\theta_c} = 
\frac{(d\eta/ds)\,v_c+\eta \,(dv_c/dt) i}
{2+(d\eta/ds)\, v_c-\eta\, (dv_c/dt) i}
=O(R^{-1}).
$$
Hence $\theta_c$ is an orientation preserving diffeomorphism onto its image
for sufficiently large $R$,
and its maximal dilatation is bounded by $1+O(R^{-1})$.
By observing $\theta_c$
through the exponential function,
we obtain a smooth $(1+O(R^{-1}))$-quasiconformal homeomorphism 
$\Theta_c: \Bar{A(R^{1/2},R)} \to \Bar{A_c}$ 
that fixes the outer boundary
and satisfies 
$\Theta_c(z^2)=f_c(\Theta_c(z))$ on the inner boundary.
Holomorphic dependence of $c \mapsto \Theta_c(z)$
for each fixed $z \in \Bar{A(R^{1/2},R)}$
is obvious by the construction of $\theta_c$.

Finally we construct the straightening map $h_c:U_c \to \C$ 
of $f_c:U_c' \to U_c=D(R)$.
Let us extend $f_c$ to a smooth quasiregular map $F_c:\C \to \C$ by setting
$$
F_c(z):=
\left\{\begin{array}{ll}
f_c(z) 
& \text{if}~z \in U_c',\\[.5em] 
\skakko{\Theta_c^{-1}(z)}^2 
& \text{if}~z \in U_c  \sminus  U_c', ~\text{and}\\[.5em]
z^2  
& \text{if}~z \in \C \sminus U_c.
\end{array}\right.
$$
We define an $F_c$-invariant 
 Beltrami coefficient $\mu_c$ (i.e., $F_c^\ast\mu_c=\mu_c$) by
$$
\mu_c(z):=
\left\{\begin{array}{ll}
0 & 
\text{if}~z \in K(f_c)~\text{or}~ z \in \C \sminus U_c,\\[.5em] 
\dfrac{(F_c)_{\Bar{z}}(z)}{(F_c)_z(z)} 
& \text{if}~z \in U_c  \sminus  U_c',~\text{and}\\[1.2em]
\disp (f_c^n)^\ast \mu_c(z)
& \text{if}~f_c^n(z) \in U_c  \sminus  U_c'~\text{for some}~ n >0,
\end{array}\right.
$$
where 
$$
\disp  (f_c^n)^\ast \mu_c(z)
=\mu_c(f_c^n(z))\frac{\Bar{(f_c^n)'(z)}}{(f_c^n)'(z)}.
$$
Then $\mu_c$ is supported on $\Bar{D(R)}$ 
and it satisfies $\norm{\mu_c}_\infty=O(R^{-1})$.
By existence of the normal solutions of the Beltrami equations \cite[Theorem 4.24]{IT Book}, 
we have a unique $(1+ O(R^{-1}))$-quasiconformal map
$h_c:\C \to \C$ that satisfies the Beltrami equation 
$(h_c)_{\Bar{z}}=\mu_c \cdot (h_c)_z$ a.e.,
$h_c(0)=0$, and  
$(h_c)_z-1 \in L^p(\C)$ for some $p >2$.
(The relation between $R$ and $p$ will be more specified in the next lemma.)
The condition $(h_c)_z-1 \in L^p(\C)$ implies $w=h_c(z)=z+b_c+O(1/z)$ 
as $z \to \infty$ for some constant $b_c \in \C$.

Since $\mu_c$ is $F_c$-invariant, 
the map $P(w):=h_c \cc F_c \cc h_c^{-1}(w)$ is a 
holomorphic map of degree 2 with 
a critical point at $h_c(0)=0$ 
and a superattracting fixed point at $h_c(\infty)=\infty.$
Hence $P(w)$ is a quadratic polynomial.
The expansion of the form $w=h_c(z)=z+b_c+O(1/z)$
implies that we actually have $b_c=0$
and $P(w)$ is of the form 
$P(w)=w^2+\chi(c)=P_{\chi(c)}(w)$. Hence the restriction $h_c|_{U_c}$ 
is our desired straightening map.
\QED

The next lemma shows that the quasiconformal 
map $h_c:\C \to \C$ constructed above is uniformly close 
to the identity on compact sets for sufficiently large $R$:

\begin{lem}\label{lem_estimate_of_h}
Fix any $p>2$ and any compact set $E \subset \C$.
If $R$ is sufficiently large, 
then 
the quasiconformal map
$h_c$ in Lemma \ref{lem_almost_conformal} 
satisfies 
$$
|h_c(z)-z|=O(R^{-1+2/p})
$$
uniformly for each $c \in D(r)$ and $z \in E$.
\end{lem}
\noindent
Indeed, the estimate is valid for any $R \ge C_0p^2$, 
where $C_0$ is a constant independent of $p$.

\paragraph{\bf Proof.}
We have 
$\norm{\mu_c}_\infty \le C/R=:k$ for 
some constant $C$ independent of $c \in D(r)$
by the construction of $h_c$.
By \cite[Theorem 4.24]{IT Book}, 
we have $(h_c)_z-1 \in L^p(\C)$ for any $p>2$ satisfying $k C_p<1$, 
where $C_p$ is the constant that appears in the Calderon-Zygmund inequality
\cite[Proposition 4.22]{IT Book}.
Gaidashev showed in \cite[Lemma 6]{Gaidashev 2007} that $C_p \le \cot^2(\pi/2p)$.
Since $\cot^2(\pi/2p) = (2p/\pi)^2(1+o(1))$ as $p \to \infty$,
the inequality $k C_p \le (C/R) \cot^2(\pi/2p) \le 1/2$ 
is established if we take $R \ge C_0 p^2$ 
for some constant $C_0$ independent of $p>2$. 
By following the proof of \cite[\S 4, Corollary 2]{IT Book}, 
we have 
$$
|h_c(z)-z|
\le 
K_p \cdot \frac{1}{1-kC_p}\norm{\mu_c}_p|z|^{1-2/p}  
$$
for any $z \in \C$, 
where $K_p>0$ is a constant depending only on $p$
and $\norm{\mu_c}_p$ is the $L^p$-norm of $\mu_c$.
Since $|\mu_c|\le C/R$ and $\mu_c$ is supported on $\overline{D(R)}$,
we have $\norm{\mu_c}_p \le (C/R)(\pi R^2)^{1/p} =C \pi^{1/p} R^{-1+2/p}$.
Hence if we take $R \ge C_0 p^2$ such that $kC_p \le 1/2$,
we have 
$$
|h_c(z)-z| \le 2 K_p C \pi^{1/p} R^{-1+2/p} |z|^{1-2/p}.
$$
This implies that $|h_c(z)-z| =O( R^{-1+2/p})$
on each compact subset $E$ of $\C$. 
\QED 

\begin{cor}\label{cor_estimate_of_chi}
Fix any $p >2$. 
If $R$ is sufficiently large,
then for each $c \in D(r)$, 
$f_c$ is hybrid equivalent to 
a quadratic polynomial $P_{\chi(c)}(w)=w^2+\chi(c)$ 
with 
$$
|\chi(c)-c| = O(\delta) + O(R^{-1+2/p}).
$$
\end{cor}

\paragraph{\bf Proof.}
We have $\chi(c)=h_c(f_c(0))$ 
since $h_c$ maps the critical value of $f_c$ 
to that of $P_{\chi(c)}$, 
For each $c \in D(r)$,
$f_c(0) = c+u_c(0)$ is contained in a compact set $\overline{D(r+\delta)}$.
By Lemma \ref{lem_estimate_of_h}, we obtain
$$
\chi(c)=h_c(f_c(0))=c+u_c(0) + O(R^{-1+2/p})
$$
for sufficiently large $R$ and this implies the desired estimate.
\QED 

\medskip

\paragraph{\bf Coordinate changes}
Under the same assumption as in Lemma \ref{lem_almost_conformal},
we assume in addition that  
\begin{enumerate}[(i)]
\item
$\delta<1$ and  $r >4$; and 
\item
$R$ is large enough such that 
${\bs{f}}=\skakko{f_c:U_c' \to U_c}_{c\,\in \,D(r)}$
is an analytic family of quadratic-like maps 
(that may not necessarily be Mandelbrot-like),
 and that (1) and (2) of Lemma \ref{lem_almost_conformal} hold.
\end{enumerate}
By (2) of Lemma \ref{lem_almost_conformal}, 
each $f_c \in \bs{f}$ is hybrid equivalent to 
some quadratic map $P_{\chi(c)}$
by the $(1 + O(R^{-1}))$-quasiconformal straightening $h_c:U_c \to h_c(U_c)$.
We say the map
$$
\chi = \chi_{\bs f}:D(r) \to \C,\quad \chi(c)=h_c(f_c(0))
$$
is the {\it straightening map} of the family ${\bs f}$
associated with the tubing $\Theta=\{\Theta_c\}_{c \,\in\, D(r)}$.
We say the map $(z,c) \mapsto (h_c(z), \chi(c))$
defined on $D(R)\times D(r)$ 
is a {\it (straightening) coordinate change}.

Now we show that the straightening map is quasiconformal with
dilatation arbitrarily close to $1$ if we take sufficiently large $R$
and small $\delta$:

\begin{lem}[\bf Almost conformal straightening of ${\bs f}$]
\label{lem_almost_conformal_M}
If $r$ and $R$ are sufficiently large
and $\delta>0$ is sufficiently small,
then the family ${\bs f}$ is associated with 
a $(1+O(\delta)+O(R^{-1}))$-quasiconformal straightening map
$$
\chi=\chi_{\bs f}:D(r) \to \chi(D(r)) \subset \C
$$
such that
\begin{enumerate}[\rm (1)]
\item $\chi(M_{\bs f})=M$, where $M_{\bs f}$ 
is the connectedness locus of ${\bs f}$;
\item $\chi|_{D(r) \sminus M_{\bs f}}$ is $(1+O(R^{-1}))$-quasiconformal; and
\item $\chi|_{M_{\bs f}}$ extends to a $(1+O(\delta))$-quasiconformal map on the plane.
\end{enumerate}
\end{lem}

\paragraph{\bf Proof.}
By slightly shrinking $r>4$ if necessary, we may assume that 
$f_c$ is defined for  $c \in \partial D(r)$. 
When $c=re^{it}~(0 \le t \le 2 \pi)$, 
we have $|f_c(0)-0|=|c+u(0,c)| \ge r-\delta > 3$ (since $\delta<1$).
Hence as $c$ makes one turn around the origin 
so does $f_c(0)$.
By \cite[p.328]{Douady-Hubbard 1985}, 
$\chi$ gives a homeomorphism between $M_{\bs f}$ and $M$.
Moreover,  
$\chi|_{D(r) \sminus M_{\bs f}}$ is $(1+O(R^{-1}))$-quasiconformal
by \cite[Proposition 20, Lemma in p.327]{Douady-Hubbard 1985},
since each $\Theta_c$ is $(1+O(R^{-1}))$-quasiconformal.

For the dilatation of $\chi|_{M_{\bs f}}$,
we follow the argument of \cite[Lemma 4.2]{McMullen 2000}:
Consider the families 
${\bs f}_t:=\braces{f_{c,t}}_{c\, \in \,D(r)}$
defined for each $t \in \D$,
where
$$
f_{c,t}(z):=z^2+c+\frac{\,t\,}{\delta}\,u(z,c).
$$
By the same argument as above, 
the connectedness locus  $M_{{\bs f}_t}$ of ${\bs f}_t$
is homeomorphic to $M$ by the straightening map 
$$
\chi_t=\chi_{{\bs f}_t}:M_{{\bs f}_t} \to M.
$$
Then the inverse $\phi_t:=\chi_t^{-1}:M \to M_{{\bs f}_t}$
gives a holomorphic family of injections over $\D$.
By Bers and Royden's theorem \cite[Theorem 1]{Bers-Royden 1986},
each of them extends to a $(1+|t|)/(1-|t|)$-quasiconformal map 
$\widetilde{\chi}_t$ on $\C$.
In particular, $\chi|_{M_{\bs f}}=\chi_\delta$ extends 
to a $(1+O(\delta))$-quasiconformal map on $\C$.
Now we apply Lemma \ref{lem_Bers} (Bers' Gluing Lemma)
to $H_1:=\chi|_{D(r) \smallsetminus M_{\bs f}}$ and $H_2=\widetilde{\chi}_\delta$.
Then $H_1$ and $H_2$ are glued along $\partial M$
and the glued map $H:D(r)\to \C$, 
which coincides with $\chi$,
is $(1+O(\delta)+O(R^{-1}))$-quasiconformal.
\QED

\section{Proof of Theorem C}

\paragraph{\bf Idea of the proof.}
The proof follows the argument of Theorem A
and uses the results in the previous section.
Recall that in the proof of Theorem A,
 we construct two families of quadratic-like maps
$\{f_c: U_c' \to U_c\}_{c \, \in \, S \, \cap \,\Lambda}$ (``the first renormalization") 
and 
${\bs G}=\{G_c: V_c' \to U_c\}_{c \, \in \, W}$ (``the second renormalization"),
and we conclude that the small Mandelbrot set
corresponding to the family ${\bs G}$
has a desired decoration.
 
In the following proof of Theorem C, 
we first take a ``thickened" family 
$\widehat{\bs f}
=\{f_c: \widehat{U}_c' \to \widehat{U}_c\}_{c \, \in \, \widehat{\Lambda}}$ 
that contains 
${\bs f}=\{f_c: U_c' \to U_c\}_{c \, \in \, \Lambda}$
as a restriction (in both dynamical and parameter planes), 
such that $U_c \Subset \widehat{U}_c$ and 
the modulus of $\widehat{U}_c\sminus\overline{U_c}$ 
is sufficiently large.
Next we construct another ``thickened" family 
$\widehat{\bs G}=\{G_c: \widehat{V}_c' \to \widehat{U}_c\}_{c \, \in \, \widehat{W}}$
that contains $\bs{G}=\{G_c: V_c' \to U_c\}_{c \, \in \, W}$
with $V_c' \Subset \widehat{V}_c'$.
Then we can apply a slightly modified versions of the lemmas 
in the previous section to the family $\bs{G}$.
Finally we conclude that the small Mandelbrot set
corresponding to the family ${\bs G}$
has a very fine decoration.

\paragraph{\bf Notation. }
We will use a conventional notation: For complex variables $\alpha$ and $\beta$,
by $\alpha \asymp \beta$ we mean $C^{-1}|\al| \le |\beta| \le C|\al|$ for an implicit constant $C>1$. 

\paragraph{\bf First renormalization.}
We start with a result by McMullen \cite[Theorem 3.1]{McMullen 2000}
(see also \cite[Chapter V]{Douady-Hubbard 1985})
applied to (and modified for) the quadratic family:

\begin{lem}[\bf Misiurewicz cascades]
\label{lem_cascade}
For any Misiurewicz parameter $m_0$ 
and any arbitrarily large $r$ and $R$, 
there exist sequences 
$\braces{s_n}_{n \ge 1}$, 
$\braces{p_n}_{n \ge 1}$,
$\braces{t_n}_{n \ge 1}$, 
and $\braces{\delta_n}_{n \ge 1}$
that satisfy the following conditions for each sufficiently large $n$: 
\begin{enumerate}[\rm (a)]
\item
$s_n$ is a superattracting parameter of period $p_n$ 
with $|s_n -m_0| \asymp \mu_0^{-n}$, 
where $\mu_0$ is the multiplier
of the repelling cycle of $P_{m_0}$ 
on which the critical orbit lands.
\item 
$t_n \in \C^\ast$ and $t_n \asymp \mu_0^{-2n}$.
\item 
$\delta_n>0$ and $\delta_n \asymp n \mu_0^{-n}$.
\item
Let $X_n:D_n:=D(s_n,\, r |t_n|) \to D(r)$ be the affine map 
defined by 
$$
C=X_n(c):= \frac{c-s_n}{t_n}.
$$
Then there exists a non-zero holomorphic function 
$c \mapsto \al(c)=\al_c$
defined for $c \in D_n$ such that $\alpha_c \asymp \mu_0^{-n}$ 
and the map
$$
Z = A_c(z) :=\frac{z}{\alpha_c}
$$
conjugates $P_c^{p_n}$ on $D(R \,|\alpha_c|)$ to 
the map $F_{C}:=A_c \cc P_c^{p_n}\cc A_c^{-1}$ on $D(R)$
of the form
\begin{equation}
F_{C}(Z) 
 = A_c \cc P_c^{p_n} \cc A_c^{-1}(Z)
=Z^2 + C + u(Z,C),
\label{eq_F_C}
\end{equation}
where $u(Z,C) = u_C(Z)$ 
is holomorphic in both $Z \in  D(R)$ and $C \in D(r)$,
and satisfies $u'_C(0) = 0$ and $|u(Z,C)| \le \delta_n$.
\end{enumerate} 
\end{lem}

\begin{figure}[htbp]
\begin{center}
\includegraphics[width=.72\textwidth]{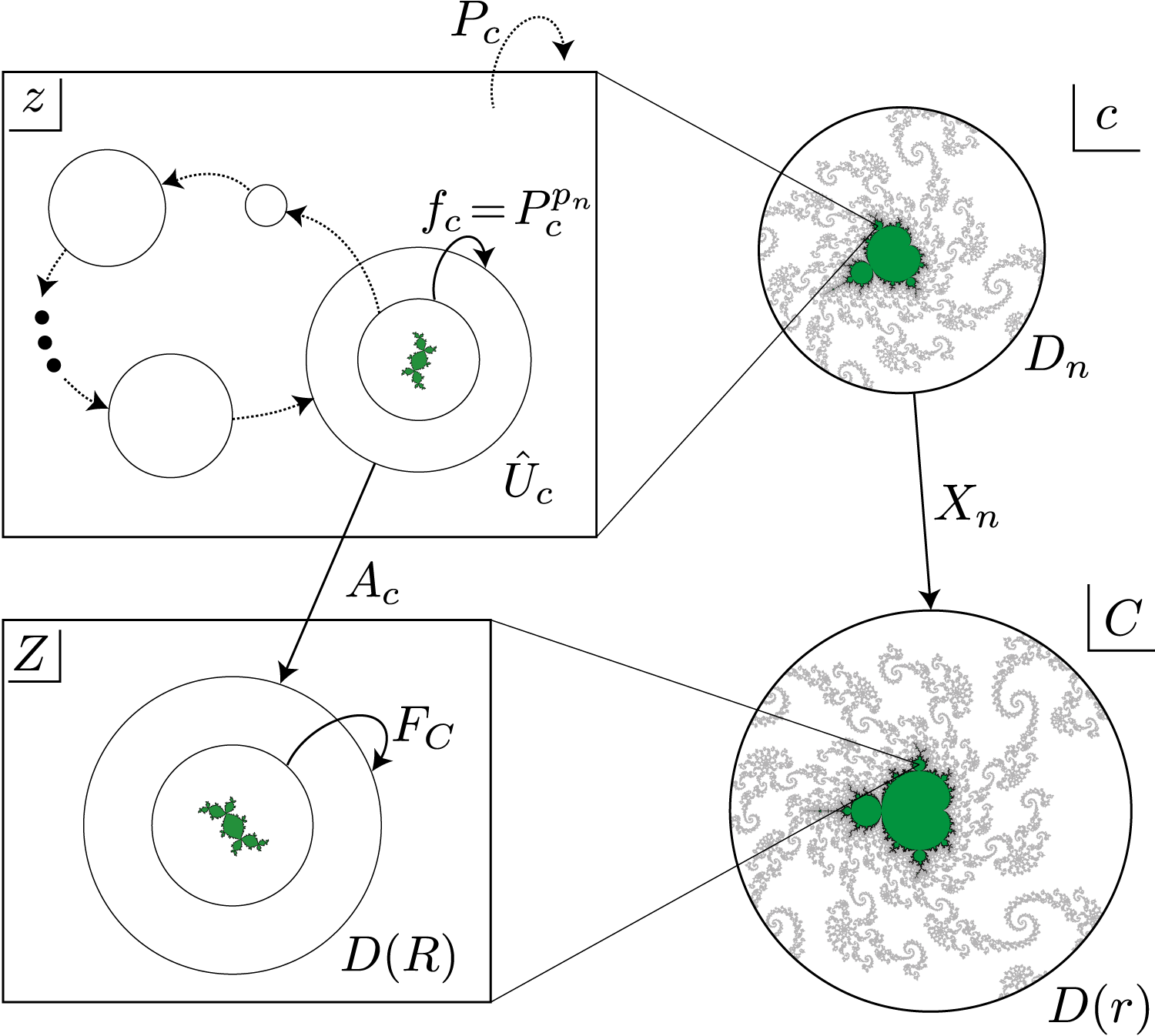}
\end{center}
\caption{\small 
An affine coordinate change.}
\label{fig_affine_coordinate_change}
\end{figure}

\paragraph{\bf Construction of the family $\widehat{\bs f}$.}
Let us fix arbitrarily small $\varepsilon>0$ 
and $\kappa>0$ as in the statement of Theorem C.
We choose any Misiurewicz parameter 
$m_0$ in $\mathrm{int}(B) \cap \partial M$,
where $B$ is the closed disk given in the statement.
(The Misiurewicz parameters are dense in $\partial M$.)

For any $r$ and $R$ bigger than $4$
(we will replace them with larger ones if necessary,
but it will happen finitely many times in what follows),
by taking a sufficiently large $n$ in Lemma \ref{lem_cascade}
such that $D_n \subset B$
\footnote{More precisely, we fix $r$ first, 
and then take a larger $R$ (and an $n$) if necessary to
apply those lemmas.},
we have an analytic family 
$\braces{F_C: D(R) \to \C}_{C \,\in \,D(r)}$
that satisfies the conditions for 
Lemma \ref{lem_almost_conformal}.
Moreover, its restriction
$$
{\bs F}_n:=
\braces{F_C: F_C^{-1}(D(R)) \to D(R)}_{C \,\in \,D(r)}
$$
is an analytic family of quadratic-like maps
that satisfies the conditions for 
Lemma \ref{lem_almost_conformal_M}.
Hence we have an associated straightening coordinate change
of the form $(Z,C) \mapsto (H_C(Z), \chi_{{\bs F}_n}(C))$.
More precisely, for each $C \in D(r)$, 
$F_C$ is hybrid equivalent to 
$Z \mapsto Z^2+\chi_{{\bs F}_n}(C)$
by a $(1+O(R^{-1}))$-quasiconformal straightening $H_C$
by Lemma \ref{lem_almost_conformal}, 
and $H_C$ satisfies the estimate of Lemma \ref{lem_estimate_of_h}.
By Lemma \ref{lem_almost_conformal_M}, 
the straightening $\chi_{{\bs F}_n}: D(r) \to \C$
is $(1+O(\delta_n)+O(R^{-1}))$-quasiconformal
and satisfies the estimate of Corollary \ref{cor_estimate_of_chi}.
Hence we may assume that $R$ and $n$ are large enough such that
both $Z \mapsto H_C(Z)$ and $C \mapsto \chi_{{\bs F}_n}(C)$
are $(1+\kappa)^{1/2}$-quasiconformal for 
$\kappa >0$ given in the statement.

Let $f_c:=P_c^{p_n}:\widehat{U}_c' \to \widehat{U}_c$
be the pull-back of $F_C:F_C^{-1}(D(R)) \to D(R)$
by the map $(z,c) \mapsto (Z,C)=(A_c(z),X_n(c))$,
which we call the {\it affine coordinate change}.
(See Figure \ref{fig_affine_coordinate_change}.
Note that $\widehat{U}_c=D(R \, |\alpha_c|)$ is a round disk.)
Set 
$$
p:=p_n,\quad s_0:=s_n,\quad \text{and} \quad 
\widehat{\Lambda}:=D_n=D(s_n,\, r |t_n| ).
$$
The quadratic-like family
$$
\widehat{{\bs f}}:=
\big\{f_c:\widehat{U}_c'\to \widehat{U}_c \big\}
_{c \,\in \,\widehat{\Lambda}}
$$
is our first family of renormalizations
whose straightening coordinate change $(z,c) \mapsto (h_c(z),\chi(c))$
is given by
$$
(h_c(z),\chi(c)):= (H_C \cc A_c(z), \chi_{{\bs F}_n} \cc X_n(c)).
$$
Note that both $h_c:\widehat{U}_c \to \C$ 
and $\chi:\widehat{\Lambda} \to \C$
are $(1+\kappa)^{1/2}$-quasiconformal.
By Lemma \ref{lem_estimate_of_h} and Corollary \ref{cor_estimate_of_chi},
if we fix any $p'>2$
and any compact subset $E$ of $\widehat{U}_c$,
then for sufficiently large $R$ we have 
$$
h_c(z) =A_c(z) +O(R^{-1+2/p'})
$$
on $E$ 
and 
$$
\chi(c)= X_n(c)+ O(\delta_n)+O(R^{-1+2/p'})
$$
on $\widehat{\Lambda}$. 
Hence the straightening coordinate change 
$(z,c) \mapsto (h_c(z),\chi(c))$
is very close to the affine coordinate change 
$(z,c) \mapsto (A_c(z),X_n(c))=(z/\al_c,(c-s_n)/t_n)$
if we take sufficiently large $R$ and $n$.

\paragraph{\bf Construction of the family ${\bs f}$.}
Let $\rho>4$ be an arbitrarily large number.
By taking sufficiently large $r$, $R$ and $n$
such that $\rho/R$ is sufficiently small,
we may assume the following:
\begin{itemize}
\item
The set $\Omega(\rho):=\skakko{C \in D(r)\st F_C(0) \in D(\rho)}$
gives a Mandelbrot-like family
$$
\bs{F}_n(\rho):=\{F_C:F_C^{-1}(D(\rho)) \to D(\rho)\}_{C \,\in \,\Omega(\rho)}.
$$
\item 
$D(\rho) \Subset F_C^{-1}(D(R))$ for any $C \in \Omega(\rho)$.
\end{itemize}
Now we define the Mandelbrot-like family
$$
\bs{f}=\{f_c:U_c' \to U_c\}_{c \,\in \,\Lambda}
$$ 
as the pull-back of $\bs{F}_n(\rho)$ by the affine coordinate change 
$(z,c) \mapsto (A_c(z),X_n(c))$ above.
More precisely,
we let 
$U_c:=A_c^{-1}(D(\rho))=D(\rho|\al_c|)$
for each $c \in \widehat{\Lambda}$,
and consider the restriction $f_c:U_c' \to U_c$ of
$f_c:\widehat{U}_c' \to \widehat{U}_c$.
Then we define the subset $\Lambda$ of $\widehat{\Lambda}$ by 
$\Lambda:=X_n^{-1}(\Omega(\rho))$ such that the family $\bs{f}$ above
becomes a Mandelbrot-like family.
Note that we have 
$U_c' \Subset U_c \Subset \widehat{U}_c' \Subset \widehat{U}_c$
for any $c \in \Lambda$.

Let $M_{\bs f} =M_{s_0} = s_0 \perp M$ be the connectedness locus of 
the family ${\bs f}$, which coincides with that of $\widehat{\bs f}$. 
Note that $\bs f$ has the same straightening coordinate change 
$(z,c) \mapsto (h_c(z),\chi(c))$ as $\widehat{\bs f}$
such that $\chi(M_{\bs f})=M$.

\paragraph{\bf Second renormalization.}
For a given Misiurewicz or 
parabolic parameter $c_0$ in the statement of Theorem C, 
we define a Misiurewicz or parabolic parameter $c_1 \in M_{\bs f}$ by
$$
c_1: = \chi^{-1}(c_0) =s_0 \perp c_0.
$$ 
Let $q_{c_1}$ be a repelling or parabolic periodic point 
of $f_{c_1}=P_{c_1}^p$ of some period $k$ that belongs to the postcritical set. 
More precisely, when $c_1$ is Misiurewicz, 
there exist minimal integers $l$ and $k$ such that 
$q_{c_1}=f_{c_1}^l(0)=f_{c_1}^{l+k}(0)$ 
and $q_{c_1}$ is repelling.
When $c_1$ is parabolic, 
the orbit of $0$ accumulates on 
a parabolic periodic point $q_{c_1}$ with $(f_{c_1}^{k\nu})'(q_{c_1})=1$, where $\nu$ is the petal number
(see \cite{Kawahira-Kisaka 2023} for more details).
Let $\Omega_{c_1}$ be the domain  
of the linearizing coordinate 
or the attracting Fatou coordinate of $q_{c_1}$.
In both cases, we may assume that 
$f_{c_1}^l(0)$ is contained in $\Omega_{c_1}$ 
 for some $l \ge 1$.

The rest of the proof of Theorem C
is divided into Claims 1 to 7 below and their proofs.
For the first two claims,
we may simply apply the argument of
Steps (M1)--(M2) in Section 4: 

\paragraph{\bf Claim 1.}
{\it  
There exists a Jordan domain $\widehat{V}_{c_1}$ 
with $C^1$ boundary and integers $N,\,j \in \N$ 
which satisfy the following:
\begin{enumerate}[\rm (1)]
\item 
$\widehat{V}_{c_1}$ is a connected component of 
$P_{c_1}^{-N}(\widehat{U}_{c_1})$.
\item 
$g_{c_1}:=P_{c_1}^N|_{\widehat{V}_{c_1}}:\widehat{V}_{c_1} \to \widehat{U}_{c_1}$ 
is an isomorphism and 
${f_{c_1}^{j}(\widehat{V}_{c_1})} 
\Subset \widehat{U}_{c_1}  \sminus  \Bar{\widehat{U}_{c_1}'}$. 
\item
$\widehat{V}_{c_1} \Subset \Omega_{c_1}$.
Also we can take $\widehat{V}_{c_1}$ 
arbitrarily close to $q_{c_1}$.
\end{enumerate}
}

\medskip

\medskip

\paragraph{\bf Claim 2.}
{\it
There exists 
a Jordan domain $\widehat{W} \Subset \Lambda \sminus M_{\bs f}
~(\subset \widehat{\Lambda} \sminus M_{\bs f})$ 
arbitrarily close to $c_1$ that satisfies the following:
\begin{enumerate}[\rm (1)]
\item
There is a holomorphic motion of $\widehat{V}_{c_1}$
over $\widehat{W}$ that generates 
a family of Jordan domains $\{\widehat{V}_{c}\}_{c \,\in \, \widehat{W}}$
with $C^1$ boundaries
such that for each $c \in \widehat{W}$, 
${f_c^j(\widehat{V}_{c} )} \Subset \widehat{U}_{c} \sminus \overline{\widehat{U}_{c}'}$
and 
$g_c:=P_{c}^{N}|_{\widehat{V}_{c}}:\widehat{V}_{c} \to \widehat{U}_{c}$
is an isomorphism.
\item
There exists an $L$ such that 
$f_{c}^{L}(0)=P_{c}^{pL}(0) \in \widehat{V}_{c}$ for any $c \in \widehat{W}$.
\item 
For $c \in \partial \widehat{W}$, 
we have $P_{c}^{pL+N}(0) \in \partial  \widehat{U}_c$.
Moreover, when $c$ makes one turn along $\partial \widehat{W}$,
then $P_{c}^{pL+N}(0)$ makes one turn around the origin.
\end{enumerate}
}
Indeed, such a domain $\widehat{W}$ and an $L$ 
are given explicitly as follows 
(though we will not use these details): 
In the Misiurewicz case 
we let $L:=l+k n'$ for some sufficiently large $n'$ and
$$
\widehat{W}
=\widehat{W}_{n'}
:=
\braces{c \in \Lambda \sminus M_{\bs f} 
\st 
f_c^{l+ki}(0) \in \Omega_{c}~\text{for 
$i=1,\,\ldots, n'-1$ and}~ 
f_c^{L}(0) \in \widehat{V}_{c}},
$$
where $\Omega_c$ is the domain of
the linearizing coordinate satisfying the conditions in
the last paragraph of Step (M2).
In the parabolic case, 
we let $L:=l+k \nu n'$ for some sufficiently large $n'$ and
$$
\widehat{W}
=\widehat{W}_{n'}
:=
\braces{c \in \Lambda \sminus M_{\bs f} 
\st 
f_c^{l+k\nu i}(0) \in \Omega_{c}^\ast~\text{for 
$i=1,\,\ldots, n'-1$ and}~ 
f_c^{L}(0) \in \widehat{V}_{c}},
$$
where $\Omega_c^\ast$ is the domain of
the perturbed Fatou coordinates.
(See \cite{Kawahira-Kisaka 2023} for more details).

\begin{figure}[htbp]
\begin{center}
\includegraphics[width=.65\textwidth]{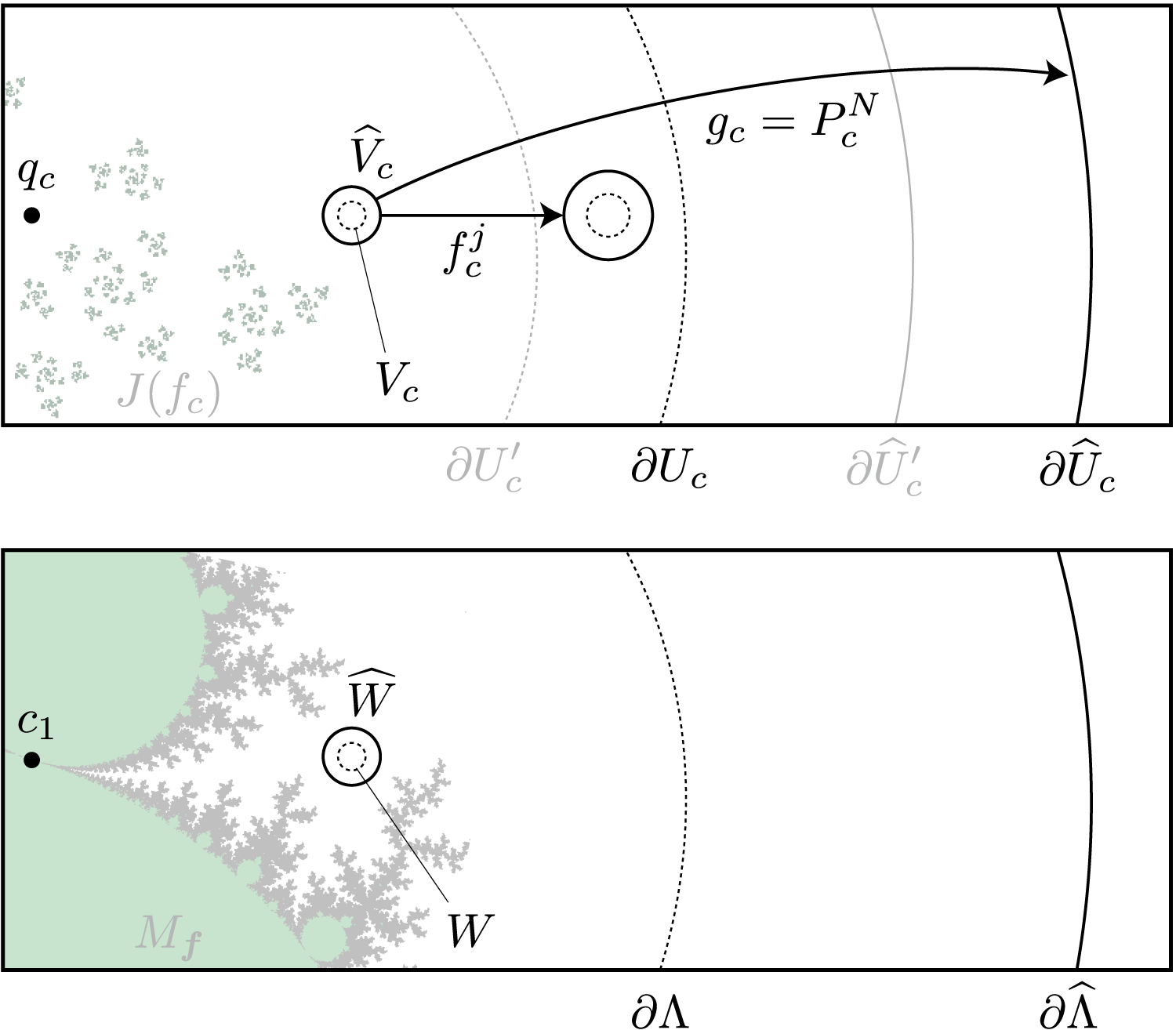}
\caption{\small
For any $c \in \widehat{W}$,
the Julia set $J(f_c)$ moves only a little from 
$J(f_s)$ of the center $s \in W \Subset \widehat{W}$.
}\label{fig_W}
\end{center}
\end{figure}

Note that if we take sufficiently large 
$N$ and $j$, we may always assume that 
$\widehat{V}_c \Subset U_c'$ as depicted in Figure \ref{fig_W}.
Moreover, the proof of Lemma 4.1 indicates that
we can choose $\widehat{V}_c$ with arbitrarily small diameter.

\medskip

\paragraph{\bf Definition of the center of $\widehat{W}$ and $W$.}
As in Step (M2), 
there exists a unique 
superattracting parameter $s \in \widehat{W}$ 
such that $P_{s}^{pL+N}(0)=0$.
We call $s$ the {\it center} of $\widehat{W}$.

By Claim 2, we can find 
a family of Jordan domains $\{V_c\}_{c \, \in \, \widehat{W}}$
with $C^1$ boundaries 
such that $V_c \Subset \widehat{V}_c$
and that $P_c^N:V_c \to U_c$ is an isomorphism for each $c \in \widehat{W}$.
Let $W \Subset \widehat{W}$ be the set of 
$c$ such that $f_c^L(0) \in V_c$.
By the same argument as in Step (M3), 
one can check that both $\widehat{W}$ and $W$ are 
Jordan domains with $C^1$ boundaries.
Note that $s$ is the center of $W$ as well.
(See Figure \ref{fig_W} again.)

Moreover, we have:

\paragraph{\bf Claim 3 (Straightening the center).}
{\it
By choosing $\widehat{V}_c$ in Claim 1 
close enough to $q_{c_1}$,
we can find an $\eta \in D(\vep)$ 
with $c_0 +\eta \notin M$
such that $f_s:\widehat{U}_s' \to \widehat{U}_s$ is 
hybrid equivalent to 
a quadratic-like restriction of $P_{c_0+\eta}$ 
 with $(1+\kappa)^{1/2}$-quasiconformal straightening map.
In particular, $J(f_s)$ is a $(1+\kappa)^{1/2}$-quasiconformal
image of $J(P_{c_0+\eta})$.
}

\paragraph{\bf Proof.}
By the construction of $\widehat{W}$ in Claim 2
(following Step (M2)), 
we can take $\widehat{W}$ arbitrarily close to $c_1$.
Since the straightening map
$\chi:\widehat{\Lambda} \to \C$
is continuous, 
we have $|\chi(s)-\chi(c_1)| = |\chi(s)-c_0|<\vep$
by taking $s \in \widehat{W}$ close enough to $c_1$.
Set $\eta:=\chi(s)-c_0$.
Then we have
$\chi(s)=c_0+\eta \in \C  \sminus M$ 
since $\widehat{W} \subset \widehat{\Lambda} \sminus M_{\bs f}$.
By the construction of the first renormalization,
$f_s$ is conjugate to $P_{\chi(s)}$ 
by the $(1+\kappa)^{1/2}$-quasiconformal straightening map $h_s$
such that $h_s(J(f_s))=J(P_{c_0+\eta})$.
\QED

\medskip

\paragraph{\bf Remark.}
Since $h_s(z)= H_{X_n(s)} \cc A_s(z)=z/\al_s + o(1)$,
 $J(f_s)$ is actually an ``almost affine" (even better than ``almost conformal"!)
 copy of $J(P_{c_0+\eta})$.

\paragraph{\bf Holomorphic motion of the Cantor Julia sets.}
Since 
$\widehat{W} \subset \widehat{\Lambda} \sminus M_{{\bs{f}}}$,
the Julia set $J(f_c)$ for each $c \in \widehat{W}$
is a Cantor set that is a $(1+\kappa)^{1/2}$-quasiconformal image of $J(P_{\chi(c)})$.
Moreover, the Julia set $J(f_c)$ moves holomorphically for $c \in \widehat{W}$:

\medskip 
\paragraph{\bf Claim 4 (Cantor Julia moves a little).}
{\it
There exists a holomorphic motion $\iota: J(f_s) \times \widehat{W} \to \C$
such that 
$\iota_c(z):=\iota(z,c)$ maps $J(f_s)$ bijectively to 
$J(f_c)$ for each $c \in  \widehat{W}$. 
Moreover, if $R$ is sufficiently large,
then $\iota_c$ extends to a $(1+\kappa)^{1/2}$-quasiconformal homeomorphism on the plane 
for each $c \in W \Subset \widehat{W}$.
}

\medskip 
A direct corollary of Claims 3 and 4 is:

\begin{cor}[\bf Julia appears in Julia]\label{cor_almost_conformal_copy_of_J}
The Julia set $J(P_c)$ of $P_c$ contains a
$(1+\kappa)$-quasiconformal copy of $J(P_{c_0+\eta})$
for any $c \in W  \Subset \widehat{W}$.
\end{cor}

\paragraph{\bf Proof of Claim 4.}
Since $J(f_c)$ is a hyperbolic set for each $c \in  \widehat{W}$,
 it has a local holomorphic motion near $c$. 
 (See \cite[p.229]{Shishikura 1998}.)
The holomorphic motion extends to that of 
$J(f_s)$ over $\widehat{W}$ 
as in the statement, 
since $\widehat{W}$ is simply connected (and isomorphic to $\D$). 

Now we consider the modulus of the annulus $\widehat{W} \sminus \overline{W}$: 
Recall that $A_c(\widehat{U}_c)=D(R)$, where $A_c(z)=z/\al_c$
and $\al_c$ depends holomorphically on $c \in \widehat{W}$ 
(Lemma \ref{lem_cascade}).
By the same argument as in Step (M2), 
for each $\zeta \in D(R)$, the equation
$$
f_{c}^L(0) = (A_c \circ  g_c)^{-1}(\zeta) \quad (\in \widehat{V}_c)
$$ 
with respect to $c$ has a unique solution $c=\check{c}(\zeta) \in \widehat{W}$ and the map $\zeta \mapsto \check{c}(\zeta)$ gives an isomorphism from $D(R)$ onto $\widehat{W}$.
In particular, we have $A_c(\widehat{U}_c \sminus \overline{U_c})=A(\rho,R)$
for any $c \in \widehat{W}$ and thus $\widehat{W} \sminus \overline{W}=\check{c}(A(\rho,R))$.
Hence we obtain
$$
\mathrm{mod}(\widehat{W} \sminus \overline{W})
= \mathrm{mod}(A(\rho,R))=\frac{\log (R/\rho)}{2 \pi}.
$$
By taking $R$ relatively larger than $\rho$,
this modulus is arbitrarily large.
Let us choose a uniformization 
$\psi:\widehat{W} \to \D$
such that $\psi(s)=0$.
For an arbitrarily small $\nu>0$,
we may assume that $\psi(W) \subset D(\nu)$
when $\mathrm{mod}(\widehat{W} \sminus \overline{W})$ 
is sufficiently large.
(See \cite[Theorems 2.1 and 2.4]{McMullen 1994}.
Indeed, it is enough to take $R$ such that $\nu \asymp \rho/R$.)
By the Bers-Royden theorem 
(\cite[Theorem 1]{Bers-Royden 1986}), 
each $\iota_c:J(f_s) \to J(f_c)$ extends to a 
$(1 + \nu)/(1 - \nu)$-quasiconformal map on $\C$.
Thus the dilatation is uniformly smaller than 
$(1+\kappa)^{1/2}$ if we choose a sufficiently large $R$.
\QED

\medskip 

\paragraph{\bf Definition of the families $\widehat{\bs G}$ 
and $\bs G$.} 
For each $c \in \widehat{W}$, 
let $\widehat{V}_c'$ be the connected component of 
$f_{c}^{-L}(\widehat{V}_{c})$ 
(or, that of $P_{c}^{-pL-N}(\widehat{U}_{c})$)
containing the critical point $0$. 
We define $G_c:\widehat{V}_c' \to  \widehat{U}_{c}$ 
by the restriction of $P_c^N \cc f_c^L = P_{c}^{pL+N}$ on $\widehat{V}_c'$. 
Then we have a family of quadratic-like maps 
$$
\widehat{\bs G}
:= \{G_c: \widehat{V}'_c \to \widehat{U}_c\}_{c  \,\in \, \widehat{W}}.
$$

Similarly, for each $c \in W$, 
let $V_c'$ be the connected component of 
$f_{c}^{-L}({V}_{c})$ 
(or, that of $P_{c}^{-pL-N}({U}_{c})$)
containing $0$.
Then we have a quadratic-like family 
$$
\bs{G}:= 
\{G_c: {V}_c' \to {U}_c\}_{c  \,\in \, W}.
$$

Note that
both the annuli $\widehat{U}_c \sminus \overline{\widehat{V}_c'}$ 
and $U_c \sminus \overline{V_c'}$ 
contain the Cantor Julia set $J(f_c)$ for each $c \in \widehat{W}$.

\paragraph{\bf Claim 5 (Extending the holomorphic motion).}~
{\it
If $R$ is sufficiently large and relatively larger than $\rho$, 
we have the following extensions of 
the holomorphic motion $\iota$ 
of the Julia set $J(f_s)$ given in Claim 4:
\begin{enumerate}[\rm (1)]
\item
An extension to the holomorphic motion of 
$J(f_s) \cup \partial \widehat{V}_s' \cup \partial U_s
\cup \partial \widehat{U}_s$ over $\widehat{W}$
that is equivariant to the action of 
$G_c:\partial \widehat{V}_c' \to \partial\widehat{U}_c$.
\item
A further extension of {\rm (1)}
to the motion of the closed annulus $\Bar{\widehat{U}_s} \sminus \widehat{V}_s'$ over $\widehat{W}$.  
\item
An extension to the holomorphic motion of 
$J(f_s) \cup \partial {V}_s' \cup \partial U_s$ 
over ${W}$
that is equivariant to the action of 
$G_c:\partial {V}_c' \to \partial{U}_c$.
\item
A further extension of {\rm (3)}
to the motion of
the closed annulus $\Bar{{U}_s} \sminus {V}_s'$ over ${W}$. 
\end{enumerate}
In particular, the quasiconformal map
$\iota_c: \Bar{{U}_s} \sminus {V}_s' \to \Bar{{U}_c} \sminus {V}_c'$
induced by {\rm (4)}
extends to a $(1+\kappa)^{1/2}$-quasiconformal map on the plane 
for each $c \in W$.
}

\begin{figure}[htbp]
\begin{center}
\includegraphics[width=.38\textwidth]{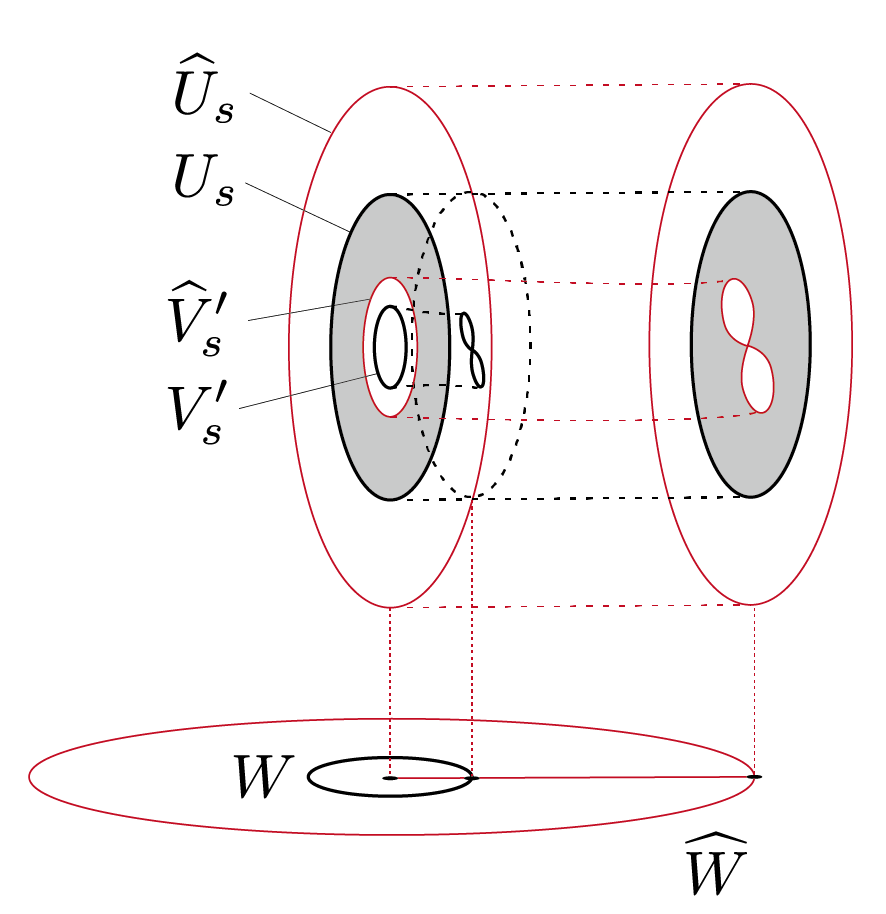}
\end{center}
\caption{\small
Extending the holomorphic motion to the closed annuli.
The motion of $J(f_s)$ is contained in the motion of shadowed annuli.
(Note that the nested annuli are extremely thick in our setting.)
}
\label{fig_extending_holo_motion}
\end{figure}

See Figure \ref{fig_extending_holo_motion}.

\medskip

\paragraph{\bf Proof.} 
(1)~
The sets $\partial U_c$, 
$\partial \widehat{V}_c'$, and 
$\partial \widehat{U}_c$ are all images of round circles
by equivariant analytic families of locally conformal injections 
over $\widehat{W}$.
In particular, they never intersect with the Julia set $J(f_c)$
for each $c \in \widehat{W}$.
Hence the extension of $\iota$ to 
$J(f_s) \cup \partial \widehat{V}_s' \cup \partial U_s
\cup \partial \widehat{U}_s$
 over $\widehat{W}$ is straightforward.\\
(2)~
By S\l odkowski's theorem (\cite{Slodkowski 1991}), 
(1) extends to the motion of $\C$,
and its restriction to 
the closed annulus $\Bar{\widehat{U}_s} \sminus \widehat{V}_s'$
is our desired motion.
Note that 
$\iota_c:\Bar{\widehat{U}_s} \sminus \widehat{V}_s'
\to \Bar{\widehat{U}_c} \sminus \widehat{V}_c'$ 
is uniformly $(1+\kappa)^{1/2}$-quasiconformal 
for $c \in W \Subset \widehat{W}$ by taking a sufficiently large $R$
that is relatively larger than $\rho$. (See the proof of Claim 4.) \\
(3)~
Similarly, $\partial V_c'$ is an image of a round circle 
$\partial U_c$ by an analytic family of injections
(where each injection is locally a univalent branch of $G_c^{-1}$) 
for $c \in W$.
Since $V_c' \Subset \widehat{V}_c'$, 
$V_c'$ never intersects with $J(f_c)$
for $c \in W$ and 
we obtain an extension of the motion of $J(f_s)$
to that of $J(f_s) \cup \partial V_s' \cup \partial U_s$ over $W$ 
which satisfies $G_c \cc \iota_c=\iota_c \cc G_s$
on $\partial V_s'$. \\
(4)~
To extend (3) to the closed annulus 
$\Bar{{U}_s} \sminus {V}_s'$, 
we divide the annulus into two annuli 
$\Bar{{U}_s} \sminus \widehat{V}_s'$ and 
$\Bar{\widehat{V}_s'} \sminus {V}_s'$.
The desired motion of $\Bar{{U}_s} \sminus \widehat{V}_s'$ over $W$
is contained in the motion given in (2).
For the annulus $\Bar{\widehat{V}_s'} \sminus {V}_s'$,
we note that the map 
$G_c: \Bar{\widehat{V}_c'} \sminus {V}_c' \to \Bar{\widehat{U}_c} \sminus {U}_c$
is a holomorphic covering of degree two.
Hence we can pull-back the motion of $\Bar{\widehat{U}_s} \sminus {U}_s$ 
over $W$ that is contained in the motion given in (2)
by these covering maps.
More precisely, we can construct  
an analytic family 
$\iota_c:\Bar{\widehat{V}_s'} \sminus \widehat{V}_s' \to 
\Bar{\widehat{V}_c'} \sminus \widehat{V}_c'$ of 
$(1+\kappa)^{1/2}$-quasiconformal maps 
that agrees with the motion of $J(f_s) \cup \partial V_s' \cup \partial U_s$,
by taking a branch of 
$G_c^{-1} \cc \iota_c\cc G_s$
with $\iota_c$ given in (2).
\QED

\medskip

\paragraph{\bf Claim 6 (Decorated tubing).}
{\it 
By taking larger $R$, $r$, and $\rho$ if necessary, 
there exist a $\rho'>0$ and a tubing
$$
\check{\Theta}
:=
\braces{
\check{\Theta}_c: \Bar{A(\check{R}, \check{R}^2)} \to \Bar{U_c} \sminus {V_c'}
}_{c\, \in \,W}
$$
of the family $\bs{G}$ with the following properties:
\begin{enumerate}[\rm (1)]
\item
$\check{R}=\rho/\rho'$ and 
$\Gamma_0(c_0+\eta)=\Gamma_0(c_0+\eta)_{\rho',\rho}$
is contained in $A(\check{R}, \check{R}^2)$. 
\\[-.9em]
\item
$\check{\Theta}_s(\check{z})
=h_s^{-1}\paren{((\rho')^{\,2}/\rho) \cdot \check{z}}$ 
for $\check{z} \in \Gamma_0(c_0+\eta)$
such that $\check{\Theta}_s$ maps $\Gamma_0(c_0+\eta)$ onto $J(f_s)$.
\\[-.7em]
\item 
Each $\check{\Theta}_c:\Bar{A(\check{R}, \check{R}^2)} \to \Bar{U_c} \sminus {V_c'}$
is a $(1+\kappa)$-quasiconformal embedding
that is compatible with the holomorphic motion of 
$\Bar{U_s}\sminus V_s'$ over $W$ given in (2) of Claim 5. 
More precisely, 
we have $\check{\Theta}_c=\iota_c \cc \check{\Theta}_s$
for each $c \in W$, 
where $\iota_c: \Bar{U_s}\sminus V_s' \to \Bar{U_c}\sminus V_c'$ is 
the quasiconformal map 
induced by the motion.
\end{enumerate}
}

We call this tubing $\check{\Theta}$ 
a {\it decorated tubing} of $\bs{G}$.

\paragraph{\bf Proof.}
For each $c \in \widehat{W}$,
the map ${G}_c=P_c^{pL+N}|_{\widehat{V}_c'}$ 
can be decomposed as
${G}_c=Q_c \cc P_0$, where $P_0(z)=z^2$ and
$Q_c:P_0(\widehat{V}_c') \to \widehat{U}_c$ is an isomorphism. 
Let 
$$
\beta_c:=Q_c'(0)
\quad \text{and} \quad 
\gamma_c:=Q_c(0).
$$
Note that $\gamma_c={G}_c(0) \in U_c$ if $c \in W$.

Since $\widehat{U}_c$ and $U_c$ are round disks 
of radii $R|\al_c|$ and $\rho|\al_c|$ respectively, 
we apply the Koebe distortion theorem to 
$Q_c^{-1}:\widehat{U}_c \to P_0(\widehat{V}_c')$
and obtain 
\begin{equation}\label{eq_Q_c_inverse}
z=Q_c^{-1}(w)=\beta_c^{-1}(w-\gamma_c)\paren{1+O(\rho/R)}
\end{equation}
for $w \in U_c$.
Indeed, by the Koebe distortion theorem (\cite[\S 2.3]{D 1983}),
we have 
$$
\abs{\frac{(Q_c^{-1})'(w)}{(Q_c^{-1})'(\gamma_c)}}=1+O(\rho/R)
\quad
\text{and}
\quad
\arg\frac{(Q_c^{-1})'(w)}{(Q_c^{-1})'(\gamma_c)}=O(\rho/R)
$$
for $w \in U_c$.
Hence we have $(Q_c^{-1})'(w)=\beta_c^{-1}(1+O(\rho/R))$ on $U_c$.
By integrating the function $(Q_c^{-1})'(w)-\beta_c^{-1}$ 
along the segment joining $\gamma_c$ to $w$ in $U_c$, 
we obtain
$$
|Q_c^{-1}(w)-\beta_c^{-1}(w-\gamma_c)|
=|w-\gamma_c||\beta_c|^{-1}O(\rho/R)
$$
that is equivalent to \eqref{eq_Q_c_inverse}.

This implies that
$$
G_c(z)=Q_c (z^2)=\gamma_c + \beta_c z^2\, \paren{1+O(\rho/R)}
$$ 
on $V_c'$.
By an affine coordinate change 
$$
\check{z}=\check{A}_c(z):=\beta_c\, z,
$$
we obtain a quadratic-like map $\check{G}_c:\check{V}_c' \to \check{U}_c$
of the form 
\begin{equation}
\check{w}=\check{G}_c(\check{z}):=\check{A}_c \cc G_c \cc \check{A}_c^{-1}(\check{z})
=\beta_c\gamma_c+\check{z}^2 
\paren{1+O(\rho/R)},
\label{eq_G_check}
\end{equation}
where $\check{V}_c':=\check{A}_c(V_c')$ and 
$\check{U}_c:=\check{A}_c(U_c)=D(\rho|\al_c||\beta_c|)$.

Now suppose that $c=s$.
Then the condition $G_s(0)=0$ implies $\gamma_s=0$.
Hence we have 
\begin{equation}\label{eq_G_check_s}
\check{w}=\check{G}_s(\check{z})=\check{z}^2(1+O(\rho/R))
\quad \text{and} \quad
\check{z}=\check{G}_s^{-1}(\check{w})=\sqrt{\check{w}(1+O(\rho/R))}.
\end{equation}
Let $\check{R}:=(\rho|\al_s||\beta_s|)^{1/2}$.
Then $\check{U}_s=D(\check{R}^2)$
and $\partial \check{V}_s'$ is
$C^1$-close to a circle $\partial D(\check{R})$.
In other words, the annulus 
$\mathcal{A}:=\check{U}_s \sminus \overline{\check{V}_s'}$
is close to a round annulus $A(\check{R}, \check{R}^2)$.
Moreover, the annulus $\mathcal{A}$ contains the compact set
$\mathcal{J}:=\check{A}_s(J(f_s))=J(f_s) \times \beta_s$.

Let us define $\rho'>0$ such that $\rho/\rho'=\check{R}$, i.e.,
$$
\rho':=\frac{\rho}{\check{R}}=\paren{\frac{\rho}{|\al_s||\beta_s|}}^{1/2}.
$$
By taking a sufficiently large $N$ in Claim 1,
we may assume that the diameter of $V_s'$ is sufficiently small
(equivalently, $|\beta_s|$ is sufficiently large,
and thus $\rho'$ is sufficiently small) such that
$$
J(P_{c_0+\eta}) \subset A(\rho',\rho). 
$$
Hence the rescaled Julia set  
$$
\mathcal{J}_0:=\Gamma_0(c_0+\eta)_{\rho',\rho}
= J(P_{c_0+\eta}) \times \frac{\rho}{(\rho')^2}
$$
is contained in the annulus $A(\check{R}, \check{R}^2)$. 

\begin{lem}\label{lem_pre_tubing}
There exists a $(1+\kappa)^{1/2}$-quasiconformal map 
$\Psi:\Bar{A(\check{R}, \check{R}^2)} 
\to \Bar{\mathcal{A}}=\Bar{\check{U}_s} \sminus {\check{V}_s'}$
 such that $\Psi(\mathcal{J}_0)=\mathcal{J}$ and 
 $\Psi(\check{z}^2)=\check{G}_s(\Psi(\check{z}))$ 
 for any $\check{z} \in \partial D(\check{R})$
 by taking sufficiently large $R,\,\rho,\,\check{R}$ with sufficiently small $\rho/R$.
\end{lem}

\paragraph{\bf Proof of Lemma \ref{lem_pre_tubing}.}
We will construct such a $\Psi$ for $\partial A(\check{R}, \check{R}^2)$
and for $\mathcal{J}_0$ separately, 
then use the Bers-Royden theorem 
to extend it to $\overline{A(\check{R}, \check{R}^2)}$. 

Let us start with the boundary of the annulus:
By (\ref{eq_G_check_s}), we have 
$$
\log \check{G}_s^{-1}(\check{w})
= \frac{1}{2}\log \check{w}+ \frac{1}{2}\log\paren{1+O(\rho/R)}
= \frac{1}{2}\log \check{w}+ O(\rho/R)
$$
near $\partial \check{U}_s=\partial D(\check{R})$.
Hence for sufficiently large $\check{R}$ and small $\rho/R$,
we may apply the same argument as the proof of 
Lemma \ref{lem_almost_conformal}
by regarding $R$ in Lemma \ref{lem_almost_conformal} as $\check{R}^2$.
Indeed, we let 
$$
\check{\Upsilon}(\check{w}):=  
\log \check{G}_s^{-1}(\check{w})-\frac{1}{2}\log \check{w}
=O(\rho/R)
$$
and $\check{v}(t):=\check{\Upsilon}(\check{w}(2t))$,
where $\check{w}(2t)=\check{R}^2e^{2t i}$
makes two turns along $\partial D(\check{R}^2)$ as $t$ varies from $0$ to $2\pi$.
Set $\check{\ell}:=\log \check{R}$ and $\check{\eta}(s):=\eta_0(s/\check{\ell}-1)$,
where $\eta_0$ is defined in the proof of Lemma \ref{lem_almost_conformal}.
Then the map
$$
\check{\theta}_\xi(s+it):=s+it + \xi \check{\eta}(s)\check{v}(t)
$$
with a parameter $\xi \in \C$ 
is defined for $(s,t) \in [\check{\ell}, 2 \check{\ell}] \times [0,2\pi]$, 
and $\check{\theta}_\xi$ is a $(1+O(|\xi|\rho/R))$-quasiconformal map
for sufficiently large $\check{R}$ and small $\rho/R$.
We fix such $\rho, R$ and $\check{R}$,
and obtain a holomorphic family of injections
$\braces{\check{\theta}_\xi}_{\xi}$
with parameter $\xi$ in a disk $D(d_0)$  of radius $d_0 \asymp R/\rho$.
(This bound comes from the estimate like (\ref{eq_injectivity})
that ensures injectivity.)
By observing the motion through the exponential map,
we obtain a holomorphic motion 
$\psi: \partial A(\check{R},\check{R}^2) \times D(d_0) \to \C$
of $\partial A(\check{R},\check{R}^2)$ over $D(d_0)$
with $\psi(\check{z},\xi):=\psi_\xi(\check{z})$.
In particular, by letting $\xi:=1$, 
the map $\Psi:= \psi_1$  
satisfies 
 $\Psi(\check{z}^2)=\check{G}_s(\Psi(\check{z}))$
 by construction.

Next we consider the Julia set:
Let $\Psi:\mathcal{J}_0 \to \mathcal{J}$ be the quasiconformal map
given by composing the four maps
$$
\mathcal{J}_0 
=\Gamma(c_0+\eta)_{\rho',\,\rho}
\stackrel{(1)}{\longrightarrow} 
J(P_{c_0+\eta})
\stackrel{(2)}{\longrightarrow} 
J(F_{X_n(s)})
\stackrel{(3)}{\longrightarrow} 
J(f_s)
\stackrel{(4)}{\longrightarrow} 
\mathcal{J},
$$
where  (1) is the affine map $\check{z} \mapsto (\rho/(\rho')^2)^{-1} \check{z}$;
(2) is the inverse of the $(1+O(R^{-1}))$-quasiconformal 
straightening $H:=H_{X_n(s)}$ of $F_{X_n(s)}$ to $P_{c_0+\eta}$;  
(3) is the inverse of the affine map 
${A}_s: z \mapsto z/\al_s$;
and 
(4) 
is the affine map $\check{A}_s:\check{z} \mapsto \beta_s \check{z}$.
The straightening map $H$ in (2) 
extends to a $(1+O(R^{-1}))$-quasiconformal map 
on the plane as in Lemma \ref{lem_almost_conformal}.
Let $\check{\mu}:=\mu_{H^{-1}}$ be the Beltrami coefficient of the {\it inverse} of 
such an extended $H$ with $\norm{\check{\mu}}_\infty=O(R^{-1}).$
Then the Beltrami equation for $\check{\mu}_\xi:=\xi \cdot \check{\mu}$ 
with a complex parameter $\xi$ 
has a solution if $\xi \in D(d_1)$ with $d_1 \asymp R$.
Let $\phi_\xi$ be the unique normalized solution 
such that $\phi_\xi(0)=0$ and $(\phi_\xi)_z-1 \in L^{p'}(\C)$ 
for some $p'>2$.
Then the map 
$\psi_\xi(\check{z})
:=\check{A}_s \cc A_s^{-1} \cc 
\phi_\xi ((\rho/(\rho')^2)^{-1} \check{z})$ 
gives a holomorphic motion
$\psi: \mathcal{J}_0 \times D(d_1) \to \C$
of $\mathcal{J}_0$ over $D(d_1)$ with $\psi(\check{z}, \xi)=\psi_\xi(\check{z})$.
In particular, by letting $\xi:=1$, the map 
\begin{equation}\label{eq_Psi}
\Psi(\check{z}):=\psi_1(\check{z})
=\check{A}_s \cc A_s^{-1} \cc 
H_{X_n(s)}^{-1} ((\rho/(\rho')^2)^{-1} \check{z})
=\check{A}_s \cc h_s^{-1} ((\rho/(\rho')^2)^{-1} \check{z})
\end{equation}
satisfies $\Psi(\mathcal{J}_0) = \mathcal{J}$.

Now by taking $R$ relatively larger than $\rho$,
 we may assume that $d_0 <d_1$.
Let us check that the unified map 
$\psi_\xi: \partial A(\check{R},\check{R}^2) \cup \mathcal{J}_0 \to \C$
gives a holomorphic family of injections for $\xi \in D(d_0)$
if $\check{R}$ is sufficiently large.
Indeed, it is enough to check that the distance between 
$\psi_\xi(\partial A(\check{R},\check{R}^2))$
and $\psi_\xi(\mathcal{J}_0)$ is bounded from below
for $\xi \in D(d_0)$.

Let us fix a constant $0<\sigma<1$ such that 
$J(P_{c_0+\eta}) \subset A(\sigma \rho,\rho)$.
Hence $\dist(0, \mathcal{J}_0) > \sigma \check{R}^2$.
Note that we can replace $\check{R}$ by an arbitrarily larger one 
with only a slight change of $\sigma$, 
because in Claim 1 we can replace $\widehat{V}_{c_1}$ 
by an arbitrarily smaller one 
such that the location of the center $s$ of $\widehat{W}$
changes only a little (relatively to the size of $\widehat{\Lambda}$). 
Hence we may assume that $\check{R}$ is large enough such 
that
$\dist(\partial D(\check{R}), \mathcal{J}_0)
\ge \dist(0, \mathcal{J}_0) -\check{R}
> \check{R}^2(\sigma -1/\check{R})
 \asymp \check{R}^2$.
By taking a sufficiently large $\rho$, 
we have $J(P_{c_0+\eta}) \subset D(\rho/2)$ 
and thus $\dist(\partial D(\check{R}^2), \mathcal{J}_0) \ge \check{R}^2/2$.
Hence we conclude that 
$\dist(\partial A(\check{R},\check{R}^2), \mathcal{J}_0) \asymp \check{R}^2$.

Now suppose that $\xi \in D(d_0)$.
Since $d_0 \asymp R/\rho$,
an explicit calculation shows that 
$\psi_\xi(\check{z})=\check{z}$ on $\partial D(\check{R}^2)$
and 
$\psi_\xi(\check{z})=\check{z}(1+\xi O(\rho/R))$
on $\partial D(\check{R})$.
Hence $\dist(0, \psi_\xi(\partial D(\check{R}))) \asymp \check{R}$.

On the other hand, 
since $\check{\mu}_\xi=O(|\xi| R^{-1}) =O(\rho^{-1})$ for $\xi \in D(d_0)$,
we have $\phi_\xi(z)=z+O(\rho^{-1+2/p'})$
 for some $p' >2$ on $J(P_{c_0+\eta})$
 (cf. Lemma \ref{lem_estimate_of_h}).
Hence 
$\dist(\mathcal{J}_0, \psi_\xi(\mathcal{J}_0)) 
\asymp O(\rho^{-1+2/p'}) \check{R}^2$ 
for sufficiently large $\rho$.
It follows that if $\rho$, $R$, and $\check{R}$ are sufficiently large
and $\rho/R$ are sufficiently small,
then we have 
\begin{align*}
&\dist(\psi_\xi(\partial A(\check{R},\check{R}^2)),
\psi_\xi(\mathcal{J}_0)) \\
\ge~ & 
\dist(\partial A(\check{R},\check{R}^2),
\mathcal{J}_0)
-\dist(\partial A(\check{R},\check{R}^2),
\psi_\xi(\partial A(\check{R},\check{R}^2)))
 -\dist(\mathcal{J}_0,  \psi_\xi(\mathcal{J}_0))\\
 \asymp~ 
 & \check{R}^2(1 - O(1/\check{R}) -O(\rho^{-1+2/p'})) 
\asymp \check{R}^2 
\end{align*}
for $\xi \in D(d_0)$.

By applying the Bers-Royden theorem
to the holomorphic motion $\psi$ of 
$\partial A(\check{R},\check{R}^2) \cup \mathcal{J}_0$
over $D(d_0)$,
the injection $z \mapsto \psi_1(z)=\psi(z,1)$
extends to a quasiconfomal map on $\C$
whose dilatation is bounded by $1+O(1/d_0)=1+O(\rho/R)$.
It is $(1+\kappa)^{1/2}$-quasiconformal 
 by taking $R$ relatively larger than $\rho$.
Thus the restriction $\Psi$ of $\psi_1$
on $\Bar{A(\check{R}, \check{R}^2)}$
is our desired map. 
\QED ( Lemma \ref{lem_pre_tubing})

\medskip

\paragraph{\bf Proof of Claim 6, continued.} 
Let $\check{\Theta}_s:=\check{A}_s^{-1}\cc \Psi$
(hence 
$\check{\Theta}_s(\check{z})=h_s^{-1} ((\rho/(\rho')^2)^{-1} \check{z})$
for $\check{z} \in \Gamma(c_0+\eta)_{\rho',\,\rho}$ by
(\ref{eq_Psi}))
and $\check{\Theta}_c:=\iota_c \cc \check{\Theta}_s$ for $c \in {W}$,
where $\iota_c: \Bar{U_s}  \sminus  {V_s'} \to \Bar{U_c}  \sminus  {V_c'}$
is a $(1+\kappa)^{1/2}$-quasiconformal map given in Claim 5.
Then $\check{\Theta}_c$ is $(1+\kappa)$-quasiconformal for each $c \in W$
with desired properties.
\QED

\medskip

\paragraph{\bf Almost conformal embedding of the model.}
We finish the proof of Theorem C by the next claim:

\paragraph{\bf Claim 7 (Almost conformal straightening).}
{\it 
The family 
$$
\bs{G}= 
\{G_c: {V}_c' \to {U}_c\}_{c  \,\in \,W} 
$$
is a Mandelbrot-like family whose straightening map
$\chi_{\bs{G}}: W \to \C$ associated with the decorated tubing 
$\check{\Theta}$ is $(1+\kappa)$-quasiconformal. 
Moreover, the inverse of $\chi_{\bs{G}}$ realizes a $(1+\kappa)$-quasiconformal embedding 
of the model $\cM(c_0+\eta)_{\rho',\rho}$.
}

\medskip

\paragraph{\bf Proof.}
By Claim 6, the family $\bs{G}$ is equipped with the decorated 
tubing $\check{\Theta}=\{\check{\Theta}_c\}_{c \, \in \, W}$ 
and hence Mandelbrot-like with connectedness locus $M_{\bs{G}}$
homeomorphic to $M$.
The straightening map $\chi_{\bs{G}}: W \to \C$ associated with 
this tubing is given by $\chi_{\bs{G}}(c):=\check{h}_c(G_c(0))$,
where $\check{h}_c:U_c \to \C$ is the $(1+\kappa)$-quasiconformal 
straightening map of $G_c$ constructed from 
the $(1+\kappa)$-quasiconformal map $\check{\Theta}_c$ 
by the same way as the proof of Lemma \ref{lem_almost_conformal}.

Let us show that the map $\chi_{\bs{G}}$ is $(1+\kappa)$-quasiconformal
by following the argument of Lemma \ref{lem_almost_conformal_M}:
It is $(1+\kappa)$-quasiconformal 
on $W \sminus  M_{\bs{G}}$ since
$\check{\Theta}_c$ is $(1+\kappa)$-quasiconformal for any $c \in W$.
By Lemma \ref{lem_Bers} (Bers' Gluing Lemma),
it is enough to show that
$\chi_{\bs{G}}: M_{\bs{G}} \to M$
extends to a $(1+\kappa)$-quasiconformal map on $\C$. 

Recall that 
the quadratic-like map $\check{G}_c:\check{V}_c' \to \check{U}_c$
given in (\ref{eq_G_check}) is of the form 
$$
\check{w}=\check{G}_c(\check{z})=
\check{A}_c \cc G_c \cc \check{A}_c^{-1}(\check{z})
=\check{z}^2 +\check{c}+\check{u}(\check{z},\check{c}),
$$
where $\check{c}=\check{X}(c):=\beta_c\gamma_c$ and 
$\check{u}(\check{z},\check{c})=\check{z}^2 \,O(\rho/R)$.
Moreover, $\check{u}(\check{z},\check{c}) =\check{u}_{\check{c}}(\check{z})$
 satisfies $\check{u}_{\check{c}}'(0)=0$,
 and the value 
$$
\check{\delta}:=
\sup \{|\check{u}(\check{z},\check{c})| \st (\check{z},\check{c}) 
\in \overline{D(4)}\times \overline{D(4)}\}
$$
is $O(\rho/R)$.
We may assume that $\rho/R$ is small enough such that $\check{\delta}<1$.
As in the proof of Lemma \ref{lem_almost_conformal_M},
we consider the analytic family
$$
\check{\bs G}(t)=
\skakko{\check{z} \mapsto \check{z}^2+\check{c}
+t \check{u}(\check{z},\check{c})/\check{\delta}}_{\check{c}\, \in \, D(4)}
$$
with parameter $t \in \D$ 
whose connectedness locus $\check{M}(t)$ is homeomorphic to $M$.
Then $\check{\bs G}(\check{\delta})=\{\check{G}_c\}_{c \,\in \, W}$
and the straightening map 
$\check{\chi}:\check{M}(\check{\delta}) \to M$ 
extends to a quasiconformal map on $\C$
with dilatation $1+O(\check{\delta})=1+O(\rho/R)$
by the Bers-Royden theorem.
Hence if $\rho/R$ is sufficiently small, 
$\check{\chi}$ is $(1+\kappa)$-quasiconformal.
Since $\chi_{\bs G}(c) =\check{\chi}(\check{X}(c))$ and
$\check{X}(c)=\beta_c \,\gamma_c$ is holomorphic near $M_{\bf G}$,
we conclude that $\chi_{\bs{G}}: M_{\bs{G}} \to M$
extends to a $(1+\kappa)$-quasiconformal map on $\C$.

As in the proof of Theorem A, the inverse of $\chi_{\bs G}(c)$ 
realizes a $(1+\kappa)$-quasiconformal embedding of 
$\cM(c_0+\eta)_{\rho',\rho}$ into $W$.
\QED

\medskip 

\section{Proof of Corollary D}
\paragraph{\bf Proof.}
We recall the setting of Theorem A. Take any small Mandelbrot set
$M_{s_0}$, where $s_0 \ne 0$ is a superattracting parameter and
take any Misiurewicz or parabolic parameter $c_0 \in \partial M$.
Then Theorem A shows that ${\mathcal M}(c_0+\eta)$ appears quasiconformally
in $M$ in a small neighborhood of $c_1 := s_0 \perp c_0$. Now let
$c$ be a parameter which belongs to the quasiconformal image of the 
decoration of ${\mathcal M}(c_0+\eta)$. This means that
$$
  G_c^k(0) \in Y:=\Theta_c(\Gamma_0(c_0+\eta)) = J(f_c) \quad 
  \text{for some} \ k \in \N.
$$
Since $G_c = P_c^{p'}$ for some $p' \in \N$, we have
$P_c^{p'k}(0) \in Y$. On the other hand, 
$Y$ is $f_c$ ($=P_c^p$)-invariant, that is, 
$P_c^{pn}(Y) \subseteq Y$ for every $n$. Then
for a fixed $0 \leq  r < p$, we have
$P_c^{pn+r}(Y) \subseteq P_c^r(Y)$ for every $n$ and each $P_c^r(Y)$ is 
apart from $0$. Therefore for every $i \geq p'k_0$ we have
$P_c^i(Y) \in \bigcup_{r=0}^{p-1} P_c^r(Y)$, which implies that
the orbit of $0$ under the iterate of $P_c$ does not accumulate 
on $0$ itself. Moreover, $P_c$ has no parabolic periodic point
since $Y = J(f_c)$ is a Cantor Julia set of a hyperbolic quadratic-like
map $f_c : U_c' \to U_c$. 
This shows that $P_c$ is semihyperbolic. 
Since Misiurewicz or parabolic parameters are dense in $\partial M$, we can
find decorations in every small neighborhood of every point in $\partial M$. 
Also there are only countably many Misiurewicz parameters. Hence it follows 
that the semihyperbolic parameters which are not Misiurewicz and 
non-hyperbolic are dense in $\partial M$.
\QED

\section{Concluding Remarks}
We have shown that we can see quasiconformal images of some Julia sets
in the Mandelbrot set $M$. But this is not satisfactory, because these 
images are all Cantor sets and disconnected. On the other hand, $M$ is 
connected and so what we have detected is only a small part of the whole
structure of $M$. From computer pictures, it is observed that the points
in these Cantor sets are connected by some complicated filament structures.
This looks like a picture which is obtained from the picture of $K(P_c)$ 
for $c \in {\rm int} (M_{s_1})$ by replacing all small filled Julia sets
with small Mandelbrot sets.
It would be interesting to explain this mathematically.
A similar phenomena as in the quadratic family are observed also in the 
unicritical family $\{ z^d+c \}_{c \in \C}$ by computer pictures. 
These phenomena should
be proved in the same manner as for the quadratic case.


\begin{thebibliography}{99}

\parskip=3pt

\bibitem[{\bf A}]{Ahlfors 1978}
L.V.~Ahlfors, 
{\it Complex Analysis, third ed.},
McGraw-Hill Book Co., New York, 1978.

\bibitem[\bf B\rm]{Beardon 1991}
A.F.~Beardon, \it
Iteration of Rational Functions\rm, 
Springer Verlag, New York, Berlin and Heidelberg, 1991.


\bibitem[\bf BH\rm]{Buff-Henriksen 2001}
X.~Buff and C.~Henriksen, Julia sets in parameter spaces, \it 
Comm. Math. Phys. \rm {\bf 220} (2001), no. 2, 333--375.


\bibitem[\bf BR\rm]{Bers-Royden 1986} 
L.~Bers and H.L.~Royden,
Holomorphic families of injections, 
{\it Acta Math.} {\bf 157} (1986), 259--286.




\bibitem[\bf CJY\rm]{Carleson-Jones-Yoccoz 1994}
L.~Carleson, P.W.~Jones and J.-C.~Yoccoz, 
Julia and John, \it
Bol. Soc. Brasil. Mat. (N.S.) {\bf 25} \rm (1994), no. 1, 1--30. 




\bibitem[\bf CRY\rm]{Cornell-Rojas-Yampolsky 2017}
D.~Cornell, C.~Rojas and M.~Yampolsky, 
Non computable Mandelbrot-like set for a one-parameter complex family, \it
Inform. and Comput. {\bf 262} \rm (2018), part 1, 110--122. 


\bibitem[\bf CT\rm]{Cui-Tan 2018}
G.~Cui and L.~Tan, 
Hyperbolic-parabolic deformations of rational maps, \it
Science China Math. {\bf 61} \rm (2018) 2157--2220.


\bibitem[{\bf D}]{D 1983} 
P.L. Duren,
{\it Univalent Functions}, 
Springer-Verlag, 1983.


\bibitem[\bf D-BDS\rm]{Douady 2000}
A. Douady, X. Buff, R. Devaney and P. Sentenac, 
Baby Mandelbrot sets are born in cauliflowers, 
in\it \ The Mandelbrot set, theme and variations\rm, London Math. Soc. Lecture 
Note Ser. {\bf 274}, Cambridge Univ. Press, Cambridge (2000), 19--36.


\bibitem[{\bf DH1}]{DH Orsay}
A. Douady and J. H. Hubbard, 
Etude dynamique des polyn\^omes complexes I \& II, 
{\it Publ. Math. Orsay, Universit\'e de Paris-Sud}, 
D\'epartement de Math\'ematiques, Orsay, 
1984/85, 84-2, 85-4.


\bibitem[\bf DH2\rm]{Douady-Hubbard 1985}
A. Douady and J. Hubbard, 
On the dynamics of polynomial-like mappings, \it
Ann. Sci. \'Ecole Norm. Sup. \rm {\bf 18} (1985), 287--343.

\bibitem[{\bf EE}]{EE 1985} 
J.-P.~Eckmann and H.~Epstein, 
Scaling of Mandelbrot sets generated by critical point preperiodicity, 
{\it Comm. Math. Phys.}, 
{\bf 101} (1985), 283--289.



\bibitem[{\bf G}]{Gaidashev 2007}
D.~Gaidashev, 
Cylinder renormalization for Siegel discs and 
a constructive measurable Riemann mapping theorem, 
{\it Nonlinearity}, {\bf 20} (2007), 713--741.


\bibitem[\bf H\rm]{Haissinsky 2000}
P. Ha\"{\i}ssinsky, 
Modulation dans l'ensemble de Mandelbrot, 
in\it \ The Mandelbrot set, theme and variations\rm, London Math. Soc. Lecture 
Note Ser. {\bf 274}, Cambridge Univ. Press, Cambridge (2000), 37--66.


\bibitem[{\bf IT}]{IT Book}
Y. Imayoshi and M. Taniguchi, 
{\it An Introduction to Teichm\" uller Spaces}. 
Springer, 1992.


\bibitem[{\bf J}]{Jung 2015}
W. Jung, personal communication, 2015. See 
\verb|http://www.mndynamics.com/indexp.html| (his web page)
and also
\verb|http://www.mrob.com/pub/muency/embeddedjuliaset.html|.

\bibitem[{\bf K}]{Kawahira 2014} 
T.~Kawahira, 
Quatre applications du lemme de Zalcman \`a la dynamique complexe, 
{\it J. Anal. Math.} {\bf 124} (2014), 309--336.

\bibitem[{\bf KK}]{Kawahira-Kisaka 2023} 
T.~Kawahira and M.~Kisaka,
Julia sets appear quasiconformally in the Mandelbrot set, II: A parabolic proof, 
{\it Preprint}. 

\bibitem[{\bf Ly1}]{Lyubich 1991} 
M. Lyubich, 
On the Lebesgue measure of the Julia set of a quadratic
polynomial, {\it Preprint} IMS at Stony Brook, \#1991/10.



\bibitem[{\bf Ly2}]{Lyubich 1999} 
M. Lyubich, 
Feigenbaum-Coullet-Tresser universality and Milnor's hairiness conjecture,
{\it Ann. of Math.} (2) {\bf 149} (1999), no. 2, 319--420.


\bibitem[{\bf Ly3}]{Lyubich Book}
M. Lyubich,
{\it 
Conformal Geometry and Dynamics of Quadratic Polynomials, vol I--II.
} Available in his web page 
(\verb|http://www.math.stonybrook.edu/~mlyubich/book.pdf|).



\bibitem[{\bf LyMin}]{Lyubich-Minsky 1997}
M.~Lyubich and Y.~Minsky, Laminations in holomorphic dynamics,
{\it J. Differential Geom.} {\bf 47} (1997), no. 1, 17--94. 






\bibitem[{\bf Ma}]{Mane 1993}
R. Ma\~n\'e, 
On a theorem of Fatou, \it
Bol. Soc. Brasil. Mat. (N.S.) \bf 24 \rm No.1 (1993), 1--11.





\bibitem[{\bf Mc1}]{McMullen 1994}
C.T. McMullen, \it
Renormalization and complex dynamics\rm, 
Ann. of Math. Studies {\bf 135}, Princeton Univ. Press, 1994. 



\bibitem[{\bf Mc2}]{McMullen 2000}
C.T. McMullen, 
The Mandelbrot set is universal, 
in\it \ The Mandelbrot set, theme and variations\rm, London Math. Soc. 
Lecture Note Ser. {\bf 274}, Cambridge Univ. Press, Cambridge (2000), 
1--17.


\bibitem[\bf Mil1\rm]{Milnor 2000}
J. Milnor, 
Periodic orbits, externals rays and the Mandelbrot set: an expository account, 
\it G\'eom\'etrie complexe et syst\`emes dynamiques \rm(Orsay, 1995), 
Ast\'erisque \bf 261 \rm (2000), xiii, 277--333.


\bibitem[\bf Mil2\rm]{Milnor 2006}
J. Milnor, \it
Dynamics in one complex variable, third ed.\rm, 
Annals of Mathematics Studies, vol. {\bf 160}, Princeton University Press,
Princeton, NJ, 2006.


\bibitem[\bf MNTU\rm]{MNTU 2000}
S. Morosawa, Y. Nishimura, M. Taniguchi and T. Ueda, \it 
Holomorphic dynamics, \rm Translated from the 1995 Japanese original 
and revised by the authors. Cambridge Studies in Advanced Mathematics 
{\bf 66}, Cambridge University Press, Cambridge, 2000.


\bibitem[\bf MTU\rm]{MTU 1995}
S. Morosawa, M. Taniguchi and T. Ueda, 
\it Holomorphic dynamics \rm (in Japanese), 
Baifuukan (publisher), 1995.


\bibitem[{\bf PS}]{Peitgen-Saupe 1988} 
H.-O.~Peitgen and D.~Saupe (eds.), 
\it The science of fractal images\rm, 
Springer-Verlag, New York, 1988.

\bibitem[\bf Riv\rm]{Rivera-Letelier 2001}
J. Rivera-Letelier, On the continuity of Hausdorff dimension of Julia sets and 
similarity between the Mandelbrot set and Julia sets, \it Fund. Math. \rm 
{\bf 170} (2001), no. 3, 287--317.


\bibitem[\bf S\rm]{Shishikura 1998}
M. Shishikura, 
The Hausdorff dimension of the boundary of the Mandelbrot set and Julia sets,
\it Ann. of Math. \rm (2) {\bf 147} (1998), no. 2, 225--267.


\bibitem[{\bf S\l}]{Slodkowski 1991}
Z. S\l odkowski, 
{\rm Holomorphic motions and polynomial hulls}, \it
Proc. Amer. Math. Soc. {\bf 111} \rm (2) (1991), 347--355.


\bibitem[\bf T\rm]{Tan Lei 1990}
L. Tan, 
Similarity between the Mandelbrot set and Julia sets, 
\it Comm. Math. Phys. 
{\bf 134} \rm (1990), no. 3, 587--617.

\end{thebibliography}
\end{document}